\documentclass[10pt]{amsart}
\usepackage{amsmath,amsfonts,amssymb}
\usepackage{amssymb,url,xspace}
\usepackage{mathtools}
\usepackage[margin=2.1cm]{geometry}
\usepackage[english]{babel}
\usepackage[T1]{fontenc}
\usepackage[sort&compress,numbers]{natbib}

\numberwithin{equation}{section}
\usepackage{amssymb,url,xspace}
\usepackage{amsmath,amsfonts,amssymb}
\usepackage{mathtools}
\usepackage{mathrsfs}
\usepackage{textcomp}
\usepackage{fancyhdr}
\usepackage{eso-pic}
\usepackage{bbold}
\usepackage{comment}
\usepackage{color}
\usepackage{xcolor}
\usepackage{stmaryrd}
\usepackage[artemisia]{textgreek}

\def\XXint#1#2#3{{\setbox0=\hbox{$#1{#2#3}{\int}$}
     \vcenter{\hbox{$#2#3$}}\kern-.5\wd0}}

\let\dem=\proof
\let\enddem=\endproof
\let\vph=\varphi
\let\eps=\varepsilon
\DeclarePairedDelimiter\abs\lvert\rvert
\DeclarePairedDelimiter\norm\lVert\rVert

\DeclareMathOperator\dvg{div}
\DeclareMathOperator\rot{curl}
\DeclareMathOperator\adj{Adj}
\newcommand\R{\mathbb{R}}
\newcommand\N{\mathbb{N}}
\newcommand\D{\mathbb{D}}

\newcommand\dpt{\partial_t}

\newcommand\cC{\mathscr{C}}

\newcommand\loc{\text{loc}}
\newtheorem{theo}{Theorem}[section]
\newtheorem{lemm}[theo]{Lemma}

\theoremstyle{definition}
\newtheorem{definition}[theo]{Definition}

\newtheorem{prop}[theo]{Proposition}

\newtheorem{rema}[theo]{Remark}
\newtheorem{coro}[theo]{Corollary}

\definecolor{mycolor}{HTML}{D35400}

\definecolor{hide}{HTML}{A0AC81}
\definecolor{myred}{HTML}{8A1538}

\usepackage[colorlinks=true,citecolor=blue,urlcolor=cyan]{hyperref}
\usepackage{cleveref}
\crefformat{prop}{#2 Proposition~#1#3}{}
\crefformat{lemma}{#2Lemma~#1#3}{}
\crefformat{rema}{#2Remark~#1#3}{}
\crefformat{theo}{#2Theorem~#1#3}{}
\crefformat{appendix}{#2Appendix~#1#3}{}
\crefformat{section}{#2Section~#1#3}{}
\crefformat{coro}{#2Corollary~#1#3}{}
\crefformat{nota}{#2Notation~#1#3}{}
\crefformat{def}{#2Definition~#1#3}{}
\begin{document}

\title{Discontinuous solutions for the Navier-Stokes equations with density-dependent viscosity}

\author{}
\address{}
\curraddr{}
\email{}
\thanks{}

\author[]{Sagbo Marcel ZODJI}
\address{Université Paris Cité and Sorbonne Université, CNRS, IMJ-PRG, F-75013 Paris, France.}
\email{marcel.zodji@imj-prg.fr}
\thanks{}


\date{\today}

\dedicatory{}


\begin{abstract}
We prove existence of a unique  global-in-time weak solutions of the Navier-Stokes equations that govern the motion of a 
compressible viscous fluid with density-dependent viscosity in two-dimensional  space. The initial velocity belongs to the Sobolev space 
$H^1(\R^2)$, and the initial  fluid density is $\alpha$-H\"older continuous on both sides of a $\cC^{1+\alpha}$-regular interface with some geometrical assumption. We prove that this configuration persists over time: the initial interface is transported by the flow to an interface that maintains the same regularity as the initial one.

Our result generalizes previous known of  Hoff \cite{hoff2002dynamics}, Hoff and Santos \cite{hoff2008lagrangean} concerning the propagation of regularity for
discontinuity surfaces by allowing more general nonlinear pressure law and density-dependent viscosity. Moreover, it supplements the work by  Danchin,  Fanelli and  Paicu \cite{danchin2020well} with global-in-time well-posedness, even for density-dependent viscosity and we achieve uniqueness in a large space. 
\end{abstract}

\maketitle


\section{Introduction}
\subsection{Presentation of the model}
In this paper, we study the problem of existence and uniqueness of global-in-time weak solutions with intermediate regularity for the Navier-Stokes equations describing the motion of compressible fluid with density-dependent viscosity in $\R^2$.  Our main interest is to generalize the result by  Hoff \cite{hoff2002dynamics}, Hoff and Santos \cite{hoff2008lagrangean} concerning the propagation of discontinuous surfaces by allowing nonlinear pressure law and density-dependent viscosity. We aim to  supplement the work by  Danchin,  Fanelli,  Paicu \cite{danchin2020well} with global-in-time well-posedness, even for density-dependent viscosity. 
Indeed, we consider the following system:
\begin{gather}\label{ep4.1}
\begin{cases}
    \dpt \rho +\dvg (\rho u)=0,\\
    \dpt (\rho u)+\dvg (\rho u\otimes u)+\nabla P(\rho)=\dvg (2\mu(\rho)\D u)+\nabla(\lambda(\rho)\dvg u)
\end{cases}
\end{gather}
describing the motion of a compressible viscous fluid at constant temperature.  Above, $\rho=\rho(t,x)\geqslant 0$ and  $u=u(t,x)\in \R^2$ are respectively the density and the velocity of the fluid and they are the  unknowns of the problem.  Meanwhile, $P=P(\rho)$, $\mu=\mu(\rho)$,  $\lambda=\lambda(\rho)$ are respectively the pressure, dynamic and 
kinematic viscosity law of the fluid and they are given $\cC^2$-regular functions of the density. The equations \eqref{ep4.1} are supplemented with initial 
data 
\begin{gather}\label{c5.1}
    \rho_{|t=0}=\rho_0\in L^\infty(\R^2)\quad \text{ and }\quad u_{|t=0}=u_0\in H^1(\R^2).
\end{gather}
We assume there exists $\widetilde\rho>0$ such that
\begin{gather}\label{c5.2}
\rho_0-\widetilde\rho\in L^2(\R^2),\quad \text{and we define}\quad \widetilde P= P(\widetilde \rho), \quad  \widetilde\mu= \mu(\widetilde\rho).
\end{gather}
 Additionally, we suppose that $\rho_0$ is upper bounded, bounded away from zero:
\begin{gather}\label{epq30}
0<\rho_{*,0}:=\inf_{x\in \R^2}\rho_0(x)\leqslant \sup_{x\in \R^2}\rho_0(x)=:\rho^*_0<\infty,\quad \mu_{0,*}:= \inf_{x\in \R^2}\mu(\rho_0(x))>0,
\end{gather}
and  $\alpha$-H\"older continuous on both sides of a  $\cC^{1+\alpha}$-regular non-self-intersecting curve $\mathcal{C}(0)$,  which is the boundary of an open, bounded and simply connected domain $D(0)\subset \R^2$. The latter regularity is defined as follows (based on \cite{zodji2023well}):

\begin{definition}\label[def]{def1}
\noindent

\begin{enumerate}
\item   \label{c3.24} 
    We say that an  interface  $\mathcal{C}$ is $\cC^{1+\alpha}$-regular and non-self-intersecting if:
    \begin{itemize}
    \item \textbf{$\cC^{1+\alpha}$-regularity}: There exist  intervals  $V_j\subset \R$, $j\in \llbracket 1,J\rrbracket$ and maps
\[
 \gamma_j\colon  V_j\mapsto \R^2\in \cC^{1+\alpha} \quad \text{such that}\quad \mathcal{C}\subset \bigcup_{j=1}^J \gamma_j (V_j),
\]
with well-defined normal vector fields.
\item \textbf{Non-self-intersection condition:} There exists $c_\gamma>0$ such that:
\begin{gather}\label{ep4.2}
    \forall\; j\in \llbracket 1, J\rrbracket,\,(s,s')\in V_j\times V_j\; \;\text{we have}\;\;\abs{\gamma_j(s)-\gamma_j(s')}\geqslant c_\gamma^{-1}\abs{s-s'}.
\end{gather}
\end{itemize}

\item Consider an open, bounded, and simply connected domain $D$ in $\R^2$. Assume that $D$ is  $\cC^{1+\alpha}$-regular. There exists a  function $\vph\colon \R^2\mapsto \R\in \cC^{1+\alpha}$ such that 
\[
D= \{x\in \R^2\colon \vph(x)>0\},\quad\text{and} \quad \abs{\nabla \vph}_{\text{inf}}:=\inf_{x\in\partial D}\abs{\nabla \vph(x)}>0. 
\]
We refer to  \cite[Section 3.1]{kiselev2016finite} for the construction of such level-set function. 
We then define:
\begin{gather}\label{v5}
\ell_{\vph}=\min\left\{1,\left(\dfrac{\abs{\nabla \vph}_{\text{inf}}}{\norm{\nabla\vph}_{\dot \cC^\alpha}}\right)^{1/\alpha}\right\}.
\end{gather}
\item  Given that
\[
\R^2= D\cup \mathcal{C}\cup (\R^2\setminus \overline{D}),
\]
we   define the space of piecewise $\alpha$-H\"older continuous  functions with respect to $\mathcal{C}$ as follows: 
\[
\dot \cC^\alpha_{pw,\gamma}(\R^2):=\dot \cC^\alpha(\overline{D})\cap \dot \cC^\alpha(\R^2\setminus D),\;\text{ with }\; \norm{g}_{\dot\cC^\alpha_{pw,\gamma}(\R^2)}= \norm{g}_{\dot\cC^\alpha(\overline{D})}+\norm{g}_{\dot\cC^\alpha(\R^2\setminus D)},
\]
and the non-homogeneous space
\[
 \cC^\alpha_{pw,\gamma}(\R^2):= L^\infty(\R^2)\cap \dot \cC^\alpha_{pw,\gamma}(\R^2), \;\text{ with }\; \norm{g}_{\cC^\alpha_{pw,\gamma}(\R^2)}=\norm{g}_{L^\infty(\R^2)}+\norm{g}_{\dot\cC^\alpha_{pw,\gamma}(\R^2)}.
\]
This spaces strictly contain the H\"older space $\cC^\alpha(\R^2)$. 
\item Given a function $g\in \cC^{\alpha}_{pw,\gamma}(\R^2)$, the jump $\llbracket g\rrbracket$ and the average $<g>$ of $g$ are defined as follows: for all $\sigma\in \mathcal{C}$,
\begin{gather}\label{v3}
\begin{cases}\displaystyle
\llbracket g\rrbracket(\sigma)=\lim_{r\to 0}\left[g(\sigma+r n_x(\sigma))-g(\sigma-r n_x(\sigma))\right],\\ \displaystyle
<g>(\sigma)=\dfrac{1}{2}\lim_{r\to 0}\left[g(\sigma+r n_x(\sigma))+g(\sigma-r n_x(\sigma))\right].
\end{cases}
\end{gather}
Above $n_x$ denotes the normal vector of $\mathcal{C}$.
\item  Consider a time-dependent interface $\mathcal{C}=\mathcal{C}(t)$, with a local parameterization $\gamma$. We assume that $\gamma_j\in \cC(I,\cC^{1+\alpha}(V_j))$ and for all $t\in I$, $\mathcal{C}(t)$ is a $\cC^{1+\alpha}$-regular, non-self-intersecting hypersurface that forms the boundary of an open, bounded, and simply connected domain $D(t)$. We define the following space:
\[
L^p( I, \;\; \cC^\alpha_{pw,\gamma}(\R^2)):=\left\{ g=g(t,x)\colon 
\begin{cases}
    \displaystyle\int_I\norm{g(t)}_{ \cC^\alpha_{pw,\gamma(t)}(\R^2)}^pdt<\infty\quad \text{ if }\quad 1\leqslant p<\infty,\\
    \displaystyle\sup_{t\in I}\text{ess} \norm{g(t)}_{ \cC^\alpha_{pw,\gamma(t)}(\R^2)}<\infty \quad \text{ if }\quad p=\infty
\end{cases} \right\}.
\]
\item 
Finally, we define the functional
\begin{gather}\label{c3.34}
\mathfrak{P}_{\gamma(t)}= \left(1+\abs{\mathcal{C}(t)}\right)\mathfrak{P}\left(\norm{\nabla\gamma(t)}_{L^\infty}+ c_{\gamma(t)}\right)\left\|\nabla \gamma(t)\right\|_{\dot \cC^{\alpha}}
\end{gather}
where $c_{\gamma(t)}$ satisfies \eqref{ep4.2}. Here, $\mathfrak{P}$ is a polynomial that is larger than those provided by \cref{ThB2} below for specific second-order Riesz operators.
\end{enumerate}
\end{definition}

In addition to the assumption on the initial density $\rho_0$ in \eqref{c5.2}-\eqref{epq30}, 
we assume the existence of an open, bounded,  and simply connected set $D(0)\subset \R^2$ such that $\mathcal{C}(0)=\partial D(0)$
(with parameterization $\gamma_0$) is a $\cC^{1+\alpha}$-regular and non-self-intersecting curve and: 
\begin{gather}\label{c5.3}
\rho_0\in \cC^\alpha_{pw,\gamma_0}(\R^2).
\end{gather}

\emph{ 
The purpose of this paper is to establish the existence of a unique  global-in-time weak solution to the system \eqref{ep4.1} with initial data \eqref{c5.1}-\eqref{c5.2}-\eqref{epq30}-\eqref{c5.3} in the spirit of the works by  Hoff \cite{hoff2002dynamics}, Hoff and Santos \cite{hoff2008lagrangean}.  
The regularity of the velocity helps to propagate 
the $\cC^{1+\alpha}$  and the non-self-intersection \eqref{ep4.2} regularities of the initial curve at all over time. We find that discontinuities in the initial density persist over time, with the jump decaying exponentially in time. The extension to density-dependent viscosity is not trivial and 
the analysis of the model is more subtle. We obtain uniqueness in a large space. On the other hand, this result supplements the work by Danchin, Fanelli and Paicu \cite{danchin2020well} with global well-posedness  even with density-dependent viscosity.}

We will now proceed with the review of known results on the propagation of discontinuity surfaces in the mathematical analysis of the Navier-Stokes equations for compressible fluids.
\subsection{Review of known results}
Since its definitive formulation in the mid-19 th century, the Navier-Stokes equations have consistently captivated the attention of numerous mathematicians. The inaugural achievement in this realm is attributed to Nash \cite{nash1962probleme} who proved a local well-posedness of strong solution in the whole $\R^3$. The density belongs to $\cC^{1+\alpha}$ while the velocity belongs 
to $\cC^{2+\alpha}$ for some $\alpha\in (0,1)$. We also refer to  Solonnikov's work  \cite{solonnikov1980solvability}, in which 
the system \eqref{ep4.1} is considered in a $\cC^2$-regular bounded domain  $\Omega\subset \R^d,\; d\in \{2,3\}$.  The initial density 
is bounded away from vacuum and belongs to $W^{1,p}(\Omega)$ for some $p>d$, whereas the 
initial velocity belongs to the Sobolev-Slobodetskii space $W^{2,1}_q(\Omega)$. Nash’s work considers heat-conducting fluids with viscosity laws that may depend on density or temperature. In contrast, Solonnikov did not account for temperature, and the viscosities are constant. The first global-in-time result is obtained by Matsumura and Nishida \cite{matsumura1979initial} for small initial data.  For constant viscosity, the initial data needs to be small in $H^3(\R^3)$,\emph{ while for non-constant viscosity, smallness in $H^4(\R^3)$ is required}. Later, the regularity requirements for the initial data were relaxed in \cite{charve2010global,Chenglobal,haspot2011existence}, allowing for small initial data in the critical Besov space \emph{in the case of constant viscosity}. However, the density remains a continuous function in space.

\emph{In the constant viscosity setting}, the classical solutions constructed in the referenced papers, come with the following energy balance: 
\begin{gather}\label{ep4.3}
    E(t)+ \mu\int_0^t \norm{\nabla u}_{L^2(\R^d)}^2+ (\mu+\lambda)\int_0^t \norm{\dvg u}_{L^2(\R^d)}= E(0)=:E_0.
\end{gather}
Above $E^H(t)$ is the energy functional defined by: 
\begin{gather}
    E(t)= \int_{\R^d} \rho(t,x) \left[\dfrac{\abs{u(t,x)}^2}{2}+ \int_{\widetilde\rho}^{\rho(t,x)} s^{-2}\left(P(s)-\widetilde P\right)ds \right]dx.
\end{gather}
In the particular case when the pressure law is of the form $P(\rho)= a \rho^\gamma$, global weak solutions
are obtained for the first time by P-L Lions \cite{lions1996mathematical}, and Feireisl, Novotn\'y, Petzeltov\'a \cite{feireisl2001existence} with some restriction on the adiabatic constant $\gamma$.  The initial data is assumed to have finite initial energy, that is $E_0<\infty$, and the solutions verify \eqref{ep4.3} with inequality.  In \cite{bresch2007existence}, Bresch and Desjardins established the existence of a global weak solution for the Navier-Stokes equation with density-dependent viscosity. \emph{However, their result requires certain Sobolev regularity assumptions on the initial density, which do not apply to our framework, as we assume the density is discontinuous across a hypersurface, with its weak gradient containing Dirac masses}.

In the last three decades, there has been significant interest in studying the propagation of discontinuity surfaces in models derived from fluid mechanics, such as the Euler or Navier-Stokes equations.
These discontinuity surfaces are sets of singularity points of certain quantities, such as vorticity for the incompressible Euler equations or density for the Navier-Stokes equations. For instance, we refer to the so-called density-patch problem proposed by P-L Lions \cite{lionsvol1} for the incompressible Navier-Stokes equations: \emph{assuming 
 $\rho_0=\mathbb{1}_{D_0}$ for some  domain  $D_0\subset \R^2$, the question is whether or not for any time $t>0$, the density is 
 $\rho(t)=\mathbb{1}_{D(t)}$, with $D(t)$ a domain with the same regularity as the initial one}.  This problem has been addressed in  \cite{danchin2012lagrangian,danchin2013incompressible,danchin2019incompressible,danchin2017persistence,gancedo2018global,gancedo2021global,liao2019global,liao2016global,liao2019globalZ,paicu2020striated}, where satisfactory solutions were provided for different regularities of $D_0$, including cases with density-dependent viscosity. 
 As far as we know, there are not enough results in the literature concerning the analogous density-patch problem 
for the Navier-Stokes equations for compressible fluids. On the one hand,  classical solutions are too regular and do not account for discontinuous initial density. On the other hand, while the weak solutions constructed by P-L Lions \cite{lions1996mathematical} or Feireisl, Novotn\'y, Petzeltv\'a \cite{feireisl2001existence} allow for discontinuous initial density,  the associated velocity is too weak to track  down density discontinuities.
As explained, for instance, in \cite{gancedo2021global}, an effective approach to tracking 
discontinuous surfaces is to construct weak 
solutions for the full model within a class that allows for the study of its dynamics. 

The initial result in this area is credited to Hoff's 2002 study \cite{hoff2002dynamics}, which is a logical follow-up to his previous results \cite{hoff1995global, hoff1995strongpoly}.
Indeed, in his pioneer work \cite{hoff1995global},  Hoff provided bounds for the following functionals (with $d\in \{2,3\}$)
\begin{gather}\label{c4.2}
\mathcal{A}_1^H(t)=\sup_{[0,t]}\sigma\norm{\nabla u}_{L^2(\R^d)}^2+\int_0^t\sigma\norm{\sqrt{\rho}\dot u}_{L^2(\R^d)}^2\quad\text{and}\quad\mathcal{A}_2^H(t)=\sup_{[0,t]}\sigma^{d}\norm{\sqrt{\rho}\dot u}_{L^2(\R^d)}^2+\int_0^t\sigma^d\norm{\nabla\dot u}_{L^2(\R^d)}^2
\end{gather}
owing to some smallness condition on the initial data. Here, $\dot{v}$ denotes the material derivative of $v$, while $\sigma$ represents a time weight. There are defined as follows:
\[
\dot v:= \dpt v+ (u\cdot \nabla) v\quad \text{and}\quad \sigma(t):=\min\{1,t\}.
\]
He observed that the so-called effective flux  defined by 
\begin{gather}
F^H:=(2\mu+\lambda)\dvg u-P(\rho)+\widetilde P
\end{gather}
and the vorticity $\rot u$ solve the following elliptic equations:
\begin{gather}\label{c4.3}
\Delta F^H=\dvg (\rho \dot u)\quad \text{and}\quad\mu\Delta \rot u= \rot (\rho \dot u).
\end{gather}
 By using the regularity of $\dot u$ provided by the functionals $\mathcal{A}_1$ and $\mathcal{A}_2$ in \eqref{c4.2}, he obtained the fact that the effective flux and vorticity belong, at least, to $L^{8/3}((1,\infty), L^\infty(\R^d))$. This finding enable  the   propagation of the $L^\infty(\R^d)$-norm of the density. As a result, he proved the existence of a global weak solution for the Navier-Stokes equations with a linear pressure law in an initial paper \cite{hoff1995global}, and later extended this to a nonlinear pressure law (gamma-law) in a subsequent paper \cite{hoff1995strongpoly}.
These weak solutions have lower regularity compared to the unique global classical solution constructed by Matsumura and Nishida in \cite{matsumura1980initial}, however they exhibit higher regularity than the solutions with finite initial energy constructed by P.-L. Lions \cite{lions1996mathematical} or by Feireisl, Novotn\'y, and Petzeltov\'a \cite{feireisl2001existence}.
Specifically, discontinuous initial densities are allowed, and the regularity of the velocity at positive times aids in tracking down the discontinuities in the density. For instance, in 2008, Hoff and Santos \cite{hoff2008lagrangean} explored the Lagrangian structure of these weak solutions. Basically, they write the velocity as sum of two terms:
\begin{align}
 u&=-\left(\dfrac{1}{2\mu+\lambda} (-\Delta)^{-1}\nabla F+ (-\Delta)^{-1}\nabla \cdot \rot u\right)+\dfrac{1}{2\mu+\lambda} (-\Delta)^{-1}\nabla (P(\rho)-\widetilde P)\nonumber\\
&=: u_F+  u_P\label{v.44}
\end{align}
The first term is at least Lipschitz at positive times, while the second one is less regular in space than the first one. Specifically, its gradient  belongs to $L^\infty((0,\infty), BMO(\R^d))$. To lower the initial time singularity of the first term, they require the initial velocity to be slightly more regular ($u_0\in H^s(\R^d)$ for $s>0$ in $d=2$ and $s>1/2$ in $d=3$). 
As a result,  the  velocity gradient belongs to $L^1_\loc([0,\infty), BMO(\R^d))$, which is sufficient to define a continuous flow map for the velocity field $u$. Consequently, continuous manifolds preserve their regularity over time.
 However, for initially $\cC^\alpha$-regular interfaces, one can only ensure  $\cC^{\alpha(t)}$-regularity, with $\alpha(t)$ decaying exponentially to zero. It is worth noting  that  they constructed a solution to the heat equation with specific initial velocity $u_0\in H^s(\R^3)$ for $s<1/2$, which has infinitely many 
integral curves approaching $x=0$ as $t$ goes to $0$.  This exponential loss of interface regularity was also observed in \cite{danchin2019compressible}. In that paper, the authors constructed global-in-time solutions with large data under the nearly incompressible assumption: the velocity divergence is assumed to be small. The velocity is relatively weak ($\nabla u\in L^1_\loc([0,\infty), BMO(\R^d))$), leading to the exponential-in-time loss of interface regularity. They proved uniqueness only for the linear pressure law case, even though the velocity field is not Lipschitz.

The decomposition of the velocity \eqref{v.44} was previously used by Hoff \cite{hoff2002dynamics} to propagate the regularity of discontinuity surfaces in $\mathbb{R}^2$. Specifically, for an initial velocity in $H^{\beta}(\mathbb{R}^2)$, both the effective flux and vorticity belong to $L^1_\loc([0,\infty), \cC^\alpha(\R^2))$, $0<\alpha<\beta<1$, and as a result, the gradient of the regular part of the velocity does as well. Assuming that initially the density is H\"older continuous on both sides of a $\cC^{1+\alpha}$ interface (with geometrical assumption \eqref{ep4.2}) across which it is discontinuous, the author showed that the second part of the velocity is at least Lipschitz. In fact, its gradient is H\"older continuous   along the tangential direction of the transported interface, guaranteeing that the latter retains the same regularity as the initial one.

 It is worth noting that the approximate density sequence is constructed  within a large space that precludes any nonlinear 
 pressure law. Recently, in \cite{zodji2023well}, we established the existence of local-in-time weak solutions for the 
 two-fluid model with density-dependent viscosity and discontinuous initial data. Notably, the regularity of these 
 local solutions is  sufficient to maintain the regularity of the interface. These solutions accommodate general 
 nonlinear pressure laws and will serve as block for the construction of global-in-time solutions in this paper.
 
In their 2020 paper, Danchin, Fanelli, and Paicu \cite{danchin2020well} proved that, assuming the initial density has tangential regularity, the less regular part of the velocity is Lipschitz. They specifically showed that $W^{2,p}$-regular
 hypersurfaces retain their regularity up to a finite time. However, they also noted that the regularity of the interface does
 not hold globally-in-time, even for small initial data. 

\emph{All of the above results pertain to the case of constant viscosity.}  When the dynamic viscosity is constant and  $\lambda(\rho)=\rho^\beta$, with $\beta>3$, the existence of a global-in-time weak solution, with no small assumption,
was pioneered by Kazhikhov and Vaigant \cite{vaigant1995existence}.  Their framework allows discontinuous density, and although not explicitly stated, the propagation of H\"older interface regularity with exponential-in-time loss also holds, with the analysis being very similar to that in \cite{hoff2008lagrangean}. However, when the viscosity $\mu$ depends on the density, there is no clear notion of effective flux, the analysis complicates and it is not even clear how can one propagate the $L^\infty$-norm of the density. In what follows, we will present some observations showing that when the viscosity $\mu$ is discontinuous, the effective flux and the vorticity lack the smoothness observed in the constant viscosity case. Specifically, we will show that these quantities are continuous at the interface only where the viscosity $\mu$ is continuous.

We first apply  the divergence and the rotational operators to the momentum equation to express:
\begin{gather}
F:=(2\mu(\rho)+\lambda(\rho))\dvg u-P(\rho)+\widetilde P= -(-\Delta)^{-1}\dvg (\rho \dot u)+ [K,\mu(\rho)-\widetilde\mu]\D u,
\end{gather}
\begin{gather}
 \mu(\rho)\rot u= -(-\Delta)^{-1}\rot(\rho \dot u)+ [K',\mu(\rho)-\widetilde\mu]\D u.
\end{gather}
We refer to \cref{prprop1} for technical details leading to these expressions. Above, $K$ and $K'$ are second-order Riesz operators  well-known to map linearly $L^p(\R^d)$ into itself for $1<p<\infty$. At the end point they map $L^\infty(\R^d)$ into the $BMO(\R^d)$ space. 

 It is straightforward to derive from the mass equation (see the computations leading to \eqref{v1} below):
\begin{gather}\label{v2}
\dpt f(\rho)+ u\cdot \nabla f(\rho)+ P(\rho)-\widetilde P= -F, \quad\text{where}\quad f(\rho)=\int_{\widetilde\rho}^\rho\dfrac{2\mu(s)+\lambda(s)}{s}ds,
\end{gather}
with the pressure term on the left-hand side interpreted as a damping term. Hence, estimates for the lower and upper bounds of the density stem from $L^\infty$-norm boundedness of $F$. A priori energy estimates (see for example \eqref{c4.2} above) provide
regularity for $\dot u$, which translates to $L^\infty(\R^2)$-norm boundedness for the first term of $F$. In contrast,  the second term of $F$, which vanishes for constant viscosity, is actually of the same order as  $\nabla u$. 
Worse, given that $K$ is not continuous over $L^\infty(\R^2)$, it is less clear whether the last term of the expression of $F$ is bounded, with the only information that $\nabla u,\, \mu(\rho)\in L^\infty(\R^2)$. This issue precludes the $L^\infty(\R^2)$- norm propagation of the density as it has been done for the isotropic case by Hoff \cite{hoff1995global}.

\vspace{0,1cm}
The second observation involves computing the  jump of the effective flux, vorticity, and velocity at the interface. The discussion follows the same lines as in \cite{novotny2016panoramas}, where the case of constant viscosity is analyzed.
First, we observe that there is a balance of forces applied to the interface, which suggests the continuity of the stress tensor in the normal direction, that is:
\begin{gather}\label{c2}
\llbracket \Pi^j\rrbracket\cdot n_x=0,\quad \text{ where }\quad \Pi^{jk}= 2\mu(\rho)\D^{jk} u+ \left(\lambda(\rho)\dvg u- P(\rho)+ \widetilde P\right)\delta^{jk},
\end{gather}
is the stress tensor, $n_x$ is the outward normal vector field of the interface and $\llbracket g\rrbracket$ denotes the jump of $g$ at the interface (see \eqref{v3}).
Next, since the velocity  is continuous in the whole space, and its gradient is continuous on both sides of the interface, then the velocity gradient is also continuous in the tangential direction of the interface. Basically,  discontinuities in  $\nabla u$ can only occur in the normal direction of the interface. In other words, there exists a vector field $\mathrm a=\mathrm a(t,\sigma)\in \R^2$ such that
\begin{gather}\label{c2.18}
\llbracket \nabla u\rrbracket= \mathrm a\cdot n^t_x.
\end{gather}
From this, one easily deduces that such vector field reads: 
\begin{gather}\label{c4.7}
\mathrm a=  (\mathrm{a}\cdot \tau_x )\tau_x+(\mathrm{a}\cdot n_x)n_x=\llbracket \rot u\rrbracket \tau_x+ \llbracket \dvg u \rrbracket n_x,
\end{gather}
where $\tau_x$ is the tangential vector field of the interface.
Using \eqref{c2.18}, we rewrite \eqref{c2} as follows:
\begin{gather}\label{eq0.1}
< \mu(\rho)>\left(\mathrm a^j +\mathrm a\cdot n_x n^j_x\right)+ 2\llbracket\mu(\rho)\rrbracket<\D^{jk} u> n^k_x+ \llbracket \lambda(\rho)\dvg u- P(\rho)\rrbracket n^j_x=0.
\end{gather}
Next, we multiply the above by $n^j_x$ before summing over $j$ to obtain:
\[
2 < \mu(\rho)> \mathrm a\cdot n_x +2\llbracket\mu(\rho)\rrbracket<\D^{jk} u> n^j_xn^k_x+  \llbracket \lambda(\rho)\dvg u- P(\rho)\rrbracket =0,
\]
and since $\mathrm a\cdot n_x= \llbracket \dvg u\rrbracket$, the jump of the effective flux reads:
\begin{gather}\label{eq0.31}
\llbracket (2\mu(\rho)+\lambda(\rho))\dvg u- P(\rho)\rrbracket= 2 \llbracket\mu(\rho)\rrbracket \left(<\dvg u>-<\D^{jk} u> n^j_xn^k_x\right).
\end{gather}
As above, we take the scalar 
product of \eqref{eq0.1} with the tangential vector $\tau_x$ and use the fact that $\llbracket\rot u\rrbracket = \mathrm a\cdot \tau_x$
 to obtain:
\begin{gather}\label{eq0.32}
\llbracket \mu(\rho)\rot u\rrbracket = \llbracket\mu(\rho)\rrbracket\left(<\rot u>-2<\D^{jk} u> n^k_x \tau^j_x\right) .
\end{gather}

It turns out that when the dynamic viscosity is continuous at the interface, for instance when it is constant, the effective flux and vorticity are also continuous at the interface. Another condition for these quantities to be continuous is  that the terms in brackets vanish, this seems not hold in general.

\emph{In view of all the above observations, it is less clear whether the effective flux and the vorticity are continuous at the interface. However, their jumps are "\textit{proportional}" to the jump in viscosity $\mu(\rho)$. As we will see in \cref{jumpdecay} below, the viscous damping of the density will cause the jump in viscosity $\mu(\rho)$ to decrease exponentially-in-time. Consequently, the jumps in effective flux, vorticity, and  velocity gradient will also decay exponentially-in-time, as observed in \cite{hoff2002dynamics,hoff2008lagrangean,novotny2016panoramas} in case of constant viscosity}. 

We are now in position to state our main result.
\subsection{Statement of the main result}

 We consider  the classical Navier-Stokes equations  \eqref{ep4.1} in two-dimensional space  with initial data  \eqref{c5.1} 
 that fulfills  \eqref{c5.2}-\eqref{epq30}-\eqref{c5.3}. The parameterization $\gamma_0$ of the  interface $\mathcal{C}(0)$  fulfills  the condition \eqref{ep4.2} with a constant $c_{\gamma_0}$. Next, we introduce some energy functionals:
 \begin{itemize}
     \item Classical energy functional: 
     \begin{gather}
     E(t)= \int_{\R^2}\left[\rho\dfrac{\abs{u}^2}{2}+ \rho\int_{\widetilde\rho}^\rho s^{-2} (P(s)-\widetilde P)ds\right](t,x)dx.
     \end{gather}
     \item Higher-order energy functionals:
     \begin{gather}\label{v.35}
     \begin{cases}
         \mathcal{A}_1(t)&=\displaystyle  \sup_{[0,t]}\norm{\nabla u}_{L^2(\R^2)}^2+\int_0^t\norm{\sqrt{\rho}\dot u}_{L^2(\R^2)}^2,\\
         \mathcal{A}_2(t)&=\displaystyle \sup_{[0,t]}\sigma\norm{\sqrt{\rho}\dot u}_{L^2(\R^2)}^2 +\int_0^t \sigma\norm{\nabla \dot u}_{L^2(\R^2)}^2,\\
         \mathcal{A}_3(t)&=\displaystyle \sup_{[0,t]}\sigma^2\norm{\nabla \dot u}_{L^2(\R^2)}^2+\int_0^t\sigma^2\norm{\sqrt\rho \ddot u}_{L^2(\R^2)}^2,
     \end{cases}
     \end{gather}
     where $\dot u= \dpt u+ (u\cdot\nabla) u$ and $\ddot u= \dpt \dot u+ (u\cdot\nabla) \dot u$ and $\sigma(t)=\min\{1,t\}$.
     \item Piecewise H\"older regularity functional:
     \begin{gather}\label{c4.37}
\text{\texttheta}(t)=\sup_{[0,t]}\norm{f(\rho)}_{\cC^\alpha_{pw,\gamma}(\R^2)}^4+\int_0^t \left[\norm{f(\rho(\tau))}_{\cC^\alpha_{pw,\gamma(\tau)}(\R^2)}^4+\sigma^{r_{\alpha}}(\tau)\norm{\nabla u(\tau)}_{\cC^\alpha_{pw,\gamma(\tau)}(\R^2)}^4\right]d\tau.
    \end{gather}
    Above,  $\gamma(\cdot)=X(\cdot) \gamma_{0}$, where $X$ is the flow of $u$;  
    \begin{gather}\label{v.36}
    r_\alpha=1+2\alpha,\quad\text{and}\quad f(\rho)=\int_{\widetilde\rho}^\rho\dfrac{2\mu(s)+\lambda(s)}{s}ds.
    \end{gather}
 \end{itemize}
 
 Recall that $P$, $\mu$, and $\lambda$  are  $\cC^2$-regular functions of the density. Additionally, we assume the existence of $a_*\in (0,\rho_{*,0}/4)$ and $a^*\in (4\rho_0^*,\infty)$ (see \eqref{epq30} for the definitions of $\rho_{0,*}$ and $\rho_0^*$) such that:
 \begin{gather}\label{v4}
P'(\rho)>0, \quad  \mu(\rho)>0, \quad \text{and}\quad\lambda(\rho)\geqslant 0, \quad \text{for all }\;\; \rho\in [a_*, a^*].
\end{gather}
 
The smallness of the initial data will be measured in the following norms:
\begin{gather}\label{ep3.1}
c_0:= \norm{u_0}_{H^1(\R^2)}^2+ \norm{\rho_0-\widetilde\rho}_{L^2(\R^2)\cap \cC^\alpha_{pw,\gamma_0}(\R^2)}^2+\norm{\llbracket \rho_0\rrbracket}_{L^4(\mathcal{C}(0))\cap L^\infty(\mathcal{C}(0))}^2.
\end{gather}

Our result reads as follows:
\begin{theo}\label[theo]{thglobal}
    Let $(\rho_0, u_0)$ be the initial data associated with the Navier-Stokes equations \eqref{ep4.1} and satisfying conditions \eqref{c5.2}, \eqref{epq30}, and \eqref{c5.3}. Additionally, assume that condition \eqref{v4} holds for the pressure and viscosity laws.

There exist constants $c>0$ and $[\mu]_0>0$ such that if:
\[
c_0\leqslant c\quad \text{ and }\quad \norm{\mu(\rho_0)-\widetilde\mu}_{\cC^\alpha_{pw,\gamma_0}(\R^2)}\leqslant [\mu]_0,
\]
then there exists a unique solution $(\rho, u)$ for the Cauchy problem associated with \eqref{ep4.1} and initial data $(\rho_0,u_0)$, satisfying:
    \begin{gather}\label{c2.10}
     E(t)+\mathcal{A}_1(t)+ \mathcal{A}_2(t)+ \mathcal{A}_3(t) +\sqrt{\text{\texttheta}(t)} \leqslant  C c_0 \quad \text{for all} \quad t\in (0,\infty).
    \end{gather}
    Above, the constant $C$ depends non-linearly on $\alpha,\rho_{0,*},\,\rho^*_0,\,\mu_{0,*},\,c_{\gamma_0},\,\norm{\nabla \gamma_0}_{\cC^\alpha},\, \norm{\nabla \vph_0}_{\cC^\alpha}$, and $\abs{\nabla\vph_0}_{\text{inf}}$.
\end{theo}
The proof of \cref{thglobal} is presented in \cref{stability} below. 
\begin{rema}
A primary challenge in this paper is deriving a Lipschitz bound for the velocity (see functional \text{\texttheta} in \eqref{c4.37}). To tackle this, we express the velocity gradient as the sum of four terms (we refer to the computations leading to \eqref{v.34}):
\begin{align}
\widetilde\mu\nabla u&=-(-\Delta)^{-1}\nabla (\rho \dot u)+ \nabla (-\Delta)^{-1}\nabla \left(\dfrac{\widetilde\mu+\lambda(\rho)}{2\mu(\rho)+\lambda(\rho)} F\right)\nonumber\\
&-\nabla (-\Delta)^{-1}\nabla \left(\dfrac{2\mu(\rho)-\widetilde\mu}{2\mu(\rho)+\lambda(\rho)}(P(\rho)-\widetilde P)\right)+\nabla(-\Delta)^{-1}\dvg ((2\mu(\rho)-\widetilde\mu)\D u)\nonumber\\
&= \nabla u_* + \nabla u_F+ \nabla u_P+ \nabla u_\delta.\label{v36}
\end{align}
\begin{itemize}
    \item The energy functionals $\mathcal{A}_1$ and $\mathcal{A}_2$, as defined in \eqref{v.35} above,
    provide sufficient regularity for $\dot u$, ensuring that the first term of the expression above, namely $\nabla u_*$,
    belongs to $L^1_\loc([0,\infty), \cC^\alpha(\R^2))$. In fact, we are not able to obtain a uniform-in-time estimate for the $L^1((1,t), L^\infty(\R^2))$-norm of $\nabla u_*$. 
    \item The other terms are less regular, and  \cref{ThB2} below
    is crucial for obtaining their  piecewise H\"older regularity. With the help of \cref{ThB2}, we obtain that 
    the last term, $\nabla u_\delta$,  is small compared to the left-hand side as long as the viscosity $\mu(\rho)$ is a small 
    perturbation of $\widetilde{\mu}$. 
    \item The piecewise H\"older norm of $F$ can be derived similarly as in the two previous steps as $F$ reads (see \eqref{v13} below):
    \[
    F=-(-\Delta)^{-1}\dvg (\rho \dot u)+ [K,\mu(\rho)-\widetilde\mu]\D u.
    \]
     Then, we make use of \cref{ThB2} to convert this bound for $F$ into a piecewise H\"older norm for $\nabla u_F$.
     \item  Given that:
     \[
     \dpt f(\rho)+ u\cdot \nabla f(\rho)+ P(\rho)-\widetilde P= -F,
     \]
     we  convert the previously obtained bound on $F$ into a piecewise H\"older norm for the pressure, and subsequently for the remaining terms of \eqref{v36}.
     \item Finally, we use a bootstrap argument to close all of these estimates, and uniform-in-time bounds are required. As mentioned in the first step, $\nabla u_*$ lacks adequate time decay, and it is unclear whether an $L^2((1,t), \cC^\alpha_{pw,\gamma}(\R^2))$-norm estimate, uniform with respect to $t>1$, can be established for $\nabla u_P$. However, we succeed in 
     obtaining a uniform-in-time estimate for $\nabla u_P$ in $L^4((0,t), \cC^\alpha_{pw,\gamma}(\R^2))$. Hence, we achieve higher time-integrability for $\nabla u_*$ by providing bound for the functional
     $\mathcal{A}_3$, since this ensures  $\nabla \dot u\in L^\infty((1, \infty), L^2(\R^2))$. A similar  functional  was derived in \cite{raphaelglobamuniquesolution} in the context of the incompressible Navier-Stokes model.
\end{itemize}
 \end{rema}
\begin{rema}
\noindent
\begin{enumerate}
    \item \cref{thglobal} generalises the works by Hoff \cite{hoff2002dynamics}, Hoff and Santos \cite{hoff2008lagrangean} by allowing nonlinear pressure law and density-dependent viscosity. We also extend the work of Danchin, Fanelli, and Paicu \cite{danchin2020well} by achieving global-in-time propagation of interface regularity.
    \item Our result  accounts for viscosity of the form $1+ \mu(\rho)$, which falls outside the Bresch-Desjardins  framework and is relevant for suspension models; see, for example, \cite{gerard2022correction}.
\end{enumerate} 
\end{rema}
\begin{rema}
     Since $r_\alpha<3$, the  velocity gradient  belongs to $L^1_\loc ([0,\infty), \cC^\alpha_{pw,\gamma}(\R^2))$,  which is sufficient to propagate the regularity of the initial curve. As a result,  the characteristics  of the interface $\gamma(t)$   exhibit exponential growth over time:
    \[
    \abs{\nabla \vph(t)}_{\text{inf}}^{-1}+c_{\gamma(t)}+\norm{\nabla\gamma(t),\, \nabla\vph(t)}_{\cC^{\alpha}}\leqslant C e^{Ct^{3/4}},
    \]
    although this growth is slower than that obtained in \cite{hoff2002dynamics}.
\end{rema}
\begin{rema}[Exponential-in-time decay of jumps]
In \cref{jumpdecay} below, we derive that $f(\rho)$, as defined in \eqref{v.36} above, verifies:
\begin{align}
\llbracket f(\rho(t,\gamma(t,s)))\rrbracket &= \llbracket f(\rho_0,\gamma_0(s))\rrbracket\nonumber\\
&\times\exp\left[\int_0^t\left[-g(\tau,s)-2h(\tau,s)
 \left(<\dvg u(\tau,\gamma(\tau,s))>-<\D^{jk} u (\tau,\gamma(\tau,s))> (n^j_xn^k_x)(\tau,s)\right)\right]d\tau \vphantom{\int_0^t}\right],
\end{align}
where $g$ and $h$ are given by:
\[
g(t,s)=\dfrac{\llbracket P(\rho(t,\gamma(t,s)))\rrbracket}{\llbracket f(\rho(t,\gamma(t,s)))\rrbracket} \quad \text{ and }\quad
h(t,s)=\dfrac{\llbracket \mu(\rho(t,\gamma(t,s)))\rrbracket}{\llbracket f(\rho(t,\gamma(t,s)))\rrbracket}.
\]
\begin{itemize}
    \item  For constant viscosity, $f$ is the logarithm function, $h=0$, and the exponential decay rate is immediate as soon as the pressure is an increasing function of the density. This observation was made by  Hoff \cite{hoff2002dynamics}, and Hoff and Santos \cite{hoff2008lagrangean}. Note that the increasing assumption on the pressure law ensures that $g$ is lower bounded away from zero.
    \item  In our context, although $h$ is no longer zero, it remains upper bounded. Therefore, by applying Young's inequality and using the $L^4((1,\infty), L^\infty(\R^2))$-norm estimate for the velocity gradient, we obtain the exponential-in-time decay  for $\llbracket f(\rho)\rrbracket$. This results in the exponential decay over time of the pressure and viscosity jump, given that $g$ and $h$ are upper bounded. This leads to the exponential decay of the vorticity and effective flux jumps over time; see \eqref{eq0.31}-\eqref{eq0.32}. As a result, the vector field $\mathrm{a}$, defined in \eqref{c4.7}, decays exponentially over time,  and so as for the jump of the velocity gradient.
    \item  Notably, if the density is initially continuous at a point $\gamma_0(s)$ of the interface, then the density, effective flux, vorticity, and velocity gradient at time $t$ are continuous at $\gamma(t,s)$ for  $t>0$.
\end{itemize}
\end{rema}
 
\paragraph*{\textbf{Outline of the paper}} The rest of this paper is structured as follows. In the next section, \cref{globalresult}, we derive an a priori estimates for local-in-time solutions. In \cref{prooofs}, we provide the proofs of the lemmas presented in the aforementioned section. The proof of the main theorem, which is a consequence of \cref{globalresult}, is the focus of \cref{stability}.
\section{Sketch of the proof of the main result}\label[section]{globalresult}
In this section, we derive a priori estimates for weak solutions for the following system:
\begin{gather}\label{epq4}
\begin{cases}
    \dpt \rho+ \dvg (\rho u)=0,\\
    \dpt (\rho u)+\dvg(\rho u\otimes u)+\nabla P(\rho)=\dvg(2\mu(\rho)\D u)+\nabla(\lambda(\rho)\dvg u).
\end{cases}
\end{gather} 
 The existence and uniqueness of such a local solution is the purpose of our recent contribution \cite{zodji2023well}, which is summarized as follows.
\begin{theo}\label[theo]{thlocal}
        Let $(\rho_0, u_0)$ be a initial data associated with the system \eqref{epq4} that satisfies the conditions \eqref{c5.2}-\eqref{epq30}-\eqref{c5.3}. Assume that the viscosity and pressure laws satisfy \eqref{v4}, and additionally, assume the compatibility condition:
    \begin{gather}\label{epq2}
        (\rho \dot u)_{|t=0}=\dvg (\Pi)_{|t=0}\in L^2(\R^2).
    \end{gather}
    There exists a positive constant $[\mu]>0$ depending only on  $\alpha$ and $\widetilde\mu$ such that if \footnote{We refer to \eqref{v5}-\eqref{c3.34} for the definition of $\ell_{\vph_0}$ and $\mathfrak{P}_{\gamma_0}$.} 
    \begin{multline}
    \left[1+\norm{\lambda(\rho_0)}_{\dot \cC^\alpha_{pw,\gamma_0}(\R^2)}+\left( \mathfrak{P}_{\gamma_0} +\ell^{-\alpha}_{\vph_0}\right)\llbracket\lambda(\rho_0)\rrbracket_{L^\infty(\mathcal{C}(0))}\right]\norm{\mu(\rho_0)-\widetilde\mu}_{\cC^\alpha_{pw,\gamma_0}(\R^2)}\\
    +\left( \mathfrak{P}_{\gamma_0} +\ell^{-\alpha}_{\vph_0}\right)\left[\norm{\llbracket\mu(\rho_0)\rrbracket}_{L^\infty(\mathcal{C}(0))}+\norm{\llbracket\mu(\rho_0)\rrbracket,\, \llbracket \lambda(\rho_0)\rrbracket}_{L^\infty(\mathcal{C}(0))}\left\|1-\dfrac{\widetilde\mu}{<\mu(\rho_0)>}\right\|_{L^\infty(\mathcal{C}(0))}\right]\leqslant \dfrac{[\mu]}{4},\label{ep3.2}
    \end{multline}
    then there exist a  time $T>0$ and a unique solution $(\rho, u)$ of the system \eqref{epq4} with initial data $(\rho_0,u_0)$, which satisfies the following:
    \begin{enumerate}
        \item  $P(\rho)-\widetilde P\in \cC([0,T], L^2(\R^2)\cap \cC^{\alpha}_{pw,\gamma}(\R^2))$, where $\gamma=\gamma(t)$ 
        is a parameterization of an $\cC^{1+\alpha}$-regular and non-self-intersecting interface $\mathcal{C}(t)$;
        \item $u\in \cC([0,T], H^1(\R^2))\cap L^\infty((0,T), \dot W^{1,6}(\R^2))\cap L^{16}((0,T), \dot W^{1,8}(\R^2))$,\, $\sigma^{\overline r/4}\nabla u\in L^4((0,T), \cC^\alpha_{pw,\gamma}(\R^2))$ for 
        \begin{gather}\label{c4.38}
            \overline r=\max\left\{\dfrac{1}{3},\;2\alpha\right\};
        \end{gather}
        \item   $\dot u\in \cC([0,T], L^2(\R^2))\cap L^2((0,T),\dot H^1(\R^2))$,\, $\sqrt{\sigma}\nabla \dot u\in L^\infty((0,T), L^2(\R^2))$,\; $\sigma^{\tfrac{1}{2}}\nabla \dot u\in L^{4}((0,T)\times \R^2)$;
        \item  $\sqrt{\sigma} \ddot u\in L^2((0,T)\times\R^2)$, $\sigma\ddot u\in L^\infty((0,T), L^2(\R^2))\cap L^2((0,T), \dot H^1(\R^2))$.
    \end{enumerate}
\end{theo}
 \begin{rema}\label[rema]{remar1}
    The velocity exhibits additional regularity. Indeed, the proof of \cref{thlocal} (specifically Remark 1.4, items 3 and 5) shows that not only is $u$  continuous throughout the entire space, but its material derivative is as well. Furthermore,
    $\ddot u$ is at least continuous across the interface $\gamma$.  
Additionally, both $\nabla u$ and $\nabla \dot u$ are H\"older  continuous on both sides of the interface $\gamma$.
\end{rema}
 \cref{thlocal} comes  with the following blow-up criterion:
\begin{coro}[Blow-up criterion]\label[coro]{thblowup}
    Let $(\rho, u)$ be the solution constructed in \cref{thlocal} defined up to a maximal time $T^*$. If 
    \begin{align}
    \limsup_{t\to T^*}&\left\{c_{\gamma(t)}+\norm{\nabla \gamma(t)}_{\cC^\alpha}+\bigg\|\dfrac{1}{\rho(t)},\;\;\dfrac{1}{\mu(\rho(t),c(t))}\bigg\|_{L^\infty(\R^2)}\right\}\nonumber\\
    +\limsup_{t\to T^*}&\left\{\norm{u(t)}_{H^1(\R^2)}+\norm{(\rho \dot u)(t)}_{L^2(\R^2)}+
    \norm{P(\rho(t))-\widetilde P}_{\cC^\alpha_{pw,\gamma(t)}(\R^2)}\right\}<\infty,\label{v.40}
    \end{align}
    and
    \begin{multline}
    \limsup_{t\to T^*}\left[1+\norm{\lambda (\rho(t))}_{\dot \cC^\alpha_{pw,\gamma(t)}(\R^2)}+\left( \mathfrak{P}_{\gamma(t)} +\ell^{-\alpha}_{\vph(t)}\right)\llbracket\lambda(\rho(t))\rrbracket_{L^\infty(\mathcal{C}(t))}\right]\norm{\mu(\rho(t))-\widetilde\mu}_{\cC^\alpha_{pw,\gamma(t)}(\R^2)}\\
    +\limsup_{t\to T^*}\left( \mathfrak{P}_{\gamma(t)} +\ell^{-\alpha}_{\vph(t)}\right)\left[\norm{\llbracket\mu(\rho(t))\rrbracket}_{L^\infty(\mathcal{C}(t))}+\norm{\llbracket\mu(\rho(t))\rrbracket,\, \llbracket \lambda(\rho(t))\rrbracket}_{L^\infty(\mathcal{C}(t))}\left\|1-\dfrac{\widetilde\mu}{<\mu(\rho(t))>}\right\|_{L^\infty(\mathcal{C}(t))}\right]<[\mu],\label{c3.38}
    \end{multline}
    then $T^*=\infty$.
\end{coro}

\cref{thlocal} and \cref{thblowup} lay the groundwork for constructing the global-in-time solution  in \cref{thglobal}.
The regularity of $u$ is sufficient in order to use $u$, $\dot u$, $\sigma \dot u$ and $\sigma^2 \ddot u$ as a test function in the following computations. Also, the a priori estimates we will derive, will involve lower regularity on the initial data; in particular, the compatibility condition \eqref{epq2} will not be required.

\subsection{Basic energy functional}
The basic energy balance is derived by taking the scalar product of the momentum equation $\eqref{epq4}_2$ with the velocity $u$, and then integrating over time and space. By doing so, we obtain:
\begin{gather}\label{eq3.3}
     E(t)+\int_0^t\int_{ \R^2}\{2\mu(\rho)\abs{\D u}^2+\lambda(\rho)\abs{\dvg u}^2\}= E(0)=E_0.
\end{gather}
Here, $E$ represents the energy functional defined as:
\[
E(t)=\int_{\R^2}\left\{\rho \dfrac{\abs{u}^2}{2}+ H_1(\rho)\right\}(t,x)dx,
\]
where $H_1$ stands for the potential energy that solves the following ODE : 
\[
\rho   H_1'(\rho)- H_1(\rho)= P(\rho)-\widetilde P\quad \text{and given by}\quad H_1(\rho)= \rho \int_{\widetilde\rho}^\rho s^{-2}(P(s)-\widetilde P)ds.
\]
More generally, as in \cite{bresch:hal-03342304}, we define  the potential energy $H_l$, $l\in (1,\infty)$,  as the  solution to  the ODE:
\[
\rho  H_l'(\rho)-H_l(\rho)=\abs{P(\rho)-\widetilde P}^{l-1}(P(\rho)-\widetilde P)
\quad
\text{which reads}
\quad
H_l(\rho)=\rho \int_{\widetilde \rho}^{\rho} s^{-2}\abs{P(s)-\widetilde P}^{l-1}(P(s)-\widetilde P)ds.
\]
These potential energies satisfy:
\begin{gather}\label{v.37}
    \dpt H_l(\rho)+\dvg (H_l(\rho) u)+ \abs{P(\rho)-\widetilde P}^{l-1}(P(\rho)-\widetilde P)\dvg u=0,
\end{gather}
and they will help in deriving the $L^{l+1}(\R^2)$-norm estimate for the pressure in the subsequent step.
\subsection{\texorpdfstring{Estimates for the functionals $\mathcal{A}_1$ and $\mathcal{A}_2$}{}}
This subsection is devoted to providing bounds for the functionals $\mathcal{A}_1$ and $\mathcal{A}_2$. These functionals yield estimates for the material acceleration $\dot u$ and, consequently, for the velocity.  We always assume the following bounds for the density and viscosity:
\begin{gather}\label{epq14}
0<\underline\rho\leqslant  \rho(t,x)\leqslant \overline\rho, \quad 0<\underline\mu\leqslant \mu(\rho(t,x))\leqslant\overline{\mu}\quad\text{and}\quad  0\leqslant \underline\lambda\leqslant \lambda(\rho(t,x))\leqslant\overline{\lambda},
\end{gather}
and we define the viscosity fluctuation:
\begin{gather}
\delta(t):=\dfrac{1}{\underline\mu}\sup_{[0,t]}\norm{\mu(\rho)-\widetilde\mu}_{L^\infty(\R^2)}.
\end{gather}
We denote by $C_*$ any constant that depends on the bounds of the density and viscosity and $\delta(t)$, and by $C_0$ any constant that depends polynomially on
\[
\norm{u_0}_{H^1(\R^2)}^2+ \norm{\rho_0-\widetilde\rho}_{L^2(\R^2)}^2
\]
These constants may change from one line to another. We derive the following estimates for functionals $\mathcal{A}_1$ and $\mathcal{A}_2$ under smallness of $\delta$.
\begin{lemm}\label[lemma]{prop1}
Suppose \eqref{epq14} holds. There exists a positive function $\text{\textkappa}=\text{\textkappa}(l)$, $l\in (1,\infty)$, such that if 
\begin{gather}\label{v.20}
\delta(t) \left(\dfrac{2\underline{\mu}+\underline\lambda}{2\overline{\mu}+\overline{\lambda}}\right)^{-\tfrac{1}{l+1}}<\text{\textkappa}(l),
\end{gather}
for $l\in \{2,3\}$,
then we have: 
\begin{gather}\label{epq16}
\mathcal{A}_1(t)
\leqslant C_* \left(C_0+\mathcal{A}_1(t)^{2}\right),
\end{gather}
\begin{gather}\label{epq17}
\mathcal{A}_2(t)\leqslant C_*\left[C_0+\mathcal{A}_1(\sigma(t))+\mathcal{A}_1(t)\left( C_0+ \mathcal{A}_1(t)\right) \right].
\end{gather}
\end{lemm}
The proof of \cref{prop1} is given in \cref{prprop1}.  The functionals $\mathcal{A}_1$ and $\mathcal{A}_2$ are under control as long as the density is bounded away 
from vacuum and upper bounded, and  the dynamic viscosity is a  small perturbation of 
the constant state $\widetilde\mu$. Achieving this control is the purpose of the subsequent steps.
\subsection{\texorpdfstring{Discussion on the propagation of the $L^\infty(\R^2)$-norm of the density}{}}
In this section, we show where the difficulty in propagating the $L^\infty$-norm of the density, as done by Hoff \cite{hoff1995global}, arises and how we circumvent this difficulty. In particular, we achieve exponential-in-time decay of jumps, which compensates for the exponential-in-time growth of the interface characteristics. 

First, we use mass equation \eqref{epq4} to derive:
\[
\dpt \log \rho +u\cdot \nabla \log \rho +\dvg u=0.
\]
Multiplying by $2\mu(\rho)+\lambda(\rho)$ and using  the expression  of the effective flux
\[
F=(2\mu(\rho)+\lambda(\rho))\dvg u- P(\rho)+\widetilde P
\]
we find (see \eqref{c4.37} for the definition of $f$) 
\[
\dpt f(\rho)+u\cdot \nabla f(\rho)+ P(\rho)-\widetilde P= -F,
\]
where the last term of the left-hand side is understood as a damping term.
As derived in \eqref{epq22} below, $F$ can also be expressed as:
\begin{gather}\label{v13}
F=- (-\Delta)^{-1}\dvg (\rho \dot u)+[K,\mu(\rho)-\widetilde \mu]\D u,
\end{gather}
where $K$ is a combination of second-order Riesz operators.   In general, the $L^\infty(\R^2)$-norm estimate of the density, or equivalently of $f(\rho)$, follows as long as we have an $L^\infty(\R^2)$-norm estimate for $F$. However,
as explained in \cite{bresch:hal-03342304}, the algebraic structure of the Navier-Stokes equations with density-dependent viscosity does not allow for such an estimate for $F$, as is often done in the isotropic case (see, for example \cite{hoff1995global}). The issue stems from the last term of $F$ (see \eqref{v13} above), which is of the same order as the velocity gradient due to the roughness of $\mu(\rho)$. Indeed, it is not clear whether the commutator $[K, \mu(\rho) - \widetilde{\mu}]$ is continuous over $L^\infty(\mathbb{R}^2)$ when the viscosity is discontinuous across a hypersurface.  In fact, this term is discontinuous even for regular velocity, as its jump corresponds exactly to the right-hand side of \eqref{eq0.31} above. In order to control this term we use the following result  which establishes that even-order singular operators are continuous over $\cC^{\alpha}_{pw,\gamma}(\R^2)$.
\begin{prop}\label[prop]{ThB2} Let the hypotheses in \cref{def1}-\cref{c3.24} hold for an interface $\mathcal{C}$ and consider a  Calder\'on-Zygmund-type singular integral operators $\mathcal{T}$ of even-order and let $p\in [1,\infty)$. There exists a constant $C=C(\alpha,p)$ such that 
for $g\in L^p(\R^2)\cap \cC^\alpha_{pw,\gamma}(\R^2)$ we have:
\begin{align}
\norm{\mathcal{T}(g)}_{\cC^\alpha_{pw,\gamma}(\R^2)}&\leqslant C\left(\norm{g}_{L^p(\R^2)}+\norm{g}_{\cC^\alpha_{pw,\gamma}(\R^2)}+\ell_\vph^{-\tfrac{1}{p}}\norm{\llbracket g\rrbracket}_{L^p(\mathcal{C})}\right)\nonumber\\
&+ C\norm{\llbracket g\rrbracket}_{L^\infty(\mathcal{C})}\left(\ell_\vph^{-\alpha}+(1+\abs{\mathcal C})\mathfrak{P}^\mathcal{T}\left(\norm{\nabla\gamma}_{L^\infty}+c_{\gamma}\right)\norm{\nabla \gamma}_{\dot \cC^\alpha}\right).
\end{align}
Above $\mathfrak{P}^\mathcal{T}$ is a polynomial depending on the kernel of $\mathcal{T}$.
\end{prop}
The proof of \cref{ThB2} follows directly from \cite[Lemma A.1 \& Lemma A.2]{zodji2023well}, so we do not present it here. 
Applying the result above with  
\begin{gather}\label{v.18}
\mathcal{C}(t)= X(t) C(0) \quad\text{and}\quad \vph(t)= \vph_0(X^{-1}(t)),
\end{gather}
 we obtain:
\begin{align}
    \norm{[K,\mu(\rho)-\widetilde \mu]&\D u(t)}_{\cC^\alpha_{pw,\gamma(t)}(\R^2)}\nonumber\\
    &\leqslant   C\norm{\mu(\rho(t))-\widetilde\mu}_{\cC^\alpha_{pw,\gamma(t)}(\R^2)}\left(\norm{\nabla u(t)}_{L^4(\R^2)}+\norm{\nabla  u(t)}_{\cC^\alpha_{pw,\gamma(t)}(\R^2)}\right)\nonumber\\
    &+C\norm{\mu(\rho(t))-\widetilde\mu}_{\cC^\alpha_{pw,\gamma(t)}(\R^2)}\left(\ell_{\vph(t)}^{-\tfrac{1}{4}}\norm{\llbracket \nabla u (t)\rrbracket}_{L^p(\mathcal{C}(t))}+\left(\ell_{\vph(t)}^{-\alpha}+\mathfrak{P}_{\gamma(t)}\right)\norm{\llbracket \nabla u(t)\rrbracket}_{L^\infty(\mathcal{C}(t))}\right)\nonumber\\
    &+C\norm{\nabla u(t)}_{L^\infty(\R^2)}\left(\ell_{\vph(t)}^{-\tfrac{1}{4}}\norm{\llbracket \mu(\rho(t)\rrbracket}_{L^p(\mathcal{C}(t))}+\left(\ell_{\vph(t)}^{-\alpha}+\mathfrak{P}_{\gamma(t)}\right)\norm{\llbracket \mu(\rho(t)\rrbracket}_{L^\infty(\mathcal{C}(t))}\right),\label{v.14}
\end{align}
for a.e $t$. The proof of
\eqref{v.14} is the purpose of the first part of  \cref{jumpdecay}. Recall  the definitions of $\ell_{\vph(t)}$ and $\mathfrak{P}_{\gamma(t)}$ given in \eqref{v5}-\eqref{c3.34}, which we will estimate in the next step.  

From \eqref{v.18}, it is straightforward to derive:
\begin{gather}\label{v.22}
\begin{cases}
\abs{\mathcal{C}(t)}&\leqslant \abs{\mathcal{C}(0)}\displaystyle \exp\left(\int_0^t \norm{\nabla u}_{L^\infty(\R^2)}\right)\\
\norm{\nabla \gamma (t)}_{L^\infty}&\leqslant \displaystyle \norm{\nabla \gamma_0}_{L^\infty} \exp\left(\int_0^t \norm{\nabla u}_{L^\infty(\R^2)}\right),\\ 
c_{\gamma(t)}&\leqslant \displaystyle c_{\gamma_0} \exp\left(\int_0^t \norm{\nabla u}_{L^\infty(\R^2)}\right),\\
\abs{\nabla\vph(t)}_{\text{inf}}&\displaystyle\geqslant  \abs{\nabla\vph_0}_{\text{inf}}\exp\left(-\int_0^t \norm{\nabla u}_{L^\infty(\R^2)}\right),
\end{cases}
\end{gather}
and 
\begin{gather}\label{v.23}
\begin{cases}
\norm{\nabla \gamma(t)}_{\dot \cC^\alpha}&\displaystyle \leqslant \left(\norm{\nabla \gamma_0}_{\dot \cC^\alpha}+\norm{\nabla\gamma_0}_{L^\infty}^{1+\alpha}\int_0^t \norm{\nabla u(\tau)}_{\dot \cC^\alpha_{pw,\gamma(\tau)}}d\tau\right)\exp\left((2+\alpha)\int_0^t\norm{\nabla u}_{L^\infty(\R^2)}\right),\\
\norm{\nabla \vph(t)}_{\dot \cC^\alpha}&\displaystyle \leqslant  \left(\norm{\nabla \vph_0}_{\dot \cC^\alpha}+\norm{\nabla\vph_0}_{L^\infty}^{1+\alpha}\int_0^t \norm{\nabla u(\tau)}_{\dot \cC^\alpha_{pw,\gamma(\tau)}}d\tau\right)\exp\left((2+\alpha)\int_0^t\norm{\nabla u}_{L^\infty(\R^2)}\right),
\end{cases}
\end{gather}
resulting in the fact that the parameters $\mathfrak{P}_{\gamma(t)}$, $\ell_{\vph(t)}^{-1}$ appearing in \eqref{v.14} grow exponentially with respect to
\begin{gather}\label{v.19}
\int_0^t \norm{\nabla u(\tau)}_{L^\infty(\R^2)}d\tau.
\end{gather}
We are unable to obtain a uniform-in-time estimate for \eqref{v.19}, and this results in exponential growth over time of  $\mathfrak{P}_{\gamma(t)}$ and $\ell_{\vph(t)}^{-1}$. To counterbalance the growth of the interface characteristics over time, the exponential-in-time decay of the viscosity and velocity gradient jumps is crucial. This leads to the following lemma.

\begin{lemm}\label[lemma]{propp3}
Suppose \eqref{epq14} holds. There are constants $0<\underline{\text{\textnu}}<\overline{\text{\textnu}}$ depending only on $\underline\rho,\, \underline\mu $ and $\overline\rho,\, \overline\mu,\, \overline\lambda$ such that the following hold true:
\begin{gather}\label{epq19}
\norm{\llbracket f(\rho(t))\rrbracket}_{L^p(\mathcal{C}(t))}\leqslant \norm{\llbracket f(\rho_0)\rrbracket}_{L^p(\mathcal{C}(0))}
\exp\left(-\underline{\text{\textnu}}\; t+(6\overline{\text{\textnu}}+1/p)\int_0^t\norm{\nabla u(\tau)}_{L^\infty(\R^2)}d \tau\right),
\end{gather}    
\begin{align}
    \norm{\llbracket \nabla u(t)\rrbracket}_{L^p(\mathcal{C}(t))}&\leqslant C_{*} \norm{\llbracket f(\rho_0)\rrbracket}_{L^p(\mathcal{C}(0))}\left(1+\norm{\nabla u(t)}_{L^\infty(\R^2)}\right)\nonumber\\
    &\times\exp\left(-\underline{\text{\textnu}}\, t+(6\overline{\text{\textnu}}+1/p)\int_0^t\norm{\nabla u(\tau)}_{L^\infty(\R^2)}d \tau\right), \label{epq20}
\end{align}
for all $1\leqslant p\leqslant \infty$. 
\end{lemm}
The proof of \cref{propp3} is given in the second part of \cref{jumpdecay}. It shows that when the pressure and viscosity laws are proportional, the constants $\underline{\text{\textnu}}$ and $\overline{\text{\textnu}}$ do not depend on $\underline\rho,\, \underline\mu $ or $\overline\rho,\, \overline\mu,\, \overline\lambda$.
 
The exponential-in-time decay of the pressure and velocity gradient jump follows immediately as long as we have a uniform-in-time $L^q((0,t), L^\infty(\R^2))$ estimate for $\sigma^{s}\nabla u$, with some $q < \infty$ and $sq'<1$. Indeed, H\"older's and Young's inequalities yield:
\begin{gather}\label{v.24}
\int_0^t\norm{\nabla u(\tau)}_{L^\infty(\R^2)}d\tau \leqslant \dfrac{\eps}{1-sq'} t+ \dfrac{1}{q(\eps q')^{q-1}} \int_0^t \sigma^{sq}\norm{\nabla u}_{L^\infty(\R^2)}^q, 
\end{gather}
for all $\eps>0$; and in virtue of \eqref{epq19}-\eqref{epq20}, we can take  
\[
\eps= (1-sq')\dfrac{\underline{\text{\textnu}}}{4 (6\overline{\text{\textnu}}+1/p)}.
\]
Achieving such a uniform-in-time estimate is the purpose of the next section.
\subsection{Final estimates}
In this section, we establish bounds for the functionals $\mathcal{A}_3$ and $\text{\texttheta}$, and we conclude by closing all the estimates. First of all, we observe that \eqref{epq4} can be rewritten as:
\begin{align*}
\widetilde\mu \Delta u&=\rho \dot u-\nabla \left((\widetilde\mu+\lambda(\rho))\dvg u-P(\rho)+\widetilde P\right)-\dvg ((2\mu(\rho)-\widetilde\mu)\D u)\\
&=\rho \dot u-\nabla\left(\dfrac{\widetilde\mu+\lambda(\rho)}{2\mu(\rho)+\lambda(\rho)} F\right) +\nabla \left(\dfrac{2\mu(\rho)-\widetilde\mu}{2\mu(\rho)+\lambda(\rho)}(P(\rho)-\widetilde P)\right)-\dvg ((2\mu(\rho)-\widetilde\mu)\D u),
\end{align*}
and therefore
\begin{align}
\widetilde\mu\nabla u&=-(-\Delta)^{-1}\nabla (\rho \dot u)+ \nabla (-\Delta)^{-1}\nabla \left(\dfrac{\widetilde\mu+\lambda(\rho)}{2\mu(\rho)+\lambda(\rho)} F\right)\nonumber\\
&-\nabla (-\Delta)^{-1}\nabla \left(\dfrac{2\mu(\rho)-\widetilde\mu}{2\mu(\rho)+\lambda(\rho)}(P(\rho)-\widetilde P)\right)+\nabla(-\Delta)^{-1}\dvg ((2\mu(\rho)-\widetilde\mu)\D u)\nonumber\\
&=: \nabla u_* + \nabla u_F+ \nabla u_P+ \nabla u_\delta.\label{v.34}
\end{align}
The last three terms above are second-order Riesz transforms  of discontinuous functions, and \cref{ThB2} plays a key role in establishing  piecewise H\"older norms estimates. The regularity of $\dot{u}$, as provided by the functionals $\mathcal{A}_1$ and $\mathcal{A}_2$, ensures that the remaining term $\nabla u_*$ is  H\"older continuous in the whole space.

 On the other hand, it is unclear whether a uniform-in-time $L^q((0,t), L^\infty(\R^2))$-norm estimate, with $q=2$, can be established for $\sigma^s\nabla u_*$ and $\sigma^s\nabla u_P$. However, such an estimate is possible for $q=4$, which explains  the time-integral in the definition of $\text{\texttheta}$, see \eqref{c4.37} above. At the same time, the functionals $\mathcal{A}_1$ and $\mathcal{A}_2$ do not provide sufficient time integrability for the material derivative $\dot{u}$ to control the $L^4((0,t), \cC^\alpha_{pw,\gamma}(\R^2))$ norm for $\sigma^s\nabla u_*$. This leads us to perform another estimate with the goal of obtaining better integrability for $\nabla \dot{u}$.  
\begin{lemm}\label[lemma]{prop2}
    Assume that the density and viscosities are bounded as in \eqref{epq14}, and suppose that \eqref{v.20}  holds for $l\in \{2,3,5\}$. Then, the following estimate holds for the functional $\mathcal{A}_3$, defined in \eqref{v.35}:
    \begin{gather}\label{epq18}
    \mathcal{A}_3(t)\leqslant C_*\left[ C_0+\mathcal{A}_3(t)^2+ \mathcal{A}_1(t)\left(1+\mathcal{A}_1(t)^3\right)+ \mathcal{A}_2(t)\left(1+\mathcal{A}_2(t)^3\right) \right].
    \end{gather}
\end{lemm}
The proof of \cref{prop2} is the focus of  \cref{prprop2} below.  The functional $\mathcal{A}_3$ provides us with 
\begin{gather}\label{epq15}
\sigma \nabla \dot u\in L^\infty((0,t), L^2(\R^2))\quad\text{and whence}\quad \dot u \in L^\infty((\sigma(t), t), L^p(\R^2))\quad \text{for all}\quad p\in (2, \infty).
\end{gather}
In \cite[Section 2.3]{hoff2002dynamics}, the author obtained the conclusion above using a different approach. Although the viscosities are constant in his analysis, he assumed a smallness condition on the kinematic viscosity $\lambda$, which does not apply in this paper.
We now proceed to derive an estimate for the functional $\text{\texttheta}$.
\begin{lemm}\label[lemma]{prop3}
Assume that \eqref{epq14} holds, and consider the functional $\text{\texttheta}$, as defined in \eqref{c4.37}. Then  we have:
\begin{align*}
\begin{cases}
    \displaystyle\int_0^t \sigma^{r_\alpha}\norm{F}_{\cC^\alpha_{pw,\gamma}(\R^2)}^4 &\displaystyle\leqslant C_*\left[ C_0+\mathcal{A}_1(t)^2(1+\mathcal{A}_1(t)^2)+\mathcal{A}_2(t)^2+\mathcal{A}_3(t)^2+ \text{\texttheta}(t)^2\right]\\
    &\displaystyle +C_*K_0e^{C_*\text{\texttheta}(t)}\norm{\llbracket f(\rho_0)\rrbracket}_{L^4\cap L^\infty(\mathcal{C}(0))}^4,\\
\text{\texttheta}(t)\displaystyle &\displaystyle\leqslant C_*\left[C_0+\mathcal{A}_1(t)^2+\left(1+\text{\texttheta}(t)+K_0
e^{C_*\text{\texttheta}(t)}\norm{\llbracket f(\rho_0) \rrbracket}_{L^4\cap L^\infty(\mathcal{C}(0))}^{4}\right)\int_0^t\sigma^{r_\alpha}\norm{F}_{\cC^\alpha_{pw,\gamma}(\R^2)}^4\right]\\
&\displaystyle+C_* e^{C_* \text{\texttheta}(t)}\left(\norm{f(\rho_0)}_{ \cC^\alpha_{pw,\gamma_0}(\R^2)}^4+K_0\norm{\llbracket f(\rho_0)\rrbracket}_{L^4\cap L^\infty(\mathcal{C}(0))}^4
+\int_0^t\sigma^{r_\alpha}\norm{F}_{\cC^\alpha_{pw,\gamma}(\R^2)}^4\right).
\end{cases}
\end{align*}
Above, $K_0$ is a constant that depends  on $c_{\gamma_0},\,\norm{\nabla \gamma_0}_{\cC^\alpha},\, \norm{\nabla \vph_0}_{\cC^\alpha}$, and $\abs{\nabla\vph_0}_{\text{inf}}$.
\end{lemm}

The proof of this proposition is given in \cref{prprop3} below.  We finally close the estimates in \cref{prop1}, \cref{prop2} and \cref{prop3} with the help of the bootstrap argument similar to the  one in \cite{bresch:hal-03342304}, and we do not present the proof. 
\begin{lemm}\label[lemma]{lem4.9}
    Let $(\rho,u)$ be a local solution to the equations  \eqref{epq4} with initial data $(\rho_0,u_0)$ that verifies \eqref{c5.2}-\eqref{epq30}-\eqref{c5.3} and the compatibility condition:
    \[
    \dvg \{2\mu(\rho_0)\D u_0+ (\lambda(\rho_0)\dvg u_0- P(\rho_0)+\widetilde P) I_d\}\in L^2(\R^2).
    \]
    Assume that the solution $(\rho,u)$ is defined up to a maximal time $T^*$.  There exist  constants $c>0$ and $[\mu]_0>0$, such that if
    \[
    C_0:=\norm{u_0}_{H^1(\R^2)}^2+ \norm{\rho_0-\widetilde\rho}_{L^2(\R^2)\cap \cC^\alpha_{pw,\gamma_0}(\R^2)}^2+\norm{\llbracket \rho_0\rrbracket}_{L^2(\mathcal{C}(0))\cap L^\infty(\mathcal{C}(0))}^2\leqslant c,
    \]
    and 
    \[
    \norm{\mu(\rho_0)-\widetilde\mu}_{\cC^\alpha_{pw,\gamma_0}(\R^2)}\leqslant [\mu]_0,
    \]
    then we have (see \eqref{v4} for the definition of $a_*$ and $a^*$): 
    \begin{gather}
    \begin{cases} \displaystyle
    0<a_*\leqslant\inf_{x\in \R^2}\rho(t,x)\leqslant \sup_{x\in \R^2}\rho(t,x)\leqslant a^*,\\
         E(t)+\mathcal{A}_1(t)+ \mathcal{A}_2(t)+ \mathcal{A}_3(t) +\sqrt{\text{\texttheta}(t)}\leqslant  C C_0,
    \end{cases}
    \end{gather}
    for all  $t\in (0, T^*)$.
    
    Above, $C$ is a constant that depends on $\alpha,\rho_{0,*},\,\rho^*_0,\,\mu_{0,*},\,c_{\gamma_0},\,\norm{\nabla \gamma_0}_{\cC^\alpha},\, \norm{\nabla \vph_0}_{\cC^\alpha}$, and  $\abs{\nabla\vph_0}_{\text{inf}}$.  
\end{lemm}
This concludes this section. The proof of \cref{thglobal} is postponed  \cref{stability} below. 
\section{Proofs}\label[section]{prooofs}
\subsection{\texorpdfstring{Proof of \cref{prop1}}{}}\label[section]{prprop1}
\dem
In this section we prove \eqref{epq16}-\eqref{epq17}.
\paragraph{\textbf {Preliminary estimates}}
The functional $\mathcal{A}_1$ appears while testing the 
momentum equation, written in the form 
\begin{gather}\label{eq3.7}
\rho \dot u^j= \dvg \Pi^j, 
\end{gather}
with $\dot u$. By doing so, we obtain:
\begin{align*}
     \int_{\R^2}\rho \abs{\dot u^j}^2&+\dfrac{d}{dt}\int_{\R^2}\left\{\mu(\rho) \abs{\D u}^2+\dfrac{\lambda(\rho)}{2}\abs{\dvg u}^2\right\}=
    - \int_{\R^2}2\mu(\rho)\D^{jk} u \partial_ku^l\partial_l u^j+ \int_{\R^2}(\rho\mu'(\rho)-\mu(\rho))\abs{\D u}^2 \dvg u\nonumber\\
    &+\dfrac{1}{2} \int_{\R^2} (\rho\lambda'(\rho)-\lambda(\rho))(\dvg u)^3- \int_{\R^2}\lambda(\rho)\dvg u\nabla u^l \partial_l u
    +\dfrac{d}{dt}\int_{\R^2}\left(\dvg u\{P(\rho)-\widetilde P\}\right)\nonumber\\
    &+ \int_{\R^2}\nabla u^l \partial_l u\left(P(\rho)-\widetilde P\right)+ \int_{\R^2}(\dvg u)^2 \left(\rho  P'(\rho)-P(\rho)+\widetilde P\right).
\end{align*}
 Integrating the above in time, we find:
\begin{gather}\label{eq3.8}
\mathcal{A}_1(t)\leqslant C_*\left(C_0+\sup_{[0,t]}\norm{P(\rho)-\widetilde P}_{L^2(\R^2)}^2+ \int_0^t \norm{\nabla u}_{L^3(\R^2)}^3+ \int_0^t\norm{P(\rho)-\widetilde P}_{L^p(\R^2)}\norm{\nabla u}_{L^{2p'}(\R^2)}^2\right),
\end{gather}
where we have used the classical energy \eqref{eq3.3} and where $p\geqslant 3$. On the other hand, we take the material derivative, $\dpt \cdot+\dvg (\cdot\, u)$,   of the momentum 
equation \eqref{eq3.7}, then multiply the resulting equation by $\dot u$,  yielding:
\begin{align*}
    \dfrac{1}{2}\dfrac{d}{dt}\int_{\R^2}\rho \abs{\dot u^j}^2+&\int_{\R^2}\left\{2\mu(\rho)\abs{\D^{jk}\dot u}^2 +\lambda(\rho)\abs{\dvg\dot u}^2\right\}=\int_{\R^2}\partial_k\dot u^j\left\{\mu(\rho)\partial_j u^l \partial_l u^k+ \mu(\rho)\partial_k u^l \partial_l u^j+2\rho\mu'(\rho)\D^{jk}u\dvg u\right\}\nonumber\\
    &+\int_{\R^2}\dvg\dot u\left\{\lambda(\rho)\nabla u^l \partial_l u+\rho\lambda'(\rho)(\dvg u)^2-\rho  P'(\rho)\dvg u\right\}-\int_{\R^2}\partial_k\dot u^j \Pi^{jk}\dvg u+\int_{\R^2}\partial_l\dot u^j \partial_k u^l \Pi^{jk}.
\end{align*}
The computations leading to the above equality can be found in \cref{Hoff2}. Next, we multiply the above 
by $\sigma(t)=\min(1,t)$ before integrating in time; then, applying H\"older's and Young's inequalities yields:
\begin{gather}\label{eq3.9}
\mathcal{A}_2(t)\leqslant C_*\left( C_0+\mathcal{A}_1(\sigma(t))+ \int_0^t\sigma \norm{\nabla u}_{L^4(\R^2)}^4+ \int_0^t\sigma \norm{P(\rho)-\widetilde P}_{L^4(\R^2)}^4\right).
\end{gather}
The remaining of this section is devoted to estimating the $L^p(\R^2)$-norms  of the gradients of the velocity and pressure as they appear in \eqref{eq3.8} and \eqref{eq3.9}. To do this, we begin by expressing the effective flux $F$ and vorticity $\rot u$ in terms of singular operators.
\paragraph{\textbf{Expression of the effective flux and vorticity}}
We apply the divergence operator to the momentum equation \eqref{eq3.7}, resulting in an elliptic equation:
\[
\dvg(\rho \dot u)=\dvg\dvg (2\mu(\rho)\D u)+\Delta\{\lambda(\rho)\dvg u-P(\rho)+\widetilde P\}
\]
from which we deduce: 
\begin{align}
\lambda(\rho)\dvg u -P(\rho)+\widetilde P&=(-\Delta)^{-1}\dvg\dvg (2\mu(\rho)\D u)-(-\Delta)^{-1}\dvg(\rho \dot u)\nonumber\\
                 &=[(-\Delta)^{-1}\dvg\dvg,\,2\mu(\rho)]\D u+2\mu(\rho)(-\Delta)^{-1}\dvg\dvg \D u-(-\Delta)^{-1}\dvg(\rho \dot u)\nonumber\\
                 &=-2\mu(\rho)\dvg u+[(-\Delta)^{-1}\dvg\dvg,\,2(\mu(\rho)-\widetilde\mu)]\D u-(-\Delta)^{-1}\dvg(\rho \dot u).\label{eq3.15}
\end{align}
Hence, we have the following expression for the effective flux:
\begin{gather}\label{epq22}
F=(2\mu(\rho)+\lambda(\rho))\dvg u-P(\rho)+\widetilde P=-(-\Delta)^{-1}\dvg(\rho \dot u)+[K,\mu(\rho)-\widetilde\mu]\D u.
\end{gather}
To express $\rot u$, we apply the rotational operator to the momentum equation  \eqref{eq3.7} to obtain:
\[
\rot_{jk} (\rho \dot u)=\rot_{jk}\dvg (2\mu(\rho)\D u),
\]
from which we deduce:
\begin{align*}
(-\Delta)^{-1}\rot_{jk} (\rho \dot u)&=(-\Delta)^{-1}\rot_{jk}\dvg(2\mu(\rho)\D u)\\
                        &=[(-\Delta)^{-1}\rot_{jk}\dvg,\, 2\mu(\rho)]\D u+2\mu(\rho)(-\Delta)^{-1}\rot_{jk}\dvg\D u\\
                        &=[(-\Delta)^{-1}\rot_{jk}\dvg,\, 2(\mu(\rho)-\widetilde\mu) ]\D u-\mu(\rho)\rot_{jk} u.
\end{align*}
It then holds that the  vorticity  reads:
\begin{gather}\label{v12}
\mu(\rho)\rot_{jk} u=-(-\Delta)^{-1}\rot_{jk} (\rho \dot u)+[K'_{jklm},\,\mu(\rho)-\widetilde\mu]\D ^{lm} u.
\end{gather}
Above, $K$ and $K'$ are combinations of second-order Riesz operators, and commutators with BMO functions are known to be continuous on $L^p(\R^2)$ for $p\in (1,\infty)$; see \cite{grafakos2009modern}. This aids in obtaining $L^p$-norm estimates for the velocity gradient and pressure in the next step.
\paragraph{\textbf{$L^p$-norm estimates for velocity gradient and pressure}}
With the help of the expressions \eqref{epq22}-\eqref{v12} above, we derive:
\begin{gather}\label{eq0.12}
\begin{cases}
\norm{F}_{L^p(\R^2)}\leqslant \kappa(p)\norm{\mu(\rho)-\widetilde\mu}_{L^\infty(\R^2)}\norm{\D u}_{L^p(\R^2)} +\norm{(-\Delta)^{-1}\dvg(\rho \dot u)}_{L^p(\R^2)},\\
\norm{\mu(\rho)\rot u}_{L^p(\R^2)}\leqslant  \kappa(p)\norm{\mu(\rho)-\widetilde\mu}_{L^\infty(\R^2)}\norm{\D u}_{L^p(\R^2)} +\norm{(-\Delta)^{-1}\rot(\rho \dot u)}_{L^p(\R^2)},
\end{cases}
\end{gather}
for all $1<p<\infty$. Consequently
\begin{align*}
\norm{\nabla u}_{L^p(\R^2)}&\leqslant  \kappa(p)\left(\norm{\dvg u}_{L^p(\R^2)}+ \norm{\rot u}_{L^p(\R^2)}\right)\\
                     &\leqslant \kappa(p)\left(\dfrac{1}{2\underline{\mu}+\underline{\lambda}}\norm{F+ P(\rho)-\widetilde P}_{L^p(\R^2)}+
                     \dfrac{1}{\underline{\mu}}\norm{\mu(\rho)\rot u}_{L^p(\R^2)}\right)\\
                     &\leqslant  \dfrac{\kappa(p)^2}{\underline{\mu}}\norm{\mu(\rho)-\widetilde\mu}_{L^\infty(\R^2)}\norm{\D u}_{L^p(\R^2)}
                     + \dfrac{\kappa(p)}{2\underline{\mu}+\underline{\lambda}}\norm{P(\rho)-\widetilde P}_{L^p(\R^2)}\\
                     &+ \dfrac{\kappa(p)}{\underline{\mu}}\left(\norm{(-\Delta)^{-1}\dvg(\rho \dot u)}_{L^p(\R^2)}+\norm{(-\Delta)^{-1}\rot(\rho \dot u)}_{L^p(\R^2)}\right).
\end{align*}
Assume that 
\begin{gather}\label{eq0.11}
\delta(t)=\dfrac{1}{\underline{\mu}}\sup_{[0,t]}\norm{\mu(\rho)-\widetilde\mu}_{L^\infty(\R^2)}<\dfrac{1}{\kappa(p)^2},
\end{gather}
 the first term of the right-hand side of the above inequality can be absorbed in the left-hand side,  yielding:
\begin{align}
\norm{\nabla u}_{L^p(\R^2)}&\leqslant  \dfrac{\kappa(p)}{\underline\mu(1-\delta \kappa(p)^2)}\left(\norm{(-\Delta)^{-1}\dvg(\rho \dot u)}_{L^p(\R^2)}+\norm{(-\Delta)^{-1}\rot(\rho \dot u)}_{L^p(\R^2)}\right)\nonumber\\
&+\dfrac{\kappa(p)}{(2\underline{\mu}+\underline{\lambda})(1-\delta \kappa(p)^2)}\norm{P(\rho)-\widetilde P}_{L^p(\R^2)}.\label{eq0.13}
\end{align}
We turn to estimating  the $L^p$-norm of the pressure. We recall that the potential energy $H_l$ satisfies (see \eqref{v.37} above):
\[
\dfrac{d}{dt}\int_{\R^2}H_l(\rho)+\int_{\R^2}\abs{P(\rho)-\widetilde P}^{l-1}(P(\rho)-\widetilde P)\dvg u=0.
\]
We substitute the divergence of the velocity in the  expression of the effective flux \eqref{epq22} as:
\[
\dvg u= (2\mu(\rho)+\lambda(\rho))^{-1}\left( F+ P(\rho)-\widetilde P\right)
\]
to obtain, after  H\"older's  inequality:
\begin{align*}
\dfrac{d}{dt}\int_{\R^2}H_l(\rho)&+ \int_{\R^2}(2\mu(\rho)+\lambda(\rho))^{-1}\abs{P(\rho)-\widetilde P}^{l+1}\\
&=-\int_{\R^2}(2\mu(\rho)+\lambda(\rho))^{-1}\abs{P(\rho)-\widetilde P}^{l-1} (P(\rho)-\widetilde P) F\\
&\leqslant \dfrac{l}{l+1}\int_{\R^2}(2\mu(\rho)+\lambda(\rho))^{-1}\abs{P(\rho)-\widetilde P}^{l+1}
+\dfrac{1}{l+1}\int_{\R^2}(2\mu(\rho)+\lambda(\rho))^{-1}\abs{F}^{l+1}.
\end{align*}
The first term in the right-hand side can be absorbed in the left-hand side and it follows: 
\begin{gather}\label{eq0.14}
    \dfrac{d}{dt}\int_{\R^2}H_l(\rho)+\dfrac{1}{l+1}\int_{\R^2}(2\mu(\rho)+\lambda(\rho))^{-1}\abs{P(\rho)-\widetilde P}^{l+1}
\leqslant \dfrac{1}{l+1}\int_{\R^2}(2\mu(\rho)+\lambda(\rho))^{-1}\abs{F}^{l+1}.
\end{gather}
To estimate the $L^{l+1}$-norm of the effective flux, we go back to \eqref{eq0.12}
and make use of \eqref{eq0.13} and obtain:
\begin{align*}
    \int_{\R^2}(2\mu(\rho)&+\lambda(\rho))^{-1}\abs{F}^{l+1}
    \leqslant \dfrac{1}{2\underline{\mu}+\underline{\lambda}}\norm{K\{(\mu(\rho)-\widetilde\mu)\D u\}}_{L^{l+1}(\R^2)}^{l+1}
    +\dfrac{1}{2\underline{\mu}+\underline{\lambda}}\int_{\R^2}\abs{(-\Delta)^{-1}\dvg(\rho \dot u)}^{l+1}\\
    &\leqslant  \dfrac{\kappa(l+1)^{l+1}}{2\underline{\mu}+\underline{\lambda}}\norm{\mu(\rho)-\widetilde\mu}_{L^\infty(\R^2)}^{l+1}\norm{\D u}_{L^{l+1}(\R^2)}^{l+1}
    +\dfrac{1}{2\underline{\mu}+\underline{\lambda}}\int_{\R^2}\abs{(-\Delta)^{-1}\dvg(\rho \dot u)}^{l+1}\\
    &\leqslant \dfrac{\kappa(l+1)^{l+1}}{2\underline{\mu}+\underline{\lambda}}\left(\dfrac{\kappa(l+1)}{2\underline{\mu}+\underline{\lambda}}\dfrac{\norm{\mu(\rho)-\widetilde\mu}_{L^\infty}}{1-\delta \kappa(l+1)^2}\right)^{l+1}\norm{P(\rho)-\widetilde P}_{L^{l+1}(\R^2)}^{l+1}\\
    &+\dfrac{\kappa(l+1)^{l+1}}{2\underline{\mu}+\underline\lambda}\left(\dfrac{\kappa(l+1)}{\underline\mu}\dfrac{\norm{\mu(\rho)-\widetilde\mu}_{L^\infty(\R^2)}}{1-\delta\kappa(l+1)^2}\right)^{l+1}\left(\norm{(-\Delta)^{-1}\dvg(\rho \dot u)}_{L^{l+1}(\R^2)}^{l+1}\right.\\
    &\left.+\norm{(-\Delta)^{-1}\rot(\rho \dot u)}_{L^{l+1}(\R^2)}^{l+1}\right)+\dfrac{1}{2\underline{\mu}+\underline{\lambda}}\norm{(-\Delta)^{-1}\dvg(\rho \dot u)}^{l+1}_{L^{l+1}(\R^2)}.
\end{align*}
We replace the above estimate  in \eqref{eq0.14} and absorb the the first term of the right-hand side  above  in the left-hand side of \eqref{eq0.14}. To achieve this, we require the following smallness assumption:
\[
\delta(t) \left(\dfrac{2\underline{\mu}+\underline\lambda}{2\overline{\mu}+\overline{\lambda}}\right)^{-\tfrac{1}{l+1}} <\dfrac{1}{3\sqrt{(l+1)}\kappa(l+1)^{2}}.
\]
In conclusion, for all $l>1$ there exists a constant $\text{\textkappa}=\text{\textkappa}(l)$ such that if 
\begin{gather}\label{eq0.15}
 \delta \left(\dfrac{2\underline{\mu}+\underline\lambda}{2\overline{\mu}+\overline{\lambda}}\right)^{-\tfrac{1}{l+1}} <\text{\textkappa}(l),
\end{gather}
then \eqref{eq0.13} holds for $p=l+1$ and additionally:
\begin{gather}\label{eq0.16}
    \dfrac{d}{dt}\int_{\R^2}H_l(\rho)+\norm{P(\rho)-\widetilde P}^{l+1}_{L^{l+1}(\R^2)}\leqslant C_*\left(\norm{(-\Delta)^{-1}\dvg(\rho \dot u)}_{L^{l+1}(\R^2)}^{l+1}+\norm{(-\Delta)^{-1}\rot(\rho \dot u)}_{L^{l+1}(\R^2)}^{l+1}\right).
\end{gather}
 To close the estimates for functionals $\mathcal{A}_1$ and $\mathcal{A}_2$, we will only need  the 
smallness condition \eqref{eq0.15} for $l\in \{2,3\}$. 
We take $l=2$ in \eqref{eq0.16} and integrate in time to obtain:
\begin{align}
    \sup_{[0,t]}\int_{\R^2}H_2(\rho)&+\int_0^t\norm{P(\rho)-\widetilde P}^{3}_{L^{3}(\R^2)}\nonumber\\
    &\leqslant C_*\left( \int_{\R^2}H_2(\rho_0)
    +\int_0^t\left(\norm{(-\Delta)^{-1}\dvg(\rho \dot u)}_{L^{3}(\R^2)}^{3}+\norm{(-\Delta)^{-1}\rot(\rho \dot u)}_{L^{3}(\R^2)}^{3}\right)\right).\label{epq23}
\end{align}
Next, we take $l=3$ in \eqref{eq0.16}, then multiply by $\sigma=\min(1,t)$ before integrating in time to obtain:
\begin{align*}
    \sup_{[0,t]}\sigma\int_{\R^2}H_3(\rho)&+\int_0^t\sigma\norm{P(\rho)-\widetilde P}^{4}_{L^{4}(\R^2)}\\
    &\leqslant C_*\left( \int_0^{\sigma(t)}\int_{\R^2}H_3(\rho)
    +\int_0^t\sigma\left(\norm{(-\Delta)^{-1}\dvg(\rho \dot u)}_{L^{4}(\R^2)}^{4}+\norm{(-\Delta)^{-1}\rot(\rho \dot u)}_{L^{4}(\R^2)}^{4}\right)\right).
\end{align*}
Combining the above estimate with \eqref{eq0.13}, Gagliardo-Nirenberg's inequality, and with the fact that 
\[
H_3(\rho)+\abs{P(\rho)-\widetilde P}^2\leqslant C_* H_1(\rho),
\]
we find:
\begin{gather}\label{eq0.24}
\int_0^t\sigma\left( \norm{\nabla u}_{L^4(\R^2)}^4+ \norm{P(\rho)-\widetilde P}_{L^4(\R^2)}^4\right)
\leqslant C_*\left[  C_0+ \mathcal{A}_1(t)\left(C_0+\mathcal{A}_1(t)\right)\right].
\end{gather}
Similarly, \eqref{eq0.13}-\eqref{epq23} imply:
\begin{gather}\label{c7}
\int_0^{\sigma(t)}\left(\norm{\nabla u}_{L^3(\R^2)}^3+ \norm{P(\rho)-\widetilde P}_{L^3(\R^2)}^3\right)\leqslant C_*\left(
 C_0+\mathcal{A}_1(t)^{1/2}\left(C_0+\mathcal{A}_1(t)\right)\right).
\end{gather}
With \eqref{eq0.24}–\eqref{c7} in hand, we can finally close the estimates for the functionals $\mathcal{A}_1$ and $\mathcal{A}_2$.
\paragraph{\textbf{Final estimates}}
We return to \eqref{eq3.9} from which we deduce: 
\[
\mathcal{A}_2(t)\leqslant C_*\left( C_0+\mathcal{A}_1(\sigma(t))+\mathcal{A}_1(t)\left( C_0+ \mathcal{A}_1(t)\right) \right).
\]
We recall the following estimate for $\mathcal{A}_1$ (see \eqref{eq3.8}):
\[
\mathcal{A}_1(t)\leqslant C_*\left( C_0+ \int_0^t \norm{\nabla u}_{L^3(\R^2)}^3+ \int_0^t\norm{P(\rho)-\widetilde P}_{L^p(\R^2)}\norm{\nabla u}_{L^{2p'}(\R^2)}^2\right)
\]
for some $p\geqslant 3$. The time integral is split  into two parts:
\[
\int_0^t=\int_0^{\sigma(t)}+\int_{\sigma(t)}^t.
\]
To bound the first term, we take $p=3$ and apply H\"older's  and Young's inequalities to obtain \eqref{c7}. For the second part,
we take  $p=4$ and similar arguments together with  \eqref{eq0.24} yield:
\begin{align}
\int_{\sigma(t)}^t\norm{\nabla u}_{L^3(\R^2)}^3 +  \int_{\sigma(t)}^t\norm{P(\rho)-\widetilde P}_{L^4(\R^2)}\norm{\nabla u}_{L^{8/3}(\R^2)}^{2}&\leqslant \int_{\sigma(t)}^t\norm{\nabla u}_{L^2(\R^2)}^2+\int_{\sigma(t)}^t\norm{\nabla u,\,P(\rho)-\widetilde P}_{L^4(\R^2)}^4\nonumber\\
&\leqslant  C_*\left( C_0+ \mathcal{A}_1(t)\left(C_0+\mathcal{A}_1(t)\right)\right). \label{c8}
\end{align}
We finally end up with:
\[
\mathcal{A}_1(t)\leqslant C_* C_0+C_*\mathcal{A}_1(t)^{1/2}\left(1+\mathcal{A}_1(t)^{1/2}\right)\left(C_0+\mathcal{A}_1(t)\right),
\]
and \eqref{epq16} follows from Young's inequality. This ends the proof of \cref{prop1}. 
\enddem
\subsection{\texorpdfstring{Proofs of \eqref{v.14}  and \cref{propp3}}{}}\label[section]{jumpdecay}
\dem[Proof of \eqref{v.14}]
We consider 
\[
\mathcal{C}(t)= X(t) \mathcal{C}(0) \quad \text{and}\quad \vph(t,x)=\vph_0(X^{-1}(t,x))
\]
where $X$ is the flow associated with the velocity $u$. 
It is clear that
\begin{align}
\norm{[K, \mu(\rho)-\widetilde\mu] \D u}_{\cC^\alpha_{pw,\gamma}(\R^2)}&\leqslant  \norm{K((\mu(\rho)-\widetilde\mu)\D u)}_{\cC^\alpha_{pw,\gamma}(\R^2)}+\norm{(\mu(\rho)-\widetilde\mu) K(\D u)}_{\cC^\alpha_{pw,\gamma}(\R^2)}\nonumber\\
&\leqslant \norm{K((\mu(\rho)-\widetilde\mu)\D u)}_{\cC^\alpha_{pw,\gamma}(\R^2)}+\norm{\mu(\rho)-\widetilde\mu}_{\cC^\alpha_{pw,\gamma}(\R^2)}\norm{ K(\D u)}_{\cC^\alpha_{pw,\gamma}(\R^2)},\label{v.15}
\end{align}
and \cref{ThB2} provides  us with:
\begin{align}
\norm{ K(\D u)(t)}_{\cC^\alpha_{pw,\gamma(t)}(\R^2)}&\leqslant C\left(\norm{\nabla u(t)}_{L^4(\R^2)}+\norm{\nabla u(t)}_{\cC^\alpha_{pw,\gamma(t)}(\R^2)}+\ell_{\vph(t)}^{-\tfrac{1}{4}}\norm{\llbracket \nabla u(t)\rrbracket}_{L^4(\mathcal{C}(t))}\right)\nonumber\\
&+ C\norm{\llbracket \nabla u(t)\rrbracket}_{L^\infty(\mathcal{C}(t))}\left(\ell_{\vph(t)}^{-\alpha}+\mathfrak{P}_{\gamma(t)}\right),\label{v.17}
\end{align}
where  $\mathfrak{P}_{\gamma(t)}$ and $\ell_{\vph(t)}$ are defined in \eqref{c3.34}-\eqref{v5}, and associated with $\mathcal{C}(t)$. We notice that $\mathfrak{P}$ is a polynomial which satisfies $\mathfrak{P}^K\leqslant \mathfrak{P}$. Similarly, we estimate the first term of the right-hand side of \eqref{v.15} as:
\begin{align}
    \norm{K((\mu(\rho)-\widetilde\mu)&\D u)(t)}_{\cC^\alpha_{pw,\gamma(t)}(\R^2)}\nonumber\\
    &\leqslant C\left(\norm{\mu(\rho(t))-\widetilde\mu}_{L^\infty(\R^2)}\norm{\nabla u(t)}_{L^4(\R^2)}+\norm{\mu(\rho(t))-\widetilde\mu}_{\cC^\alpha_{pw,\gamma(t)}(\R^2)}\norm{\nabla  u(t)}_{\cC^\alpha_{pw,\gamma(t)}(\R^2)}\right)\nonumber\\
    &+C\ell_{\vph(t)}^{-\tfrac{1}{4}}\left(\norm{\mu(\rho(t))-\widetilde\mu}_{L^\infty(\R^2)}\norm{\llbracket \nabla u (t)\rrbracket}_{L^p(\mathcal{C}(t))}+\norm{\llbracket \mu(\rho(t)\rrbracket}_{L^p(\mathcal{C}(t))}\norm{\nabla u(t)}_{L^\infty(\R^2)}\right)\nonumber\\
    &+C\left(\ell_{\vph(t)}^{-\alpha}+\mathfrak{P}_{\gamma(t)}\right)\left(\norm{\mu(\rho(t))-\widetilde\mu}_{L^\infty(\R^2)}\norm{\llbracket \nabla u(t)\rrbracket}_{L^\infty(\mathcal{C}(t))}+\norm{\llbracket \mu(\rho(t)\rrbracket}_{L^\infty(\mathcal{C}(t))}\norm{\nabla u(t)}_{L^\infty(\R^2)}\right).\label{v.16}
\end{align}
Finally \eqref{v.14} follows by summing  \eqref{v.15}-\eqref{v.17}-\eqref{v.16}.
\enddem
\dem[Proof of \cref{propp3}]
This section aims at proving \eqref{epq19}-\eqref{epq20}. We first rewrite the mass equation $\eqref{epq4}_1$ as:
\[
\dpt \log \rho + u\cdot \nabla\log \rho +\dvg u= 0.
\]
Then, we multiply the above by $2\mu(\rho)+\lambda(\rho)$ and  substitute the last term with the help of the effective flux, see \eqref{epq22}, to obtain:
\begin{gather}\label{v1}
\dpt f(\rho)+ u\cdot \nabla f(\rho)+ P(\rho)-\widetilde P= -F
\end{gather}
where $f(\rho)$ is:
\begin{gather}\label{c1.2}
f(\rho)=\int_{\widetilde\rho}^\rho \dfrac{2\mu(s)+\lambda(s)}{s}ds.
\end{gather}
 Then, for all $x\in \R^2$, we  have:
\begin{gather}\label{epq26}
\dfrac{d}{dt} f(\rho(t,X(t,x)))+ P(\rho(t,X(t,x)))-\widetilde P=-F(t,X(t,x)).
\end{gather}
In particular, along the interface $\mathcal{C}(t)$, which is parameterized by $\gamma(t,s)=X(t,\gamma_0(s))$, we have:
\begin{gather*}
\dfrac{d}{dt} f(\rho(t,\gamma(t,s)))+ P(\rho(t,\gamma(t,s)))-\widetilde P=-F(t,\gamma(t,s))
\end{gather*}
and then, by taking the jump at $\gamma(t,s)$, it holds:
\begin{align*}
\dfrac{d}{dt} \llbracket f(\rho(t,\gamma(t,s)))&\rrbracket+ g(t,s)\llbracket f(\rho(t,\gamma(t,s)))\rrbracket=-\llbracket F(t,\gamma(t,s))\rrbracket\\
&=-2 h(t,s)\llbracket f(\rho(t,\gamma(t,s)))\rrbracket \left(<\dvg u>-<\D^{jk} u> n^j_xn^k_x\right)(t,\gamma(t,s)).
\end{align*}
Above, we utilized the expression for the effective flux jump as derived in \eqref{eq0.31} above. Additionally, the functions $g$ and $h$ are defined as:
\[
g(t,s):=\dfrac{\llbracket P(\rho(t,\gamma(t,s)))\rrbracket}{\llbracket f(\rho(t,\gamma(t,s)))\rrbracket} \quad \text{ and }\quad h(t,s):=\dfrac{\llbracket \mu(\rho(t,\gamma(t,s)))\rrbracket}{\llbracket f(\rho(t,\gamma(t,s)))\rrbracket}.
\]
To achieve an exponential-in-time decay for the jump of $f(\rho)$, and subsequently for the pressure jump, we require that $g$ is both  upper bounded and bounded away from zero, while $h$ needs  simply to be upper bounded.  Specifically, there might exist two constants $\underline{\text{\textnu}}$ and $\overline{\text{\textnu}}$, potentially dependent on $\underline\rho,\, \underline\mu $ and $\overline\rho,\, \overline\mu,\, \overline\lambda$, such that for all $0 < \underline{\rho} \leqslant \rho, \, \rho' \leqslant \overline{\rho}$:
\begin{gather}\label{c1.0}
0<\underline{\text{\textnu}}\leqslant \dfrac{P(\rho)-P(\rho')}{f(\rho)-f(\rho')}\leqslant \overline{\text{\textnu}}\quad
\text{and}\quad\bigg|\dfrac{\mu (\rho)-\mu(\rho')}{f(\rho)-f(\rho')}\bigg|\leqslant\overline{\text{\textnu}}.
\end{gather}
 We observe that when the pressure and viscosity laws are proportional, the constants $\underline{\text{\textnu}}$ and $\overline{\text{\textnu}}$ in \eqref{c1.0} do not depend on the bounds of the density.  The (strict) positivity of $\underline{\text{\textnu}}$ arises from the fact that both the pressure and $f(\rho)$ are increasing functions of $\rho$.
 It then follows that 
\begin{gather*}
\dfrac{d}{dt}\left\{e^{\int_0^tg(\tau,s)d\tau}\llbracket f(\rho(t,\gamma(t,s)))\rrbracket\right\} = -2h(t,s)
 e^{\int_0^tg(\tau,s)d\tau}\llbracket f(\rho(t,\gamma(t,s)))\rrbracket\left(<\dvg u>-<\D^{jk} u> n^j_xn^k_x\right)(t,\gamma(t,s))
\end{gather*}
and whence:
\[
\abs{\llbracket f(\rho(t,\gamma(t,s)))\rrbracket}\leqslant \abs{\llbracket f(\rho_0(\gamma_0(s)))\rrbracket}\exp\left(-\underline{\text{\textnu}} t+6\overline{\text{\textnu}}\int_0^t\norm{\nabla u(\tau)}_{L^\infty(\R^2)}d \tau\right).
\]
Therefore:
\[
\norm{\llbracket f(\rho(t))\rrbracket}_{L^\infty(\mathcal{C}(t))}\leqslant \norm{\llbracket f(\rho_0,)\rrbracket}_{L^\infty(\mathcal{C}(0))}
\exp\left(-\underline{\text{\textnu}} t+6\overline{\text{\textnu}}\int_0^t\norm{\nabla u(\tau)}_{L^\infty(\R^2)}d \tau\right),
\]
and furthermore, for all $1\leqslant p< \infty$, we have:
\begin{gather}\label{eq0.43}
\norm{\llbracket f(\rho(t))\rrbracket}_{L^p(\mathcal{C}(t))}\leqslant \norm{\llbracket f(\rho_0)\rrbracket}_{L^p(\mathcal{C}(0))}
\exp\left(-\underline{\text{\textnu}} t+(6\overline{\text{\textnu}}+1/p)\int_0^t\norm{\nabla u(\tau)}_{L^\infty(\R^2)}d \tau\right).
\end{gather}
Given \eqref{c1.0}, a similar estimate applies to the jumps in viscosity $\mu(\rho)$ and pressure $P(\rho)$. Assuming that $\lambda$ also satisfies the same condition as $\mu$ in \eqref{c1.0}, a similar estimate to \eqref{eq0.43} also applies to $\llbracket \lambda(\rho) \rrbracket$.  This proves \eqref{epq19}.
Returning to \eqref{eq0.31}-\eqref{eq0.32}, we express the jumps in $\dvg u$ and $\rot u$ as follows:
\[
<2\mu(\rho)+\lambda(\rho)>\llbracket \dvg u\rrbracket=\llbracket P(\rho)\rrbracket-\llbracket\lambda(\rho)\rrbracket<\dvg u>-2\llbracket\mu(\rho)\rrbracket<\D^{jk} u> n_x^j n_x^k,
\]
\[
<\mu(\rho)>\llbracket \rot u\rrbracket=-2\llbracket\mu(\rho)\rrbracket <\D^{jk} u> n_x^k\tau_x^j.
\]
From these expressions and \eqref{eq0.43}, we deduce an estimate of the  $L^p(\mathcal{C}(t))$-norm  of $\llbracket \dvg u\rrbracket $ and $\llbracket \rot u\rrbracket $ before using \eqref{c2.18}-\eqref{c4.7} to derive:
\begin{align}
    \norm{\llbracket \nabla u(t)\rrbracket}_{L^p(\mathcal{C}(t))}&\leqslant C_* \norm{\llbracket f(\rho_0)\rrbracket}_{L^p(\mathcal{C}(0))}\left(1+\norm{\nabla u(t)}_{L^\infty(\R^2)}\right)\nonumber\\
 &\times\exp\left(-\underline{\text{\textnu}} t+(6\overline{\text{\textnu}}+1/p)\int_0^t\norm{\nabla u(\tau)}_{L^\infty(\R^2)}d \tau\right),\; p\in [1,\infty]. \label{eq0.44}
\end{align}
This ends the proof of \cref{propp3}.
\enddem
\subsection{\texorpdfstring{Proof of \cref{prop2}}{}}\label[section]{prprop2}
\dem

In this section, we will provide estimate for $\mathcal{A}_3$ as defined in \eqref{v.35}.
\paragraph{\textbf{Preliminary estimates}}
 We apply the material derivative $\dpt \cdot+ \dvg (\, \cdot\, u)$ to the momentum equations $\eqref{epq4}_2$ and find  that $\dot u$ satisfies:
\begin{gather}\label{eq0.17}
\dpt (\rho \dot u^j)+\dvg(\rho \dot u^j u)= \partial_k (\dot \Pi^{jk})+\partial_k(\Pi^{jk}\dvg u)-\dvg(\partial_k u \Pi^{jk}).
\end{gather}
We then use $\sigma^2 \ddot u$ as a test function to obtain:
\begin{align}
    \int_0^t\sigma^2\norm{\sqrt \rho \ddot u}_{L^2(\R^2)}^2&+\sigma^2(t)\int_{\R^2}\left\{\mu(\rho)\abs{\D \dot u}^2+\dfrac{\lambda(\rho)}{2}\abs{\dvg \dot u}^2\right\}= 2\int_0^{\sigma(t)}\sigma\int_{\R^2}\left\{\mu(\rho)\abs{\D \dot u}^2+\dfrac{\lambda(\rho)}{2}\abs{\dvg \dot u}^2\right\}\nonumber\\
        &-\sigma^2\int_{\R^2}\rho  P'(\rho)\dvg \dot u\dvg u+\int_0^{\sigma(t)}\sigma\int_{\R^2}\rho  P'(\rho)\dvg \dot u\dvg u+\int_0^t\sigma^2\int_{\R^2}\rho  P'(\rho)(\dvg \dot u)^2\nonumber\\
        &+\int_0^t\sigma^2\int_{\R^2} \dvg \dot u(\dvg u)^2(\rho  P'(\rho)- P(\rho)+\widetilde P)-\int_0^t\sigma^2\int_{\R^2}\partial_l\dot u^j \dvg u \partial_j u^l (\rho  P'(\rho)-P(\rho)+\widetilde P)\nonumber\\
        &+\sigma^2(t)I_1(t)-2\int_0^{\sigma(t)} \sigma I_1(s)ds+ \int_0^t\sigma^2I_2(s)ds+ \int_0^t\sigma^2I_3(s)ds,\label{c5}
\end{align}
where terms $I_1,\; I_2,\; I_3$ are
\begin{gather}\label{epq41}
I_1= \int_{\R^2}\vph(\rho)\partial_{j_1} \dot u^{j_2}\partial_{j_3} u^{j_4} \partial_{j_5} u^{j_6}
+\int_{\R^2}\partial_{j_1} \dot u^{j_2}\partial_{j_3} u^{j_4} (P(\rho)-\widetilde P),
\end{gather}
\begin{gather}\label{epq42}
I_2= \int_{\R^2}\vph(\rho)\partial_{j_1} \dot u^{j_2}\partial_{j_3}\dot u^{j_4} \partial_{j_5} u^{j_6}
+\int_{\R^2}\partial_{j_1} \dot u^{j_2}\partial_{j_3}\dot u^{j_4} (P(\rho)-\widetilde P)
+\int_{\R^2} \psi(\rho) \partial_{j_1} \dot u ^{j_2}  \partial_{j_3}  u ^{j_4} \partial_{j_5} \dot u ^{j_6},
\end{gather}
\begin{gather}\label{epq43}
I_3= \int_{\R^2}\vph(\rho)\partial_{j_1} \dot u^{j_2}\partial_{j_3} u^{j_4} \partial_{j_5} u^{j_6}\partial_{j_7} u^{j_8}
+\int_{\R^2}\partial_{j_1} \dot u^{j_2}\partial_{j_3} u^{j_4} \partial_{j_5} u^{j_6} (P(\rho)-\widetilde P),
\end{gather}
and where $\vph$ is either the viscosity $\mu$, $\lambda$, $\rho  \mu'$ , $\rho  \lambda'$, $\rho^2  \mu''$ or $\rho^2  \lambda''$;
whereas $\psi$ is either  $\rho  P'$ or $\rho^2 P''$. The computations leading to \eqref{c5} can be found in \cref{THoff}.  
In the following, we will estimate the terms appearing in the left-hand side of \eqref{c5}.

\paragraph{\textbf{Estimates for the lower-order terms}}
The first term on the right-hand side of \eqref{c5} is bounded by:
\[
2\int_0^{\sigma(t)}\sigma\int_{\R^2}\left\{\mu(\rho)\abs{\D \dot u}^2+\dfrac{\lambda(\rho)}{2}\abs{\dvg \dot u}^2\right\}
\leqslant C_*\mathcal{A}_2(\sigma(t)),
\]
while the subsequent term is estimated as:
\[
\sigma^2(t)\bigg|\int_{\R^2}\rho  P'(\rho)\dvg \dot u\dvg u\bigg|\leqslant \eta \mathcal{A}_3(t) + \dfrac{C_*}{\eta} \mathcal{A}_1(t),
\]
where $\eta$ is a small positive constant. Next, the third and fourth terms on the right-hand side of \eqref{c5} are controlled by:
\[
\bigg|\int_0^{\sigma(t)}\sigma\int_{\R^2}\rho  P'(\rho)\dvg \dot u\dvg u\bigg|+\bigg|\int_0^t\sigma^2\int_{\R^2}\rho  P'(\rho)(\dvg \dot u)^2\bigg|\leqslant  E_0+ C_* \mathcal{A}_2(t),
\]
and the next two terms are bounded by:
\begin{align}
    \bigg| \int_0^t\sigma^2\int_{\R^2} \dvg \dot u(\dvg u)^2(\rho  P'(\rho)- P(\rho)+\widetilde P)\bigg|
    &+\bigg|\int_0^t\sigma^2\int_{\R^2}\partial_l\dot u^j \dvg u \partial_j u^l (\rho  P'(\rho)-P(\rho)+\widetilde P)\bigg|\nonumber\\
    &\leqslant C_*\left[\int_0^t\sigma\norm{\nabla \dot u}_{L^2(\R^2)}^2\right]^{1/2}\left[\int_0^t\sigma\norm{\nabla u}_{L^4(\R^2)}^4\right]^{1/2}\nonumber\\
    & \leqslant C_* \mathcal{A}_2(t)^{1/2}\left(C_0+ \mathcal{A}_1(t)\left(C_0+\mathcal{A}_1(t)\right)\right)^{1/2}\nonumber\\
    &\leqslant C_* \left( C_0+ \mathcal{A}_2(t)\right)+ C_*\mathcal{A}_1(t)\left(C_0+\mathcal{A}_1(t)\right)\label{epq24},
\end{align}
where we have used \eqref{eq0.24}. Similar argument leads to:
\[
\bigg|\int_0^{\sigma(t)} sI_1(s)ds\bigg|\leqslant C_* \left(C_0+ \mathcal{A}_2(t)\right)+ C_*\mathcal{A}_1(t)\left(C_0+\mathcal{A}_1(t)\right).
\]
\paragraph{\textbf{ Estimates for $\sigma^2(t) I_1(t)$, $\displaystyle\int_0^t\sigma^2(s) I_j(s)ds$, $j\in \{2,3\}$}} 
H\"older's  and Young's inequalities imply:
\[
\bigg|\sigma^2(t)I_1(t)\bigg|\leqslant \eta \sigma^2(t) \norm{\nabla \dot u}_{L^2(\R^2)}^2 +\dfrac{C_*}{\eta} \sigma^2(t)\left(\norm{\nabla u}_{L^4(\R^2)}^4+\norm{P(\rho)-\widetilde P}_{L^4(\R^2)}^4\right).
\]
From $\abs{P(\rho)-\widetilde P}^4\leqslant C_* H_1(\rho)$, classical energy balance \eqref{eq3.3}, \eqref{eq0.13}, and Gagliardo-Nirenberg's inequality, we have:
\begin{align}
\sup_{[0,t]}\sigma\left(\norm{\nabla u}^4_{L^4(\R^2)}+\norm{P(\rho)-\widetilde P}^4_{L^4(\R^2)}\right)&\leqslant C_* C_0+ C_*\sup_{[0,t]}\sigma\left(\norm{(-\Delta)^{-1}\dvg(\rho \dot u)}_{L^4(\R^2)}^4+\norm{(-\Delta)^{-1}\rot(\rho \dot u)}_{L^4(\R^2)}^4\right)\nonumber\\
&  \leqslant C_*C_0+C_*\sup_{[0,t]}\sigma\norm{\rho \dot u}_{L^2(\R^2)}^2\left(\norm{\nabla u}_{L^2(\R^2)}^2+\norm{P(\rho)-\widetilde P}_{L^2(\R^2)}^2\right)\nonumber\\
&\leqslant C_*C_0 + C_*\mathcal{A}_2(t)\left(C_0+\mathcal{A}_1(t)\right),\label{c9}
\end{align}
yielding
\begin{gather*}
\abs{\sigma^2(t) I_1(t)}\leqslant \eta \mathcal{A}_3(t)+\dfrac{C_*}{\eta}\left(C_0 + \mathcal{A}_2(t)\left(C_0+\mathcal{A}_1(t)\right)\right).
\end{gather*}

Next, interpolation, H\"older's inequalities and \eqref{eq0.13} yield:
\begin{align}
    \bigg|\int_0^t\sigma^2I_2(s)ds\bigg|&\leqslant  C_* \int_0^t \sigma^2\norm{\nabla \dot u}_{L^{8/3}(\R^2)}^2
    \left(\norm{\nabla u}_{L^4(\R^2)}+\norm{P(\rho)-\widetilde P}_{L^4(\R^2)}\right)\nonumber\\
    &\leqslant C_* \int_0^t \sigma^2\norm{\nabla \dot u}_{L^{3}(\R^2)}^{3/2}\norm{\nabla \dot u}_{L^{2}(\R^2)}^{1/2}
    \left(\norm{\nabla u}_{L^4(\R^2)}+\norm{P(\rho)-\widetilde P}_{L^4(\R^2)}\right)\nonumber\\
    &\leqslant C_* \left[\int_0^t\sigma\left(\norm{\nabla u}_{L^4(\R^2)}^4+\norm{P(\rho)-\widetilde P}_{L^4(\R^2)}^4\right)\right]^{1/4}\left[\int_0^t\sigma\norm{\nabla \dot u}_{L^2(\R^2)}^2\right]^{1/4}\left[\int_0^t\sigma^{3}\norm{\nabla\dot u}_{L^3(\R^2)}^3\right]^{1/2}\nonumber\\
    &\leqslant C_* \mathcal{A}_2(t)^{1/4}\left(C_0+ \mathcal{A}_1(t)\left(C_0+\mathcal{A}_1(t)\right)\right)^{1/4}\left[\int_0^t\sigma^{3}\norm{\nabla\dot u}_{L^3(\R^2)}^3\right]^{1/2}\nonumber\\
    &\leqslant C_* \left(C_0+\mathcal{A}_2(t)\right)+C_* \mathcal{A}_1(t)\left(C_0+\mathcal{A}_1(t)\right)+C_*\int_0^t\sigma^{3}\norm{\nabla\dot u}_{L^3(\R^2)}^3.\label{v.21}
\end{align}
 Finally, owing to the H\"older's inequality  the remaining term is bounded as:
\begin{align*}
    \bigg|\int_0^t\sigma^2I_3(s)ds\bigg|&\leqslant C_* \int_0^t\sigma^2\norm{\nabla \dot u}_{L^2(\R^2)}\left(\norm{\nabla u}_{L^6(\R^2)}^3+\norm{P(\rho)-\widetilde P}_{L^6(\R^2)}^3\right)\\
    &\leqslant C_* \left[\int_0^t\sigma\norm{\nabla \dot u}_{L^2(\R^2)}^2\right]^{1/2}\left[\int_0^t\sigma^3\left(\norm{\nabla u}_{L^6(\R^2)}^6+\norm{P(\rho)-\widetilde P}_{L^6(\R^2)}^6\right)\right]^{1/2}\\
    &\leqslant C_* \mathcal{A}_2(t)+ C_*\int_0^t\sigma^3\left(\norm{\nabla u}_{L^6(\R^2)}^6+\norm{P(\rho)-\widetilde P}_{L^6(\R^2)}^6\right).
\end{align*}
Assuming the smallness condition for the viscosity in \eqref{eq0.15} holds for $l=5$, we multiply \eqref{eq0.16} by $\sigma$, integrate over time, and, using \eqref{eq0.13}, we arrive at:
\begin{align}
    \int_0^t\sigma\left(\norm{\nabla u}_{L^6(\R^2)}^6+\norm{P(\rho)-\widetilde P}_{L^6(\R^2)}^6\right)
    &\leqslant C_*C_0+ C_*\int_0^t\sigma\left(\norm{(-\Delta)^{-1}\dvg(\rho \dot u)}_{L^6(\R^2)}^6+\norm{(-\Delta)^{-1}\rot(\rho \dot u)}_{L^6(\R^2)}^6\right)\nonumber\\
    &\leqslant C_* C_0+C_*\int_0^t\sigma \norm{\rho \dot u}_{L^2(\R^2)}^4\left(\norm{\nabla u}_{L^2(\R^2)}^2+\norm{P(\rho)-\widetilde P}_{L^2(\R^2)}^2\right)\nonumber\\
    &\leqslant C_* C_0+ C_*\mathcal{A}_1(t)\mathcal{A}_2(t)\left(\mathcal{A}_1(t)+ C_0\right).\label{c6}
\end{align}
With this, the paragraph is concluded, and the next step is to derive an $L^3((0,t)\times \R^2)$-norm estimate for $\sigma\nabla\dot{u}$ as it appears in \eqref{v.21}.
\paragraph{\textbf{$L^3((0,t)\times\R^2)$-norm estimate for $\sigma\nabla\dot u$}}
 The approach is similar to what was done previously to estimate the $L^4((0,T)\times \R^2))$-norm of $\sigma^{1/4}\nabla u$. We rewrite \eqref{eq0.17}: 
\begin{gather}\label{eq0.18}
\rho \ddot u^j= \partial_k (\dot \Pi^{jk})+\partial_k(\Pi^{jk}\dvg u)-\dvg(\partial_k u \Pi^{jk}),
\end{gather}
and by applying the divergence operator, we express   $F_*$, defined as:
\[
 F_*:=(2\mu(\rho)+\lambda(\rho))\dvg \dot u-\lambda(\rho)\nabla u^l \partial_l u-\rho  \lambda'(\rho)(\dvg u)^2+\rho   P'(\rho)\dvg u
\]
in the following form:
\begin{align}
F_*&=-(-\Delta)^{-1}\dvg (\rho \ddot u)-(-\Delta)^{-1}\partial_{jk}\{\mu(\rho)\partial_j u^l \partial_l u^k
+\mu(\rho)\partial_k u^l \partial_l u^j+2\rho \mu'(\rho)\D^{jk}u\dvg u\}\nonumber\\
&+(-\Delta)^{-1}\partial_{jk}(\Pi^{jk}\dvg u)-(-\Delta)^{-1}\partial_{j}\dvg(\partial_k u \Pi^{jk})+ [K,\,\mu(\rho)-\widetilde\mu]\D \dot u.\label{eq0.21}
\end{align}
On the other hand, by applying the  rotational operator to \eqref{eq0.18}, we obtain that $\rot \dot u$  reads: 
\begin{multline}\label{eq0.22}
\mu(\rho)\rot \dot u=-(-\Delta)^{-1}\rot (\rho \ddot u) +K'\{(\mu(\rho)-\widetilde\mu)\D \dot u\}+(-\Delta)^{-1}\rot\partial_{k}(\Pi^{jk}\dvg u)\\-(-\Delta)^{-1}\rot\dvg(\partial_k u \Pi^{jk})
-(-\Delta)^{-1}\rot\partial_{k}\{\mu(\rho)\partial_j u^l \partial_l u^k
+\mu(\rho)\partial_k u^l \partial_l u^j+2\rho \mu'(\rho)\D^{jk}u\dvg u\}.
\end{multline}
Similarly to the argument leading to \eqref{eq0.13}, we deduce: 
\[
\norm{\nabla \dot u}_{L^3(\R^2)}\leqslant C_* \norm{(-\Delta)^{-1}\dvg (\rho \ddot u),\,(-\Delta)^{-1}\rot (\rho \ddot u)}_{L^3(\R^2)}
+C_*\norm{\nabla u,\,P(\rho)-\widetilde P}_{L^6(\R^2)}^2+ C_*\norm{\nabla u}_{L^3(\R^2)},
\]
provided that \eqref{eq0.11} holds true for $p=3$ and hence
\begin{align}
    \int_0^t\sigma^{5/2}\norm{\nabla \dot u}_{L^3(\R^2)}^3&\leqslant C_* \int_0^t\sigma^{5/2}\left(\norm{(-\Delta)^{-1}\dvg (\rho \ddot u)}_{L^3(\R^2)}^3+\norm{(-\Delta)^{-1}\rot (\rho \ddot u)}_{L^3(\R^2)}^3\right)\nonumber\\
   & +C_*\int_0^t\sigma^{5/2}\norm{\nabla u}_{L^3(\R^2)}^3
    +C_*\int_0^t\sigma^{5/2}\left(\norm{\nabla u}_{L^6(\R^2)}^6+\norm{P(\rho)-\widetilde P}_{L^6(\R^2)}^6\right).\label{eq0.20}
\end{align}
The last two terms in the inequality above are bounded in \eqref{c6}-\eqref{c7}-\eqref{c8}. Using Gagliardo-Nirenberg's inequality along with \eqref{eq0.21} and \eqref{eq0.22}, to estimate $\displaystyle 
\norm{(-\Delta)^{-1}\dvg (\rho \ddot u)}_{L^2(\R^2)}\quad \text{ and }\quad \norm{(-\Delta)^{-1}\rot (\rho \ddot u)}_{L^2(\R^2)}$, as well as \eqref{eq0.24}-\eqref{c9}, we obtain:
\begin{align*}
    \int_0^t\sigma^{5/2}\left(\right.&\left.\norm{(-\Delta)^{-1}\dvg (\rho \ddot u)}_{L^3(\R^2)}^3+\norm{(-\Delta)^{-1}\rot (\rho \ddot u)}_{L^3(\R^2)}^3\right)\\
    &\leqslant C_* \int_0^t\sigma^{5/2}\norm{\rho \ddot u}_{L^2(\R^2)}\left(\norm{(-\Delta)^{-1}\dvg (\rho \ddot u)}_{L^2(\R^2)}^2+\norm{(-\Delta)^{-1}\dvg (\rho \ddot u)}_{L^2(\R^2)}^2\right)\\
    &\leqslant C_*\int_0^t\sigma^{5/2}\norm{\rho \ddot u}_{L^2(\R^2)}\left(\norm{\nabla \dot u}_{L^2(\R^2)}^2+\norm{\nabla u,\,P(\rho)-\widetilde P}_{L^4(\R^2)}^4+\norm{\nabla u}_{L^2(\R^2)}^2\right)\\
    &\leqslant \eta \int_0^t\sigma^2\norm{\sqrt{\rho}\ddot u}_{L^2(\R^2)}^2+ \dfrac{C_*}{\eta}\int_0^t\sigma^3\left(\norm{\nabla \dot u}_{L^2(\R^2)}^4+\norm{\nabla u,\,P(\rho)-\widetilde P}_{L^4(\R^2)}^8+\norm{\nabla u}_{L^2(\R^2)}^4\right)\\
    &\leqslant \eta \int_0^t\sigma^2\norm{\sqrt{\rho}\ddot u}_{L^2(\R^2)}^2+\dfrac{C_*}{\eta}\left[\mathcal{A}_3(t)\mathcal{A}_2(t)+C_0\mathcal{A}_1(t)\right]\\
    &+\dfrac{C_*}{\eta}\left[C_0 + \mathcal{A}_2(t)\left(C_0+\mathcal{A}_1(t)\right)\right]\left[C_0 + \mathcal{A}_1(t)\left(C_0+\mathcal{A}_1(t)\right)\right].
\end{align*}
Finally, \eqref{epq18} is derived by summing all the preceding computations and choosing $\eta$ small.
\enddem

\subsection{\texorpdfstring{Proof of \cref{prop3}}{}}\label[section]{prprop3}
\dem
In this section, we will derive an estimate for the functional $\text{\texttheta}$. We start by estimating the characteristics of the interface.
\paragraph{\textbf{Estimates for the $\ell_{\vph(t)}^{-1}$ and $\mathfrak{P}_{\gamma(t)}$}} We recall the definition of $\ell_{\vph(t)}$:
\[
\ell_{\vph(t)}=\min\left\{1,\left(\dfrac{\abs{\nabla \vph(t)}_{\text{inf}}}{\norm{\nabla\vph(t)}_{\dot \cC^\alpha}}\right)^{1/\alpha}\right\},
\]
where the level-set function $\vph=\vph(t)$ satisfies:
\[
\begin{cases}
\dpt \vph+ u\cdot \nabla \vph&=0,\\
\vph_{|t=0}&= \vph_0.
\end{cases}
\]
From this, we directly derive $\eqref{v.22}_4$ and $\eqref{v.23}_2$. Moreover, \eqref{v.24} implies:
\[
\abs{\nabla \vph(t)}_{\text{inf}}\geqslant \abs{\nabla\vph_0}_{\text{inf}}\exp\left(-\eps t-\dfrac{C}{\eps}\text{\texttheta}(t)\right),
\]
and 
\begin{align*}
\norm{\nabla \vph(t)}_{\dot \cC^\alpha}&\leqslant \left[\norm{\nabla \vph_0}_{\dot \cC^\alpha}+\norm{\nabla\vph_0}_{L^\infty}\text{\texttheta}(t)^{\tfrac{1}{4}}\left(\int_0^t \sigma^{-\tfrac{r_\alpha}{3}}\right)^{3/4}\right]\exp\left(\eps t+\dfrac{C}{\eps}\text{\texttheta}(t)\right)\\
&\leqslant \left[\norm{\nabla \vph_0}_{\dot \cC^\alpha}+\norm{\nabla\vph_0}_{L^\infty}\left(\eps\, t+\dfrac{C}{\eps}\text{\texttheta}(t)\right)\right]\exp\left(\eps t+\dfrac{C}{\eps}\text{\texttheta}(t)\right)\\
&\leqslant \norm{\nabla \vph_0}_{\cC^\alpha}\exp\left(2\eps\, t +\dfrac{2C}{\eps}\text{\texttheta}(t)\right).
\end{align*}
In total, there exists a constant $C_{\vph_0}>0$ depending on the regularity of $\vph_0$, such that:
\begin{gather}\label{v.25}
\ell_{\vph(t)}^{-\alpha}  \leqslant C_{\vph_0}\exp\left(3\eps t+\dfrac{3C}{\eps}\text{\texttheta}(t)\right).
\end{gather}
This completes the estimate for $\ell_{\vph(t)}^{-1}$, and we now proceed to the estimate for $\mathfrak{P}_{\gamma(t)}$. First, recall that:
\[
\mathfrak{P}_{\gamma(t)}=\left(1+\abs{\mathcal{C}(t)}\right)\mathfrak{P}\left(\norm{\nabla\gamma(t)}_{L^\infty}+ c_{\gamma(t)}\right)\left\|\nabla \gamma(t)\right\|_{\dot \cC^{\alpha}},
\]
where $\mathfrak{P}$ is a given polynomial. By combining \eqref{v.22}-\eqref{v.23} with the computations  leading to \eqref{v.25}, we obtain:
\begin{gather*}
\mathfrak{P}_{\gamma(t)}\displaystyle\leqslant  \left(1+\abs{\mathcal{C}(0)}\right)\mathfrak{P}\left(\norm{\nabla\gamma_0}_{L^\infty}+ c_{\gamma_0}\right)  \left(\norm{\nabla \gamma_0}_{\dot \cC^\alpha}+\norm{\nabla\gamma_0}_{L^\infty}^{1+\alpha}\int_0^t \norm{\nabla u(\tau)}_{\dot \cC^\alpha_{pw,\gamma(\tau)}}d\tau\right)\exp\left(C\int_0^t\norm{\nabla u}_{L^\infty(\R^2)}\right),
\end{gather*}
which simplifies to:
\begin{gather}\label{v.26}
\mathfrak{P}_{\gamma(t)} \leqslant C_{\gamma_0}\exp\left(3\eps t+\dfrac{3C}{\eps}\text{\texttheta}(t)\right),
\end{gather}
where $C_{\gamma_0}$ is a constant that depends on the regularity of $\gamma_0$.  We now turn to the estimate for the effective flux.
\paragraph{\textbf{Estimate for $F$}}
Let us begin by recalling:
\begin{gather}\label{v.30}
F=- (-\Delta)^{-1}\dvg (\rho \dot u)+[K,\mu(\rho)-\widetilde \mu]\D u,
\end{gather}
along with the estimate for the last term above derived in \eqref{v.14}, which implies: 
\begin{align}
    \int_0^t\sigma^{r_\alpha}\norm{[K,\mu(\rho)-\widetilde \mu]&\D u}_{\cC^\alpha_{pw,\gamma}(\R^2)}^4
    \leqslant   C_*\text{\texttheta}(t)\left(\text{\texttheta}(t)+\int_0^t\sigma\norm{\nabla u}_{L^4(\R^2)}^4\right)\nonumber\\
    &+C_*\text{\texttheta}(t)\int_0^t\left[\ell_{\vph(\tau)}^{-1}\norm{\llbracket \nabla u (\tau)\rrbracket}_{L^4(\mathcal{C}(\tau))}^4+\left(\ell_{\vph(\tau)}^{-4\alpha}+\mathfrak{P}_{\gamma(\tau)}^4\right)\norm{\llbracket \nabla u(\tau)\rrbracket}_{L^\infty(\mathcal{C}(\tau))}^4\right]d\tau\nonumber\\
    &+C\text{\texttheta}(t)\sup_{[0,t]}\left[\ell_{\vph}^{-1}\norm{\llbracket \mu(\rho)\rrbracket}_{L^4(\mathcal{C})}^4+\left(\ell_{\vph}^{-4\alpha}+\mathfrak{P}_{\gamma}^4\right)\norm{\llbracket \mu(\rho)\rrbracket}_{L^\infty(\mathcal{C})}^4\right].\label{v.27}
\end{align}
From \eqref{epq19}, we find:
\[
\begin{cases}
\norm{\llbracket \mu(\rho(\tau))\rrbracket}_{L^4(\mathcal{C}(\tau))}&\leqslant \displaystyle  C_*\norm{\llbracket f(\rho_0)\rrbracket}_{L^4(\mathcal{C}(0))}
\exp\left(-\dfrac{\underline{\text{\textnu}}}{2}\; \tau+ C_*\text{\texttheta}(\tau)\right),\\
\norm{\llbracket \mu(\rho(\tau))\rrbracket}_{L^\infty(\mathcal{C}(\tau))}&\leqslant \displaystyle C_* \norm{\llbracket f(\rho_0)\rrbracket}_{L^\infty(\mathcal{C}(0))}
\exp\left(-\dfrac{\underline{\text{\textnu}}}{2}\; \tau+ C_*\text{\texttheta}(\tau)\right), 
\end{cases}
\]
and hence, see \eqref{v.25}-\eqref{v.26}:
\begin{align*}
    \left[\ell_{\vph(\tau)}^{-1}\norm{\llbracket \mu(\rho(\tau))\rrbracket}_{L^4(\mathcal{C}(\tau)}^4\right.&\left.+\left(\ell_{\vph(\tau)}^{-4\alpha}+\mathfrak{P}_{\gamma(\tau)}^4\right)\norm{\llbracket \mu(\rho(\tau)\rrbracket}_{L^\infty(\mathcal{C}(\tau)}^4\right]\\
    &\leqslant C_*C_{\vph_0}^{\tfrac{1}{\alpha}}\norm{\llbracket f(\rho_0)\rrbracket}_{L^4(\mathcal{C}(0))}^4\exp\left[\left(\dfrac{3\eps}{\alpha}-2\underline{\text{\textnu}}\right)\tau+\dfrac{C_*}{\eps}\text{\texttheta}(\tau)\right]\\
    &+C_*\left(C_{\vph_0}^4+C_{\gamma_0}^4\right)\norm{\llbracket f(\rho_0)\rrbracket}_{L^\infty(\mathcal{C}(0))}^4\exp\left[\left(12\eps-2\underline{\text{\textnu}}\right)\tau+\dfrac{C_*}{\eps}\text{\texttheta}(\tau)\right].
\end{align*}
By setting $\eps=\alpha\underline{\text{\textnu}}/6$, it follows that:
\begin{gather}\label{v.28}
\sup_{[0,t]}\left[\ell_{\vph}^{-1}\norm{\llbracket \mu(\rho)\rrbracket}_{L^4(\mathcal{C})}^4+\left(\ell_{\vph}^{-4\alpha}+\mathfrak{P}_{\gamma}^4\right)\norm{\llbracket \mu(\rho)\rrbracket}_{L^\infty(\mathcal{C})}^4\right]\leqslant C_*K_0e^{C_*\text{\texttheta}(t)}\norm{\llbracket f(\rho_0)\rrbracket}_{L^4\cap L^\infty(\mathcal{C}(0))}^4,
\end{gather}
where $K_0$ depends polynomially on $C_{\gamma_0}$ and $C_{\vph_0}$. Following the same computations, we find from \eqref{epq20}-\eqref{v.25}-\eqref{v.26}:
\begin{gather}\label{v.29}
\int_0^t\left[\ell_{\vph(\tau)}^{-1}\norm{\llbracket \nabla u (\tau)\rrbracket}_{L^4(\mathcal{C}(\tau))}^4+\left(\ell_{\vph(\tau)}^{-4\alpha}+\mathfrak{P}_{\gamma(\tau)}^4\right)\norm{\llbracket \nabla u(\tau)\rrbracket}_{L^\infty(\mathcal{C}(\tau))}^4\right]d\tau\leqslant C_*K_0e^{C_*\text{\texttheta}(t)}\norm{\llbracket f(\rho_0)\rrbracket}_{L^4\cap L^\infty(\mathcal{C}(0))}^4.
\end{gather}
Summing up \eqref{v.27}-\eqref{v.28}-\eqref{v.29} and \eqref{eq0.24}, we have the following estimate for the last term of \eqref{v.30}:
\begin{gather}\label{v.32}
    \int_0^t\sigma^{r_\alpha}\norm{[K,\mu(\rho)-\widetilde \mu]\D u}_{\cC^\alpha_{pw,\gamma}(\R^2)}^4
    \leqslant C_*\text{\texttheta}(t)\left(C_0+\text{\texttheta}(t)+ \mathcal{A}_1(t)^2\right)+C_*K_0\text{\texttheta}(t)
    e^{C_*\text{\texttheta}(t)}\norm{\llbracket f(\rho_0)\rrbracket}_{L^4\cap L^\infty(\mathcal{C}(0))}^4.
\end{gather}
Next, we estimate the first term of \eqref{v.30} in terms of the functionals $\mathcal{A}_1$, $\mathcal{A}_2$, and $\mathcal{A}_3$.  The embedding inequality implies that:
\begin{align*}
    \int_0^t\sigma^{1+2\alpha}\norm{(-\Delta)^{-1}\dvg (\rho \dot u)}_{\dot\cC^\alpha(\R^2)}^4&\leqslant C\int_0^t\sigma^{1+2\alpha}\norm{\rho \dot u}_{L^{2/(1-\alpha)}(\R^2)}^{4}\\
    &\leqslant C_*\int_0^t\sigma^{1+2\alpha}\norm{\nabla \dot u}_{L^{2}(\R^2)}^{4\alpha}\norm{\dot u}_{L^{2}(\R^2)}^{4(1-\alpha)}\\
    &\leqslant C_*\left(\int_{0}^t\sigma^{3} \norm{\nabla\dot u}_{L^2(\R^2)}^{4}\right)^{\alpha}\left(\int_0^t \sigma\norm{\dot u}_{L^2(\R^2)}^{4}\right)^{1-\alpha},
\end{align*}
and whence:
\begin{gather}\label{v.33}
\int_0^t\sigma^{1+2\alpha}\norm{(-\Delta)^{-1}\dvg (\rho \dot u)}_{\dot\cC^\alpha(\R^2)}^4\leqslant C_*\left(\mathcal{A}_1(t)^2+\mathcal{A}_2(t)^2+\mathcal{A}_3(t)^2\right).
\end{gather}
Secondly, we have:
\begin{align}
    \int_0^t\sigma\norm{(-\Delta)^{-1}\dvg (\rho \dot u)}_{L^\infty(\R^2)}^4&\leqslant C\int_0^t\sigma\norm{\rho \dot u}_{L^3(\R^2)}^3\norm{(-\Delta)^{-1}\dvg (\rho \dot u)}_{L^2(\R^2)}\nonumber\\
    & \leqslant C_*\int_0^t\sigma\norm{\nabla \dot u}_{L^{2}(\R^2)}\norm{\dot u}_{L^{2}(\R^2)}^2\norm{\nabla u,\, P(\rho)-\widetilde P}_{L^{2}(\R^2)}\nonumber\\
     &\leqslant C_*\left[\int_0^t\sigma\norm{\nabla \dot u}_{L^2(\R^2)}^2\right]^{\tfrac{1}{2}}\left[\int_0^t\sigma\norm{\dot u}_{L^2(\R^2)}^4\norm{\nabla u,\, P(\rho)-\widetilde P}_{L^{2}(\R^2)}^2\right]^{\tfrac{1}{2}}\nonumber\\
     &\leqslant C_*\left( C_0+\mathcal{A}_1(t)^2+\mathcal{A}_2(t)^2\right).\label{v.31}
\end{align}
Summing up \eqref{v.30}-\eqref{v.32}-\eqref{v.33}-\eqref{v.31}, we conclude that for $r_\alpha=1+2\alpha$:
\begin{align}
    \int_0^t \sigma^{r_\alpha}\norm{F}_{\cC^\alpha_{pw,\gamma}(\R^2)}^4&\leqslant C_*\left[ C_0+\mathcal{A}_1(t)^2+\mathcal{A}_2(t)^2+\mathcal{A}_3(t)^2+ \text{\texttheta}(t)\left(C_0+\text{\texttheta}(t)+ \mathcal{A}_1(t)^2\right)\right]\nonumber\\
    &+C_*K_0
    e^{C_*\text{\texttheta}(t)}\norm{\llbracket f(\rho_0)\rrbracket}_{L^4\cap L^\infty(\mathcal{C}(0))}^4.\label{c2.5}
\end{align}
This completes the estimate for the effective flux, from which we derive the pressure estimate.
\paragraph{\textbf{Estimate for the pressure}}
We begin by recalling equation \eqref{epq26}, which gives us the following expression:
\[
f(\rho(\tau,X(\tau,x))= f(\rho_0(x))e^{-\int_0^\tau g_1(\tau',x)d\tau'}-\int_0^\tau e^{-\int_{\tau'}^\tau g_1(\tau'',x)d\tau''} F(\tau',X(\tau',x))d\tau',
\]
where $g_1$ is defined as:
\[
g_1(t,x):=\dfrac{P(\rho(t,X(t,x)))-\widetilde P}{f(\rho(t,X(t,x)))}\in [\underline{\text{\textnu}},\,\overline{\text{\textnu}}].
\]
This leads to the following bound:
\begin{gather}\label{eq0.35}
\norm{f(\rho(\tau))}_{L^\infty(\R^2)}\leqslant  e^{-\underline{\text{\textnu}}\, \tau}\norm{f(\rho_0)}_{L^\infty(\R^2)} 
+ \int_0^\tau e^{\underline{\text{\textnu}}\,(\tau'-\tau)} 
\norm{F(\tau')}_{L^\infty(\R^2)} d\tau'.
\end{gather}
We express the last term above as:
\[
\int_0^\tau e^{\underline{\text{\textnu}}\,(\tau'-\tau)} \norm{F(\tau')}_{L^\infty(\R^2)} d\tau'=\int_0^{\tau} e^{\underline{\text{\textnu}}\,(\tau'-\tau)} \norm{F(\tau')}_{L^\infty(\R^2)}\mathbb{1}_{\{\tau'< 1\}} d\tau'+\int_{0}^\tau e^{\underline{\text{\textnu}}\,(\tau'-\tau)} \norm{F(\tau')}_{L^\infty(\R^2)} \mathbb{1}_{\{\tau'> 1\}}d\tau',
\]
where the first term is bounded as:
\[
\int_0^{\tau}e^{\underline{\text{\textnu}}\,(\tau'-\tau)} \norm{F(\tau')}_{L^\infty(\R^2)}\mathbb{1}_{\{\tau'< 1\}} d\tau' 
\leqslant \left(\dfrac{3}{3-r_\alpha}\right)^{3/4}\norm{\sigma^{r_\alpha/4} F}_{L^4((0,\sigma(t)), L^\infty(\R^2))},
\]
and the second term is bounded as:
\[
\int_0^{\tau}e^{\underline{\text{\textnu}}\,(\tau'-\tau)} \norm{F(\tau')}_{L^\infty(\R^2)}\mathbb{1}_{\{\tau'> 1\}} d\tau' 
\leqslant \left(\dfrac{3}{4\underline{\text{\textnu}}}\right)^{3/4}\norm{ F}_{L^4((\sigma(t),t), L^\infty(\R^2))}.
\]
Thus, we obtain:
\begin{gather}\label{c1.6}
\sup_{[0,t]}\norm{f(\rho)}_{L^\infty(\R^2)}^4\leqslant  \norm{f(\rho_0)}_{L^\infty(\R^2)}^4 + C_*\int_0^t\sigma^{r_\alpha}\norm{F}_{L^\infty(\R^2)}^4.
\end{gather}
Moreover, by applying Young's inequality for convolution, we have:
\[
\bigg\| \int_0^{\tau}e^{\underline{\text{\textnu}}\,(\tau'-\tau)} \norm{F(\tau')}_{L^\infty(\R^2)}\mathbb{1}_{\{\tau'< 1\}} d\tau'\bigg\|_{L^4(0,t)}\leqslant \left(\dfrac{1}{4\underline{\text{\textnu}}}\right)^{1/4}\left(\dfrac{3}{3-r_\alpha}\right)^{3/4}\norm{\sigma^{r_\alpha/4} F}_{L^4((0,\sigma(t)), L^\infty(\R^2))},
\]
and
\[
\bigg\|\int_0^{\tau}e^{\underline{\text{\textnu}}\,(\tau'-\tau)} \norm{F(\tau')}_{L^\infty(\R^2)}\mathbb{1}_{\{\tau'> 1\}} d\tau' \bigg\|_{L^4(0,t)}
\leqslant \dfrac{1}{\underline{\text{\textnu}}}\norm{ F}_{L^4((\sigma(t),t), L^\infty(\R^2))}.
\]
As a result, we obtain the following estimate:
\begin{gather}\label{c2.2}
\int_{0}^t \norm{f(\rho(\tau))}_{L^\infty(\R^2)}^4 d\tau\leqslant C_*\left( \norm{f(\rho_0)}_{L^\infty(\R^2)}^4
+\int_{0}^t\sigma^{r_\alpha}\norm{F(\tau)}_{L^\infty(\R^2)}^4d\tau\right).
\end{gather}
With the $L^\infty(\R^2)$-norm estimate for $f(\rho)$ now complete, we proceed to estimating the $\dot \cC^\alpha_{pw,\gamma}(\R^2)$-norm of $f(\rho)$. To this end, we consider two points $x_i^0$, $i \in \{1, 2\}$, located on the same side of the interface, and define $x_i(t) = X(t, x_i^0)$. We infer from \eqref{epq26} that:
\begin{gather}\label{epq27}
\dfrac{d}{dt} f(\rho(t,x_i(t)))\Bigg|_{i=1}^{i=2} + f( \rho(t,x_i(t)))\Bigg|_{i=1}^{i=2}g_2(t,x_1(t),x_2(t))=-F(t,x_i(t))\Bigg|_{i=1}^{i=2},
\end{gather}
where $g_2$ is: 
\begin{gather}\label{epq29}
 g_2(t,x,y):=\dfrac{P(\rho(t,x))-P(\rho(t,y))}{f(\rho(t,x))-f(\rho(t,y))}\in [\underline{\text{\textnu}},\, \overline{\text{\textnu}}].
\end{gather}
Integrating this differential equation yields:
\begin{align}
f(\rho(\tau,x_i(\tau)))\bigg|_{i=1}^{i=2}&=  f(\rho_0(x_i^0))\bigg|_{i=1}^{i=2} e^{-\int_0^\tau g_2(\tau',x_2(\tau'),x_1(\tau'))d\tau'}\nonumber\\
&-\int_0^\tau e^{-\int_{\tau'}^\tau g_2(\tau'',x_2(\tau''),x_1(\tau''))d\tau''}F(\tau',x_i(\tau' ))\bigg|_{i=1}^{i=2}d\tau'.\label{epq28}
\end{align}
It is straightforward to obtain, for all $0\leqslant \tau'\leqslant \tau$: 
\[
\abs{x_2(\tau')-x_1(\tau')}\leqslant e^{\int_{\tau'}^{\tau}\norm{\nabla u(\tau'')}_{L^\infty(\R^2)}d\tau''} \abs{x_2(\tau)-x_1(\tau)},
\]
and  hence \eqref{epq28} implies:
\begin{align}
\norm{f(\rho(\tau))}_{\dot \cC^\alpha_{pw,\gamma(\tau)}(\R^2)}&\leqslant e^{-\underline{\text{\textnu}}\tau +\int_0^{\tau}\norm{\nabla u(\tau')}_{L^\infty(\R^2)}d\tau'}\norm{f(\rho_0)}_{\dot \cC^\alpha_{pw,\gamma_0}(\R^2)}\nonumber\\
&+\int_0^{\tau} e^{-\underline{\text{\textnu}}(\tau-\tau')+\int_{\tau'}^{\tau}\norm{\nabla u(\tau'')}_{L^\infty(\R^2)}d\tau''}\norm{F(\tau')}_{\dot \cC^\alpha_{pw,\gamma(\tau')}(\R^2)}d\tau'.\label{c1.3}
\end{align}
Given that for all $0 \leqslant s' < s$, the following inequality holds:
\begin{gather}\label{eq0.45}
\int_{s'}^s\norm{\nabla u(\tau'')}_{\cC^\alpha_{pw,\gamma(\tau'')}(\R^2)} d\tau'' \leqslant  \dfrac{1}{2\underline{\text{\textnu}}} (s-s')+ C_* \text{\texttheta}(s), 
\end{gather}
we deduce
\[
\norm{f(\rho(\tau))}_{\dot \cC^\alpha_{pw,\gamma(\tau)}(\R^2)}\leqslant e^{C_* \text{\texttheta}(\tau)}\left[e^{-\tfrac{\underline{\text{\textnu}}}{2}\tau}\norm{f(\rho_0)}_{\dot \cC^\alpha_{pw,\gamma_0}(\R^2)}+ \int_0^{\tau} e^{-\tfrac{\underline{\text{\textnu}}}{2}(\tau-\tau')}\norm{F(\tau')}_{\dot \cC^\alpha_{pw,\gamma(\tau')}(\R^2)}d\tau'\right],
\]
From this, we infer, following the computations leading to \eqref{c1.6} and \eqref{c2.2}, that:
\begin{gather}\label{c1.7}
 \sup_{[0,t]}\norm{f(\rho)}_{\dot \cC^\alpha_{pw,\gamma}(\R^2)}^4+\int_{0}^t \norm{f(\rho)}_{\dot \cC^\alpha_{pw,\gamma}(\R^2)}^4\leqslant C_* e^{C_* \text{\texttheta}(t)}\left(\norm{f(\rho_0)}_{\dot \cC^\alpha_{pw,\gamma_0}(\R^2)}^4
+\int_0^t\sigma^{r_\alpha}\norm{F}_{\dot \cC^\alpha_{pw,\gamma}(\R^2)}^4\right).
\end{gather}
Finally, the estimate for the pressure follows from \eqref{c1.6}, \eqref{c2.2}, and \eqref{c1.7}, and we now turn to the final step  devoted to the velocity gradient.
\paragraph{\textbf{Final estimates}}
We start with the following expression of the velocity gradient derived in \eqref{v.34}:
\[
\nabla u= \nabla u_* + \nabla u_F+ \nabla u_P+ \nabla u_\delta.
\]
The estimates derived for $(-\Delta)^{-1}\dvg (\rho \dot{u})$ in \eqref{v.33}-\eqref{v.31} and for $[K,\mu(\rho)-\widetilde{\mu}]\D u$ in \eqref{v.32} apply to $\nabla u_\delta$ and $\nabla u_*$, respectively. We now turn our focus to the
$L^4((0,t), \cC^\alpha_{pw,\gamma}(\R^2))$-norms estimates of $\nabla u_P$ and $\nabla u_F$. We first recall:
\[
\nabla u_P= -\nabla (-\Delta)^{-1}\nabla \left(\text{\textpsi}_1(\rho)(P(\rho)-\widetilde P)\right),\quad\text{with}\quad \text{\textpsi}_1(\rho)= \dfrac{2\mu(\rho)-\widetilde\mu}{2\mu(\rho)+\lambda(\rho)}.
\]
By applying \cref{ThB2}, we obtain:
\begin{align*}
\norm{\nabla u_P}_{\cC^\alpha_{pw,\gamma}(\R^2)}&\leqslant C\left(\norm{\text{\textpsi}_1(\rho)(P(\rho)-\widetilde P)}_{L^4(\R^2)}+\norm{\text{\textpsi}_1(\rho)(P(\rho)-\widetilde P)}_{\cC^\alpha_{pw,\gamma}(\R^2)}\right)\nonumber\\
&+ C\ell_\vph^{-\tfrac{1}{4}}\norm{\llbracket \text{\textpsi}_1(\rho)(P(\rho)-\widetilde P)\rrbracket}_{L^4(\mathcal{C})}+C\norm{\llbracket \text{\textpsi}_1(\rho)(P(\rho)-\widetilde P)\rrbracket}_{L^\infty(\mathcal{C})}\left(\ell_\vph^{-\alpha}+\mathfrak{P}_{\gamma}\right).
\end{align*}
Next, using \eqref{eq0.24} along with the previous step, we obtain:
\begin{align*}
\int_0^t \sigma^{r_\alpha} &\left[\norm{\text{\textpsi}_1(\rho)(P(\rho)-\widetilde P)}_{L^{4}(\R^2)}^4+\norm{\text{\textpsi}_1(\rho)(P(\rho)-\widetilde P)}_{\cC^\alpha_{pw,\gamma}(\R^2)}^4\right]d\tau \\
&\leqslant C_*\left( C_0+\mathcal{A}_1(t)^2\right)+C_*(1+\text{\texttheta}(t))e^{C_* \text{\texttheta}(t)}\left(\norm{f(\rho_0)}_{\dot \cC^\alpha_{pw,\gamma_0}(\R^2)}^4
+\int_0^t\sigma^{r_\alpha}\norm{F}_{\dot \cC^\alpha_{pw,\gamma}(\R^2)}^4\right)\\
&\leqslant C_*\left[ C_0+\mathcal{A}_1(t)^2+e^{C_* \text{\texttheta}(t)}\left(\norm{f(\rho_0)}_{\cC^\alpha_{pw,\gamma_0}(\R^2)}^4
+\int_0^t\sigma^{r_\alpha}\norm{F}_{\cC^\alpha_{pw,\gamma}(\R^2)}^4\right)\right].
\end{align*}
Following the computations leading to \eqref{v.29}, we have: 
\begin{align*}
\int_0^t\sigma^{r_\alpha}\left[\ell_\vph^{-1}\norm{\llbracket \text{\textpsi}_1(\rho)(P(\rho)-\widetilde P)\rrbracket}_{L^4(\mathcal{C})}^4\right.&\left.+\norm{\llbracket \text{\textpsi}_1(\rho)(P(\rho)-\widetilde P)\rrbracket}_{L^\infty(\mathcal{C})}^4\left(\ell_{\vph}^{-4\alpha}+\mathfrak{P}_{\gamma}^4\right)\right]\\
&\leqslant C_*K_0e^{C_*\text{\texttheta}(t)}\norm{\llbracket f(\rho_0)\rrbracket}_{L^4\cap L^\infty(\mathcal{C}(0))}^4.
\end{align*}
As a result, we obtain:
\begin{multline}\label{c2.8}
\int_0^t\sigma^{r_\alpha}\norm{\nabla u_P(\tau)}_{\cC^\alpha_{pw,\gamma(\tau)}(\R^2)}^4d\tau \\
\leqslant C_*\left[ C_0+\mathcal{A}_1(t)^2+e^{C_* \text{\texttheta}(t)}\left(\norm{f(\rho_0)}_{\cC^\alpha_{pw,\gamma_0}(\R^2)}^4
+K_0\norm{\llbracket f(\rho_0)\rrbracket}_{L^4\cap L^\infty(\mathcal{C}(0))}^4+\int_0^t\sigma^{r_\alpha}\norm{F}_{\cC^\alpha_{pw,\gamma}(\R^2)}^4\right)\right].
\end{multline}
We now proceed to estimate $\nabla u_F$, which is given by:
\[
\nabla u_F= \nabla (-\Delta)^{-1}\nabla \left( \text{\textpsi}_2(\rho) F\right), \quad\text{with}\quad \text{\textpsi}_2(\rho)=\dfrac{\widetilde\mu+\lambda(\rho)}{2\mu(\rho)+\lambda(\rho)}.
\]
Once again, we apply \cref{ThB2} to obtain:
\begin{align*}
\norm{\nabla u_F}_{\cC^\alpha_{pw,\gamma}(\R^2)}&\leqslant C\left(\norm{\text{\textpsi}_2(\rho) F}_{L^4(\R^2)}+\norm{\text{\textpsi}_2(\rho)F}_{\cC^\alpha_{pw,\gamma}(\R^2)}\right)\nonumber\\
&+ C\ell_\vph^{-\tfrac{1}{4}}\norm{\llbracket \text{\textpsi}_2(\rho)F\rrbracket}_{L^4(\mathcal{C})}+C\norm{\llbracket \text{\textpsi}_2(\rho)F\rrbracket}_{L^\infty(\mathcal{C})}\left(\ell_\vph^{-\alpha}+\mathfrak{P}_{\gamma}\right).
\end{align*}
Straightforwardly, we derive:
\begin{gather*}
\int_{0}^t\sigma^{r_\alpha} \left[\norm{\text{\textpsi}_2(\rho) F}_{L^{4}(\R^2)}^{4}+\norm{\text{\textpsi}_2(\rho) F}_{ \cC^\alpha_{pw,\gamma}(\R^2)}^4\right] \leqslant C_*\left[ C_0+ \mathcal{A}_1(t)^2+\left(1+\text{\texttheta}(t)\right)\int_0^t\sigma^{r_\alpha}\norm{F}_{\cC^\alpha_{pw,\gamma}(\R^2)}^4\right],
\end{gather*}
and (see the expression of $\llbracket F\rrbracket$ in \eqref{eq0.31})
\[
\norm{\llbracket \text{\textpsi}_2(\rho)F\rrbracket}_{L^4\cap L^\infty(\mathcal{C})}\displaystyle\leqslant C_*\norm{\llbracket f(\rho)\rrbracket}_{L^4\cap L^\infty (\mathcal{C})}\left(\norm{\nabla u}_{L^\infty(\R^2)}+\norm{F}_{L^\infty(\R^2)}\right).
\]
Following the computations leading to \eqref{v.28}, we have: 
\begin{align*}
    \int_0^t\sigma^{r_\alpha}\left[\ell_\vph^{-1}\norm{\llbracket \text{\textpsi}_2(\rho)F\rrbracket}_{L^4(\mathcal{C})}^4\right.&\left.+\norm{\llbracket \text{\textpsi}_2(\rho)F\rrbracket}_{L^\infty(\mathcal{C})}^4\left(\ell_\vph^{-4\alpha}+\mathfrak{P}_{\gamma}^4\right)\right]\\
    &\leqslant C_*K_0e^{C_*\text{\texttheta}(t)}\norm{\llbracket f(\rho_0)\rrbracket}_{L^4\cap L^\infty(\mathcal{C}(0))}^4\left(\text{\texttheta}(t)+\int_{0}^t\sigma^{r_\alpha}\norm{F}_{L^\infty(\R^2)}^{4}\right),
\end{align*}
 and finally:
\begin{align}\label{c2.7}
    \int_0^t\sigma^{r_\alpha}\norm{\nabla u_F}_{\cC^\alpha_{pw,\gamma}(\R^2)}^4&\leqslant C_* \left( C_0+\mathcal{A}_1(t)^2+K_0
e^{C_*\text{\texttheta}(t)}\norm{\llbracket f(\rho_0) \rrbracket}_{L^4\cap L^\infty(\mathcal{C}(0))}^{4}\right)\nonumber\\
&+C_*\left[1+\text{\texttheta}(t)+K_0
e^{C_*\text{\texttheta}(t)}\norm{\llbracket f(\rho_0) \rrbracket}_{L^4\cap L^\infty(\mathcal{C}(0))}^{4}\right]\int_0^t\sigma^{r_\alpha}\norm{F}_{\cC^\alpha_{pw,\gamma}(\R^2)}^4.
\end{align}
 \cref{prop3} just follows by summing \eqref{c2.5}, \eqref{c2.7} and \eqref{c2.8}.
\enddem
\subsection{\texorpdfstring{Proof of \cref{thglobal}}{}}\label[section]{stability}
This section is devoted to the proof of the main result of this paper. It is structured into two steps: the first part, \cref{existence}, focuses  on the construction of a solution $(\rho,u)$ to the Navier-Stokes equations \eqref{ep4.1}, while the second part, \cref{uniqueness}, establishes uniqueness within a large space.
\subsubsection{Proof of the existence}\label[section]{existence}
This section is dedicated to constructing a solution for the Navier-Stokes equations \eqref{ep4.1}. It starts with the construction of an approximate sequence $(\rho^\delta,u^\delta)$ and goes up to the convergence of this sequence to a limit $(\rho,u)$ that solves the equations \eqref{ep4.1}. Usually, $(\rho^\delta, u^\delta)$ corresponds to the solution to the Cauchy problem \eqref{epq4} with initial data $(\rho_0^\delta, u_0^\delta)$ where $(\rho_0^\delta,u_0^\delta)$ is obtained by smoothing $(\rho_0,u_0)$. In our case, this does not seem a good idea, as smoothing the initial data would result in the loss of the density discontinuity. This motivates our local-in-time well-posedness result in \cite{zodji2023well}. Although a compatibility condition is required, the solution exhibits similar regularity to that of \cref{thglobal}. In particular, the density and velocity gradient are discontinuous. We now begin the proof of existence by constructing the sequence of initial data $(\rho_0^\delta, u^\delta_0)$.

\paragraph{\textbf{Step 1: Construction of $(\rho_0^\delta, u^\delta_0)$}}

We initiate by identifying which initial quantity needs to be smoothed and which not. The initial interface $\gamma_0$, density and velocity, possess the necessary regularity as outlined in the local-in-time theorem \cref{thlocal}. Therefore, there is no need to smooth these quantities. However, the stress tensor at the initial time only  fulfills:
\[
\dvg(\Pi_0)=\dvg (2\mu(\rho_0)\D u_0+ (\lambda(\rho_0)\dvg u_0-P(\rho_0)+\widetilde P)I_2)\in H^{-1}(\R^2).
\]
To preserve the discontinuity in the initial density, we consider $(\rho_0,u_0^\delta)$ (namely $\rho_0^\delta=\rho_0)$ 
as a sequence of initial data, where $u_0^\delta$ solves  the following elliptic equation:
\begin{gather}\label{epq31}
-\dvg (2\mu(\rho_0)\D u_0^\delta +(\lambda(\rho_0)\dvg u_0^\delta- P(\rho_0)+\widetilde P)I_2)+c^\delta u_0^\delta= -\dvg (w_\delta*\Pi_0).
\end{gather}
Above, $w_\delta:=\delta^{-2}w(\cdot/\delta)$, $\delta\in (0,1)$, where  $w$ is a smooth non negative function supported in the unit ball centered at the origin, and whose integral equals $1$. Additionally, the constant $c^\delta$ is defined as:
\[
c^\delta:=\norm{w_\delta*\Pi_0-\Pi_0}_{L^2(\R^2)}\quad\text{and satisfies}\quad c^\delta\xrightarrow[]{\delta\to 0} 0.
\]
Since the viscosity $2\mu(\rho_0)+\lambda(\rho_0)$ is bounded away from vacuum and from above, and the pressure $P(\rho_0) - \widetilde P$ belongs to $L^2(\R^2)$, the existence of a unique solution $u_0^\delta \in H^1(\R^2)$ of \eqref{epq31} follows from the Lax-Milgram theorem. Moreover, the sequence $(u_0^\delta)_\delta$ satisfies the following estimate:
\begin{gather}\label{epq32}
C_*^{-1} \norm{\nabla u_0^\delta}_{L^2(\R^2)}^2+ c^\delta \norm{u_0^\delta}_{L^2(\R^2)}^2\leqslant C_* \left(\norm{\Pi_0}_{L^2(\R^2)}^2+\norm{P(\rho_0)-\widetilde P}_{L^2(\R^2)}^2\right) .
\end{gather}
We now move on to proving that $(u^\delta_0)_\delta$ converges strongly to $u_0$ in $H^1(\R^2)$.

We add  $\dvg (\Pi_0)$ to both sides of  \eqref{epq31},  obtaining:
\[
-\dvg\{2\mu(\rho_0)\D (u^\delta_0-u_0)+ \lambda(\rho_0)\dvg (u_0^\delta-u_0) I_2\}+ c^\delta u_0^\delta=-\dvg (w_\delta*\Pi_0-\Pi_0).
\]
Using  $u_0^\delta -u_0$ as a test function yields:
\begin{gather*}
2\int_{\R^2}\mu(\rho_0)\abs{\D (u_0^\delta -u_0)}^2+\int_{\R^2}\lambda(\rho_0)\abs{\dvg (u_0^\delta-u_0)}^2+c^\delta \int_{\R^2} u_0^\delta (u_0^\delta -u_0)=\int_{\R^2}\nabla (u_0^\delta-u_0)\colon (w_\delta*\Pi_0- \Pi_0),
\end{gather*}
from which,  with the help of H\"older's and Young's inequalities, we deduce:
\begin{gather}\label{epq33}
C_*^{-1}\norm{\nabla (u_0^\delta-u_0)}_{L^2(\R^2)}^2+\dfrac{1}{2}c^\delta \norm{u_0^\delta}_{L^2(\R^2)}^2\leqslant C_*\norm{w_\delta*\Pi_0-\Pi_0}_{L^2(\R^2)}^2+\dfrac{1}{2}c^\delta\norm{u_0}_{L^2(\R^2)}^2.
\end{gather}
Now, given that $c^\delta\xrightarrow[]{\delta \to 0} 0$ (see  \eqref{epq32} above), we immediately obtain:
\[
u_0^\delta \quad\xrightarrow[]{\delta \to 0} \quad u_0 \quad\text{ in }\quad \dot H^1(\R^2).
\]
Furthermore, \eqref{epq33} implies that:
\[
\limsup_{\delta\to 0}\norm{u_0^\delta}_{L^2(\R^2)}^2\leqslant  \norm{u_0}_{L^2(\R^2)}^2,
\]
and whence:
\[
u_0^\delta \quad\xrightarrow[]{\delta \to 0} \quad u_0 \quad\text{ in }\quad L^2(\R^2).
\]
This proves the strong convergence of  $(u_0^\delta)$ to $u_0$ in $H^1(\R^2)$. Additionally, as intended, we obtain from \eqref{epq31}:
\begin{gather}\label{epq34}
\dvg (2\mu(\rho_0)\D u_0^\delta +(\lambda(\rho_0)\dvg u_0^\delta- P(\rho_0)+\widetilde P)I_2)=-c^\delta u_0^\delta +\dvg (w_\delta*\Pi_0)\in L^2(\R^2).
\end{gather}

Finally, considering the regularity of the initial density and interface $\gamma_0$, the small perturbation assumption of the initial viscosity $\mu(\rho_0)$ around the constant state $\widetilde{\mu}$, along with the regularity of $u_0^\delta$ and \eqref{epq34}, \cref{thlocal} ensures the existence of a unique solution $(\rho^\delta, u^\delta)$ for equations \eqref{epq4} with the initial conditions: 
\[
\rho^\delta_{|t=0}=\rho_0\quad \text{ and }\quad u^\delta_{|t=0}= u^\delta_0.
\]
The solution is defined up to a maximal time $T_\delta>0$ and enjoys the regularity in \cref{thlocal}. This regularity is sufficient for the computations carried out in the previous sections to make sense. As a consequence, \cref{lem4.9} holds true for solution $(\rho^\delta,u^\delta)$, as well as the first condition of the blow-up criterion, \eqref{v.40}. For the second condition, \eqref{c3.38}, we use the exponential decay in time for jumps to derive: for all $t\in (0,T_\delta)$:
\begin{align*}
    &\left[1+\norm{\lambda (\rho(t))}_{\dot \cC^\alpha_{pw,\gamma(t)}(\R^2)}+\left( \mathfrak{P}_{\gamma(t)} +\ell^{-\alpha}_{\vph(t)}\right)\llbracket\lambda(\rho(t))\rrbracket_{L^\infty(\mathcal{C}(t))}\right]\norm{\mu(\rho(t))-\widetilde\mu}_{\cC^\alpha_{pw,\gamma(t)}(\R^2)}\\
    +&\left( \mathfrak{P}_{\gamma(t)} +\ell^{-\alpha}_{\vph(t)}\right)\left[\norm{\llbracket\mu(\rho(t))\rrbracket}_{L^\infty(\mathcal{C}(t))}+\norm{\llbracket\mu(\rho(t))\rrbracket,\, \llbracket \lambda(\rho(t))\rrbracket}_{L^\infty(\mathcal{C}(t))}\left\|1-\dfrac{\widetilde\mu}{<\mu(\rho(t))>}\right\|_{L^\infty(\mathcal{C}(t))}\right]\\
    &\leqslant C_*\left[\left(1+\text{\texttheta}(t)^{\tfrac{1}{4}}+ K_0 e^{C_* \text{\texttheta}(t)}\llbracket f(\rho_0)\rrbracket_{L^\infty(\mathcal{C}(0))}\right)\text{\texttheta}(t)^{\tfrac{1}{4}}+K_0 e^{C_* \text{\texttheta}(t)}\llbracket f(\rho_0)\rrbracket_{L^\infty(\mathcal{C}(0))}\right]. 
\end{align*}
Due to the smallness of $c_0$ (see \eqref{ep3.1} above), \eqref{c3.38} is satisfied, leading to $T_\delta=+\infty$. 
We now proceed to the final step, which focuses on showing that the sequence $(\rho^\delta, u^\delta)$ converges to a pair 
$(\rho,u)$ that solves \eqref{ep4.1}.
\paragraph{\textbf{Step 2: Convergence of the approximate sequence $(\rho^\delta, u^\delta)$}}
We recall that for all $\delta>0$, the pair $(\rho^\delta, u^\delta)$ satisfies the following system:
\begin{gather}\label{c2.16}
\begin{cases}
    \dpt \rho^\delta +\dvg (\rho^\delta u^\delta)=0,\\
    \dpt (\rho^\delta u^\delta)+\dvg(\rho^\delta u^\delta \otimes u^\delta)+\nabla P(\rho^\delta)=\dvg (2\mu(\rho^\delta)\D u^\delta)+\nabla (\lambda(\rho^\delta)\dvg u^\delta).
\end{cases}
\end{gather}
Additionally,  for all $T\in (0,\infty)$ and $\delta\in (0,1)$, we have:
\begin{align}
\norm{\rho^\delta-\widetilde\rho}_{L^\infty((0,T), L^2(\R^2))}^2&+\norm{\rho^\delta-\widetilde\rho}_{L^\infty((0,T)\times\R^2)}^2+\norm{u^\delta}_{L^\infty((0,T), H^1(\R^2))}^2\nonumber\\
&+\norm{\nabla u^\delta}_{L^3((0,T)\times \R^2)}^2+ \norm{\dot u^\delta}_{L^2((0,T)\times \R^2)}^2\leqslant C_{*,0}.\label{epq35}
\end{align}
Hereafter $C_{*,0}$ is a positive constant that depends on $C_*$ and $c_0$. Sometimes we write $C_{*,0}(T)$ (resp. $C_{*,0}(n)$) to emphasize the additional dependence on $T>0$ (resp. $n\in \N^*$).
 From \eqref{epq35}, there exist $\rho -\widetilde\rho\in L^\infty((0,\infty), L^2(\R^2))\cap L^\infty((0,\infty)\times\R^2)$, and $u\in L^\infty((0,\infty), H^1(\R^2))$ such that: 
\begin{gather}
    \begin{cases}
        \rho^\delta -\widetilde\rho\;\; &\rightharpoonup^*\;\;  \rho -\widetilde\rho \quad \text{ in }\quad L^\infty((0,T)\times \R^2),\\
        \rho^\delta-\widetilde\rho \;\;&\xrightarrow{}\quad\rho -\widetilde\rho \quad\text{ strongly in }\quad \cC([0,T], L^2_w(\R^2)),\\
        u^\delta  \;\; &\rightharpoonup^*\;\; u \quad \text{ in } L^\infty((0,T), H^{1}(\R^2)),\\
        u^\delta  &\xrightarrow{}\quad u \quad\text{ strongly in } \cC([0,T], L^2_\loc(\R^2)).
    \end{cases}
\end{gather}
Additionally, by interpolation, we have: 
\[
u^\delta \xrightarrow{} u \quad \text{ strongly in } \quad L^\infty_\loc((0,\infty),L^p_\loc(\R^2)),\quad \text{for every}\quad p\in [2,\infty).
\]
The initial interface $\mathcal{C}(0)$ is transported by the flow $X^\delta$ associated with the velocity $u^\delta$ into an interface 
\begin{gather}\label{c2.15}
\mathcal{C}^\delta(t)=X^\delta(t) \mathcal C(0) \quad\text{ with parameterization}\quad \gamma^\delta (t,s)=X^\delta(t,\gamma_0(s)).
\end{gather}
Given that the velocity sequence $(u^\delta)$ satisfies:
\begin{gather}\label{epq37}
    \sup_{(0,\infty)} \sigma \norm{\nabla u^\delta }_{L^4(\R^2)}^4+\int_0^\infty\sigma^{r_\alpha}\norm{\nabla u^\delta }_{\cC^\alpha_{pw,\gamma^\delta}(\R^2)}^4\leqslant C_{*,0},
\end{gather}
 we infer, together with \eqref{epq35}, that  for any $T>0$ and $n\in \N\setminus\{0\}$, 
\[
\sup_{[0,T]}\norm{\nabla \gamma^\delta}_{L^\infty}\leqslant C_{*,0}(T),\quad\text{and}\quad \sup_{[1/n,T]}\norm{\dpt \gamma^\delta}_{L^\infty}\leqslant C_{*,0}(n), \quad \text{for all} \quad \delta>0.
\] 
Hence, up to a subsequence, $(\gamma^\delta)_\delta$ convergences uniformly on $[0,T]\times V$ to some $\gamma\in W^{1,\infty}((0,T)\times V)$. From \eqref{epq35} and the embedding $W^{1,3}(\R^2)\hookrightarrow \cC^{1/3}(\R^2)$, the velocity sequence satisfies:
\[
\sup_{\delta} \norm{u^\delta}_{L^3((0,\infty), \cC^{1/3}(\R^2))}\leqslant C_{*,0}.
\]
Therefore, taking the limit as $\delta \to 0$ in \eqref{c2.15}, we find that the limit parameterization $\gamma$ satisfies:
\begin{gather}\label{epq36}
\gamma(t,s)=\gamma_0(s)+ \int_0^t u(\tau, \gamma(\tau,s))d\tau.
\end{gather}
Furthermore, from \eqref{c2.15}-\eqref{epq37}, we obtain:
\[
\sup_{\delta}\int_0^T\sigma^{r_\alpha}\norm{\dpt \nabla \gamma^\delta}_{\cC^\alpha}^4\leqslant C_{*,0}(T).
\]
Since $r_\alpha<3$, there exists a  $q\in (1, 4/(1+r_\alpha))$ such that:
\[
\sup_{\delta}\int_0^T\norm{\dpt \nabla \gamma^\delta}_{\cC^\alpha}^q\leqslant  C_{*,0}(T,q).
\]
Therefore, the lower semi-continuity of norms implies that $\dpt \nabla \gamma \in L^q_{\loc}([0,\infty), \cC^{\alpha})$, 
and as a result,
\[
\gamma \in \cC([0,\infty), \cC^{1+\alpha}(V)).
\]
Once more, \eqref{epq37} combined with the lower semi-continuity of norms implies that $\nabla u\in L^q_\loc ([0,\infty), L^\infty(\R^2))$. This guarantees the uniqueness of $\gamma$ that satisfies \eqref{epq36}. 
We now move on to the proof of the strong convergence of the density sequence  $(\rho^\delta)_\delta$.

Arguing as above, we obtain that the sequence $(X^\delta)_\delta$ converges uniformly on compact sets  in $[0,\infty)\times \R^2$ to $X\in \cC([0,\infty), \cC^{1+\alpha}(\R^2))$, the flow of the limit velocity $u$. We observe that the a priori estimate yields the following bound:
\[
\sup_{\delta}\norm{\rho^\delta-\widetilde\rho}_{L^\infty((0,\infty), \cC^\alpha_{pw,\gamma^\delta}(\R^2))}^4\leqslant C_{*,0}.
\]
From this, combined with \eqref{epq37}, we infer that for any $T>0$:
\[
\sup_{\delta} \norm{\varrho^\delta-\widetilde\rho}_{L^\infty((0,T), \cC^\alpha_{pw,\gamma_0}(\R^2))}^4+\sup_{\delta}\norm{\dpt \varrho^\delta}_{L^q ((0,T), L^\infty(\R^2))}^4\leqslant C_{*,0}(T,q),
\]
where $\varrho^\delta$ is given by:
\[
\varrho^\delta(t,x)=\rho^\delta(t,X^\delta(t,x)).
\]
Since $\cC^\alpha_{pw,\gamma_0}(\R^2)$ embeds compactly in $L^\infty_\loc (\R^2)$, we deduce 
from Aubin-Lions Lemma that  $(\varrho^\delta)_\delta$ converges uniformly on compact sets in $\{(t,x)\colon x\in \R^2\setminus \mathcal{C}(0)\}$  to some $\varrho \in\cC([0,\infty), \cC^\alpha_{pw,\gamma_0}(\R^2))$. The final step is to prove that 
$(\rho^\delta)_\delta$ converges strongly to $\overline{\varrho}$  in $L^2_\loc([0,\infty)\times \R^2)$, where 
\[
\overline{\varrho}(t,x)=\varrho(t, X^{-1}(t,x)).
\]
 Above, $X^{-1}\in \cC([0,\infty), \cC^{1+\alpha}(\R^2))$  satisfies:
\begin{gather}\label{epq38}
X^{-1}(t,x)= x-\int_0^t u(\tau,X^{-1}(\tau,x))d\tau,\quad \;X(t,X^{-1}(t,x))=x\quad\text{ and }\quad X^{-1}(t,X(t,x))=x.
\end{gather}
Let $T>0$ and $B$ be an arbitrary bounded subset of $\R^2$. We have:
\begin{multline}\label{c3.10}
\int_0^T\int_{B}\abs{\rho^\delta(t,x)-\overline{\varrho}(t,x)}^2dtdx=\int_0^T\int_{X^\delta(t)B}\abs{\rho^\delta(t,X^\delta(t,x))-\overline{\varrho}(t,X^\delta(t,x))}^2 J^\delta(t,x)dtdx\\
\leqslant C_* \int_0^T\int_{X^\delta(t)B} \abs{\varrho^\delta(t,x)-\varrho(t,x)}^2dtdx+C_*\int_0^T\int_{X^\delta(t) B}\abs{\overline{\varrho}(t,X^\delta(t,x))-\varrho(t,x)}^2dtdx,
\end{multline}
where $J^{\delta}$ comes from the change of variable $x\mapsto X^\delta(t,x)$ and satisfies (see \cite[Lemma 3.2]{hoff1995global}):
\[
\sup_{[0,T]}\norm{J^{\delta}}_{L^\infty(\R^2)}\leqslant C_*.
\]
\begin{itemize}
    \item Given that  $(X^\delta)_\delta$ converges to $X$ uniformly on compact sets in $[0,\infty)\times \R^2$ , there exists $\delta_0>0$ such that for all $\delta<\delta_0$, we have:
\[
X^\delta(t) B\subset X(t) B+ B(0,1).
\]
Furthermore, because  $\varrho$ is uniformly continuous on both sides of $\mathcal{C}(0)$, and using \eqref{epq38}, we obtain:
\[
\overline{\varrho}(t,X^\delta(t,x))=\varrho (t,X^{-1}(t,X^\delta(t,x)))\quad\xrightarrow{\delta \to 0}\quad\varrho(t,x) \quad\text{a.e}\quad x\in \R^2.
\]
Therefore, we conclude:
\begin{gather}\label{epq39}
    \lim_{\delta\to 0}\int_0^T\int_{X^\delta(t) B}\abs{\overline{\varrho}(t,X^\delta(t,x))-\varrho(t,x)}^2dtdx=0.
\end{gather}
\item  Additionally, since the sequence $(\varrho^\delta)_\delta$ converges to $\varrho$ uniformly on compact sets, we have: 
\begin{gather}\label{epq40}
\lim_{\delta\to 0}\int_0^T\int_{X^\delta(t)B} \abs{\varrho^\delta(t,x)-\varrho(t,x)}^2dtdx=0.
\end{gather}
\end{itemize}
Finally, the strong convergence of $(\rho^\delta)_\delta$ to $\overline\varrho$ in $L^2_\loc ((0,\infty)\times\R^2)$ follows from \eqref{c3.10}, \eqref{epq39}, and \eqref{epq40}. In particular, we have $\rho=\overline{\varrho}\in \cC([0,\infty),\cC^\alpha_{pw,\gamma}(\R^2))$. Additionally, by interpolation, it follows that  $(\rho^\delta)_\delta$ converges to $\rho$ strongly in $L^p_\loc ((0,\infty)\times \R^2)$ for any $1\leqslant p<\infty$.

By combining the strong and weak convergences of  $(\rho^\delta,u^\delta)_\delta$ with standard arguments, we take the limit as $\delta\to 0$  in \eqref{c2.16}  and we conclude that $(\rho,u)$ satisfies \eqref{ep4.1}. This completes the current section. We  now turn to the proof of uniqueness.
\subsubsection{Proof of the uniqueness}\label[section]{uniqueness}
The uniqueness of the solution constructed above is immediately implied by the following result.
\begin{prop}\label[prop]{propv1}
   Consider the system \eqref{ep4.1} and assume that the pressure and viscosity laws are $W^{1,\infty}$-regular functions of the density. To simplify, we assume that $\lambda$ is  nonnegative on $[0,\infty)$. 
   Let $(\rho_0, u_0)$ be the initial data associated with \eqref{ep4.1}, satisfying the following conditions:
    \begin{gather}\label{v.41}
         \rho_0,\;\dfrac{1}{\rho_0}\in L^\infty(\R^2),\quad \text{and}\quad u_0\in L^2(\R^2).
    \end{gather}
    
    Let $T>0$. On the time interval $[0,T]$, there exists at most one solution to the Cauchy problem associated with \eqref{ep4.1} and initial data $(\rho_0, u_0)$ satisfying:
    \begin{gather}\label{v.42}
            \dfrac{1}{\mu(\rho)}\in L^{\infty}((0,T)\times\R^2),\;\,\nabla u\in L^1((0,T), L^\infty(\R^2), \;\, \text{and}\;\,\sqrt{\sigma}\nabla u\in L^2((0,T), L^\infty(\R^2)).
    \end{gather}
\end{prop}
\cref{propv1} establishes uniqueness for the system \eqref{ep4.1} within a broader framework than that of  \cref{thglobal}. In particular, neither piecewise H\"older continuity for the density or velocity gradient, nor smallness conditions on the initial data or viscosity fluctuations are required.
\dem[Proof of \cref{propv1}]
Let $(\rho,u)$ and $(\varrho,v)$ be two solutions to the Cauchy problem associated with \eqref{ep4.1} and with initial data $(\rho_0,u_0)$ satisfying \eqref{v.41}. Additionally, we assume that  $(\rho,u)$ and $(\varrho,v)$ satisfy the conditions in \eqref{v.42}. As a consequence, for any  $k_0\in (0,1)$, there exists a time $T_0>0$ such that  
\begin{gather}\label{v.38}
\int_0^{T_0} \norm{\nabla u}_{L^\infty(\R^2)}<k_0\quad \text{and}\quad  \int_0^{T_0} \norm{\nabla v}_{L^\infty(\R^2)}<k_0.
\end{gather}
The regularity of $u,v$ is sufficient to recast the equations they satisfy in Lagrangian coordinates:
\begin{gather}\label{v.39}
    \begin{cases}
    \rho_0\dpt \overline u-\dvg \left[\adj(D X_{u})\left(2\mu(\rho_0 J_u^{-1})\D_{A_{u}} \overline u+\left(\lambda(\rho_0 J_u^{-1})\dvg_{A_{u}} \overline u-P(\rho_0 J_u^{-1})+\widetilde P\right) I\right)\right]=0,\\
    \rho_0\dpt \overline v-\dvg \left[\adj(D X_{v})\left(2\mu(\rho_0 J_v^{-1})\D_{A_{v}} \overline v+\left(\lambda(\rho_0 J_v^{-1})\dvg_{A_{v}} \overline v-P(\rho_0 J_v^{-1})+\widetilde P\right)I\right)\right]=0.
    \end{cases}
\end{gather}
Here, $X_w$ is the flow associated with the velocity $w$. We define $\overline{w}$ by:
\[
\overline{w}(t,y)= w(t,X_w(t,y)),\;\, \text{so that}\;\, X_w(t,y)= y+\int_0^t\overline{w}(\tau,y)d\tau.
\]
The Jacobian matrix of $X_w$ is denoted by $D X_w$, and we define $J_w= \det (D X_w)$. By \eqref{v.42}, the matrix $D X_w$
is invertible, with its inverse denoted by  $A_w$. The matrix of cofactors of $D X_w$, also known as the adjugate matrix, is denoted by $\adj( D X_w)$. Finally, the operators $\D_{A_{v}}$ and $\dvg_{A_w}$ are defined as follows:
\[
\D_{A_w}z= \dfrac{1}{2}\left(D z\cdot A_w+A_w^T\cdot \nabla z\right),\;\,\text{and}\;\, \dvg_{A_w} z = D z : A_{w} = A^T_{w}:\nabla z.
\]
The computations leading to \eqref{v.39} are standard and can be found, for instance, in \cite{danchinLagrangianapproach}. Additionally, with a slight modification, the following bounds can be derived from Lemmas A.3 and A.4 of the same reference.
\begin{lemm}\label[lemma]{lem5}
There exists a constant $C_{k_0}$, depending only on $k_0$, such that the following estimates hold for all $p\in [1,\infty]$, $t\in [0,T_0]$, and $w\in \{u,v\}$:
\[
  \begin{cases}
   \norm{\adj( D X_w(t)) \D_{A_w(t)} z-\D z }_{L^p(\R^2)} &\leqslant C_{k_0} \norm{\nabla w}_{L^1((0,t), L^\infty(\R^d))}\norm{D z}_{L^p(\R^2)},\\
   \norm{\adj(D X_w(t))\dvg_{A_w(t)} z-\dvg z I_d}_{L^p(\R^2)}&\leqslant C_{k_0}\norm{\nabla w}_{L^1((0,t), L^\infty(\R^2))}\norm{D z}_{L^p(\R^2)},
   \end{cases}
   \]
   and
\begin{gather*}
\begin{cases}
    \norm{A_{u}(t)-A_{v}(t)}_{L^p(\R^2)}&\leqslant C_{k_0} \norm{\nabla \delta\overline{u}}_{L^1((0,t),L^p(\R^2))},\\
    \norm{\adj(D X_{u}(t))-\adj(D X_{v}(t))}_{L^p(\R^2)}&\leqslant C_{k_0} \norm{\nabla \delta\overline{ u}}_{L^1((0,t),L^p(\R^2))},\\
    \norm{J^{\pm 1}_{u}(t)-J^{\pm 1}_{v}(t)}_{L^p(\R^2)}&\leqslant C_{k_0} \norm{\nabla \delta \overline u}_{L^1((0,t),L^p(\R^2))},
\end{cases}     
    \end{gather*}
    where $\delta\overline{u}=\overline{u}-\overline{v}$.
\end{lemm}
We now take the difference in \eqref{v.39} and obtain:
\begin{gather}\label{v.43}
    \rho_0\dpt \delta\overline u-\dvg \left(2\mu(\rho_0 J^{-1}_{u})\D \delta\overline{u}\right)-\nabla \left(\lambda(\rho_0J^{-1}_{u}) \dvg \delta \overline{u}\right)= \dvg(\mathcal{I}_1)+ \dvg(\mathcal{I}_2)+\dvg (\mathcal{I}_3),
\end{gather}
where 
\[
\begin{cases}
\mathcal{I}_1&=\left(\mu(\rho_0 J_{u}^{-1})-\mu(\rho_0 J_{v}^{-1})\right)  \left(\adj(D X_{u})\D_{A_{u}} \overline u-\D \overline{u}\right)+\mu(\rho_0 J_{v}^{-1})\left(\adj(D X_{u})\D_{A_{u}} \delta \overline u-\D \delta \overline{u}\right)\\
&+\mu(\rho_0 J_{v}^{-1}) \left(\adj(D X_{u})-\adj(D X_{v})\right) \D_{A_{u}}\overline{v}+\mu(\rho_0 J_{v}^{-1}) \adj(D X_{v})\left( \D_{A_{u}}- \D_{A_{v}}\right)\overline{v};\\
\mathcal{I}_2&=\left(\lambda(\rho_0 J_{u}^{-1})-\lambda(\rho_0 J_{v}^{-1})\right)  \left(\adj(D X_{u})\dvg_{A_{u}} \overline u-\dvg \overline u I \right)+\lambda(\rho_0 J_{v}^{-1})\left(\adj(D X_{u})\dvg_{A_{u}} \delta \overline u-\dvg \delta \overline{u} I\right)\\
&+\lambda(\rho_0 J_{v}^{-1}) \left(\adj(D X_{u})-\adj(D X_{v})\right) \dvg_{A_{u}}\overline{v}+\lambda(\rho_0 J_{v}^{-1}) \adj(D X_{v})\left( \dvg_{A_{u}}- \dvg_{A_{v}}\right)\overline{v};\\
\mathcal{I}_3&= (P(\rho_0 J_{v}^{-1})-\widetilde P)\left(\adj(D X_{v})-\adj(D X_{u})\right)+ \adj(D X_{u}) \left(P(\rho_0 J_{v}^{-1})-P(\rho_0 J_{u}^{-1})\right).
\end{cases}
\]

We fix  $k_0\in (0,1)$  and denote by $C_{k_0}^*$ a constant that may depend on $C_{k_0}$ (see \cref{lem5} above), as well as on the lower and upper bounds of the density and viscosity (see \eqref{v.41}-\eqref{v.42} above). Note that this constant may change from one line to the next. We now perform energy estimate for \eqref{v.43}: we use $ \delta \overline{u}$ as a test function and it follows (recall $\delta\overline{u}_{|t=0}=0$):
\begin{align}
\norm{\sqrt{\rho_0}\delta \overline{u}(t)}_{L^2(\R^2)}^2+\int_0^t\int_{\R^2}\left[2\mu(\rho_0 J_u^{-1}) \abs{\D (\delta\overline{u})}^2\right.&\left.+\lambda(\rho_0 J_u^{-1}) (\dvg (\delta \overline u))^2\right]\nonumber\\
&\leqslant \int_0^t\norm{\mathcal{I}_1,\,\mathcal{I}_2,\,\mathcal{I}_3}_{L^2(\R^2)}\norm{\nabla \delta \overline{u}}_{L^2(\R^2)}.
\end{align}
Next, applying \cref{lem5}, we derive the following estimates:
\[
\begin{cases}
    \norm{\mathcal{I}_1(t),\,\mathcal{I}_2(t)}_{L^2(\R^2)}&\leqslant C_{k_0}^*\left( \norm{\nabla\delta\overline{u}}_{L^1((0,t),L^2(\R^2))} \norm{\nabla u(t)}_{L^\infty(\R^2)}+ \norm{\nabla\delta\overline{u}(t)}_{L^2(\R^2))} \right),\\
    \norm{\mathcal{I}_3(t)}_{L^2(\R^2)}&\leqslant C_{k_0}^*  \norm{\nabla\delta\overline{u}}_{L^1((0,t),L^2(\R^2))},
\end{cases}
\]
which leads to:
\[
\norm{\sqrt{\rho_0}\delta \overline{u}(t)}_{L^2(\R^2)}^2+\int_0^t\norm{\nabla(\delta\overline{u})}_{L^2(\R^2)}^2\leqslant  C_{k_0}^*\int_0^t\left[1+\tau\left(1+\norm{\nabla u(\tau)}_{L^\infty(\R^2)}^2\right)\right]\left[\int_0^\tau\norm{\nabla \delta \overline u}_{L^2(\R^2)}^2\right]d\tau.
\]
As a result, we have:
\[
\mathscr{E}(t)\leqslant \int_0^t \left[1+\tau\left(1+\norm{\nabla u(\tau)}_{L^\infty(\R^2)}^2\right)\right] \mathscr{E}(\tau)d\tau
\quad\text{where} 
\quad
\mathscr{E}(t)= \norm{\sqrt{\rho_0}\delta \overline{u}(t)}_{L^2(\R^2)}^2+\int_0^t\norm{\nabla(\delta\overline{u})}_{L^2(\R^2)}^2.
\]
Invoking the assumption \eqref{v.42} and applying Gronwall's Lemma, we conclude that $\mathscr{E}\equiv 0$ on $[0,T_0]$, which ensures uniqueness on $[0,T_0]$. The uniqueness on $[0,T]$ follows by a standard continuation argument, thereby proving \cref{propv1} and, ultimately, \cref{thglobal}.
\enddem
\section*{Acknowledgment}
This project has received funding
from the European Union’s Horizon 2020 research and innovation
program under the Marie Skłodowska-Curie grant agreement No 945332.
I am grateful for the support of the SingFlows project grant (ANR-18- CE40-0027) of the French National Research Agency (ANR).
This work has been partially supported by the project CRISIS (ANR-20-CE40-0020-01), operated by the French National Research Agency (ANR).  I would like to acknowledge my PhD advisors Cosmin Burtea and David G\'erard-Varet for the fruitful discussions and  careful reading of this work. 


\appendix
\section{Energy computations}
In this section, we will provide details of the computations of two estimates for solution $(\rho,u)$ of the system: 
\begin{gather}\label{eqA.1}
\begin{cases}
    \dpt \rho +\dvg(\rho u)=0,\\
    \dpt(\rho u)+\dvg(\rho u\otimes u)=\dvg \Pi,\\
    \llbracket\Pi\rrbracket\cdot n=0\quad \text{ on }\quad  \mathcal{C},
\end{cases}
\end{gather}
where the stress tensor $\Pi$ is given by:
\[
\Pi= 2\mu(\rho)\D u+ (\lambda(\rho)\dvg u-P(\rho)+ \widetilde P) I_d.
\]
\subsection{Second Hoff energy}\label{Hoff2}
We  start  by investigating the Rankine Hugoniot conditions for $\dot u$ at the surface of discontinuity $\mathcal{C}$. We first notice that after applying the 
derivative $\mathcal{A}$ to the momentum equation, we obtain that $\dot u$ solves the equation 
\begin{gather}\label{eqA.15}
\dpt (\rho \dot u^j)+\dvg(\rho \dot u^j u)= \partial_k (\dot \Pi^{jk})+\partial_k(\Pi^{jk}\dvg u)-\dvg(\partial_k u \Pi^{jk}).
\end{gather}
So by Rankine Hugoniot conditions, 
\[
\llbracket\rho \dot u^j\rrbracket n_t+\llbracket \rho \dot u^j u^k\rrbracket n_x^k=\llbracket \dot \Pi^{jk}\rrbracket n_x^k+ \llbracket \Pi^{jk}\dvg u\rrbracket n_x^k-\llbracket \partial_k u^l \Pi^{jk}\rrbracket n_x^l.
\]
Also, thanks to the Rankine Hugoniot condition applied, this time,  to the mass equation $\eqref{eqA.1}_1$ the following jump condition holds true: 
\[
\llbracket\rho\rrbracket n_t+\llbracket \rho u^k\rrbracket n_x^k=0
\]
and since the material derivative of the velocity is continuous, we finally obtain that: 
\begin{gather}\label{eqA.19}
    \llbracket \dot \Pi^{jk}\rrbracket n_x^k+ \llbracket \Pi^{jk}\dvg u\rrbracket n_x^k-\llbracket \partial_k u^l \Pi^{jk}\rrbracket n_x^l=0.
\end{gather}
This relation will be used in the subsequent computations. We recall that the second Hoff estimate consists in  multiplying \eqref{eqA.15} by the material derivative of the velocity before integrating in space. By doing so, we have:
\begin{gather}\label{eqA.16}
    \int_{\R^2}\dot u^j\{\dpt (\rho \dot u^j)+\dvg(\rho \dot u^j u)\}
    =\int_{\R^2}\dot u^j\{\partial_k (\dot \Pi^{jk})+\partial_k(\Pi^{jk}\dvg u)-\dvg(\partial_k u \Pi^{jk})\}.
\end{gather}
The right-hand side of the above equality, is:
\begin{align}
    \int_{\R^2}\dot u^j\{\dpt (\rho \dot u^j)+\dvg(\rho \dot u^j u)\}&=\int_{\R^2}\dpt (\rho \abs{\dot u^j}^2)
    -\int_{\R^2}(\rho \dot u^j\dpt \dot u^j)+\int_{\R^2}\dvg(\rho \abs{\dot u^j}^2 u)(s,x)-\int_{\R^2}\rho \dot u^j u\cdot\nabla \dot u^j \nonumber\\
    &=\dfrac{1}{2}\int_{\R^2}\dpt (\rho \abs{\dot u^j}^2)+\dfrac{1}{2}\int_{\R^2}\abs{\dot u^j}^2\dpt \rho+\dfrac{1}{2}\int_{\R^2}\dvg(\rho \abs{\dot u^j}^2 u)+\dfrac{1}{2}\int_{\R^2}\abs{\dot u^j}^2\dvg(\rho u)\nonumber\\
    &=\dfrac{1}{2}\dfrac{d}{dt}\int_{\R^2}\rho \abs{\dot u^j}^2,\label{eqA.17}
\end{align}
where we have used the Liouville transport equation and the mass equation $\eqref{eqA.1}_1$. We turn now to the computations of the right-hand side of \eqref{eqA.16}:
\begin{align*}
    \int_{\R^2}\dot u^j\{\partial_k (\dot \Pi^{jk})+\partial_k(\Pi^{jk}\dvg u)-\dvg(\partial_k u \Pi^{jk})\}
    &=\int_{\R^2}\left\{\partial_k \{\dot u^j \dot \Pi^{jk}\}+ \partial_k\{\dot u^j \Pi^{jk}\dvg u\}
    -\partial_l\{\dot u^j \Pi^{jk}\partial_k u^l \right\}\\
    &-\int_{\R^2}\partial_k\dot u^j \dot \Pi^{jk}-\int_{\R^2}\partial_k\dot u^j \Pi^{jk}\dvg u+\int_{\R^2}\partial_l\dot u^j \partial_k u^l \Pi^{jk}.
\end{align*}
Since 
\begin{align}
\dot \Pi^{jk}&=2\mu(\rho)\D ^{jk}\dot u -\mu(\rho)\partial_j u^l \partial_l u^k-\mu(\rho)\partial_k u^l \partial_l u^j-2\rho \mu' (\rho)\D^{jk}u\dvg u\nonumber\\
&+\left(\lambda(\rho)\dvg \dot u-\lambda(\rho)\nabla u^l \partial_l u-\rho  \lambda' (\rho)(\dvg u)^2+\rho   P'(\rho)\dvg u\right)\delta^{jk},\label{eqA.20}
\end{align}
then, the terms in the right-hand side of \eqref{eqA.16} are:
\begin{align}
    \int_{\R^2}\dot u^j\{\partial_k (\dot \Pi^{jk})&+\partial_k(\Pi^{jk}\dvg u)-\dvg(\partial_k u \Pi^{jk})\}\nonumber\\
    &=-\int_{\R^2}2\mu(\rho)\abs{\D^{jk}\dot u}^2
    +\int_{\R^2}\left\{\partial_k \{\dot u^j \dot \Pi^{jk}\}+ \partial_k\{\dot u^j \Pi^{jk}\dvg u\}
    -\partial_l\{\dot u^j \Pi^{jk}\partial_k u^l\} \right\}
    \nonumber\\
    &+\int_{\R^2}\partial_k\dot u^j\left\{\mu(\rho)\partial_j u^l \partial_l u^k+ \mu(\rho)\partial_k u^l \partial_l u^j+2\rho \mu' (\rho)\D^{jk}u\dvg u\right\}-\int_{\R^2}\lambda(\rho)\abs{\dvg\dot u}^2\nonumber\\
    &+\int_{\R^2}\dvg\dot u\left\{\lambda(\rho)\nabla u^l \partial_l u+\rho \lambda' (\rho)(\dvg u)^2-\rho  P'(\rho)\dvg u\right\}-\int_{\R^2}\partial_k\dot u^j \Pi^{jk}\dvg u+\int_{\R^2}\partial_l\dot u^j \partial_k u^l \Pi^{jk}.\label{eqA.18}
\end{align}
We can combine \eqref{eqA.16}, \eqref{eqA.17} and \eqref{eqA.18} and use the jump condition \eqref{eqA.19} and the continuity of $\dot u$ in order to obtain:
\begin{align*}
    \dfrac{1}{2}\dfrac{d}{dt}\int_{\R^2}\rho \abs{\dot u^j}^2&+\int_{\R^2}\left\{2\mu(\rho)\abs{\D^{jk}\dot u}^2 +\lambda(\rho)\abs{\dvg\dot u}^2\right\}=\int_{\R^2}\partial_k\dot u^j\left\{\mu(\rho)\partial_j u^l \partial_l u^k+ \mu(\rho)\partial_k u^l \partial_l u^j+2\rho\mu' (\rho)\D^{jk}u\dvg u\right\}\nonumber\\
    &+\int_{\R^2}\dvg\dot u\left\{\lambda(\rho)\nabla u^l \partial_l u+\rho \lambda' (\rho)(\dvg u)^2-\rho  P'(\rho)\dvg u\right\}-\int_{\R^2}\partial_k\dot u^j \Pi^{jk}\dvg u+\int_{\R^2}\partial_l\dot u^j \partial_k u^l \Pi^{jk}.
\end{align*}
\subsection{Third Hoff estimate}\label{THoff}
While computing the second Hoff energy, one notices that the material derivative of the velocity solves a parabolic equation like 
the velocity. The goal is to perform the first Hoff energy to this equation \eqref{eqA.15} just by testing with the material derivative of $\dot u$. For this purpose, we write
\eqref{eqA.15}  as follows:
\begin{gather}\label{eqA.5.1}
\rho \ddot u^j =\partial_k (\dot \Pi^{jk})+\partial_k(\Pi^{jk}\dvg u)-\dvg(\partial_k u \Pi^{jk})
\end{gather}
where $\ddot u$ is the material derivative of $\dot u$, that is:
\[
\ddot u^j =\dpt \dot u^j +(u\cdot\nabla ) \dot u^j .
\]
One then multiplies the above by $\ddot u^j$ in order to obtain the following:
\begin{align}
\int_{\R^2}\rho \abs{\ddot u}^2&=\int_{\R^d} \ddot u^j \partial_k (\dot \Pi^{jk})+\int_{\R^d} \ddot u^j\partial_k(\Pi^{jk}\dvg u)-\int_{\R^d} \ddot u^j\dvg(\partial_k u \Pi^{jk})\nonumber\\
&=\int_{\Gamma}\left\{\llbracket\ddot u^j (\dot \Pi^{jk}+\Pi^{jk}\dvg u)\rrbracket n^k_x-\llbracket\ddot u^j\partial_k u^l \Pi^{jk} \rrbracket n_x^l\right\}\nonumber\\
&-\int_{\R^d} \partial_k \ddot u^j \dot \Pi^{jk}-\int_{\R^d} \partial_k \ddot u^j \Pi^{jk}\dvg u+\int_{\R^d} \partial_l\ddot u^j\partial_k u^l \Pi^{jk}.\label{eqA.22}
\end{align}
The first term in the right-hand side above vanishes since $\ddot u$ is continuous through the interface and due to \eqref{eqA.19}. Next, the second term is, thanks to \eqref{eqA.20}:
\begin{align}
    -\int_{\R^d} \partial_k\ddot u^j  \dot \Pi^{jk}&=-\int_{\R^d} \partial_k\ddot u^j \left\{2\mu(\rho)\D ^{jk}\dot u -\mu(\rho)\partial_j u^l \partial_l u^k-\mu(\rho)\partial_k u^l \partial_l u^j-2\rho \mu' (\rho)\D^{jk}u\dvg u\nonumber\right.\nonumber\\
&\left.+\left(\lambda(\rho)\dvg \dot u-\lambda(\rho)\nabla u^l \partial_l u-\rho  \lambda' (\rho)(\dvg u)^2+\rho  P'(\rho)\dvg u\right)\delta^{jk}\right\}.\label{eqA.21}
\end{align}
The first term in the right-hand side above is:
\begin{align*}
  -\int_{\R^d} 2\mu(\rho)\partial_k\ddot u^j \D ^{jk}\dot u  &=-2\int_{\R^d}\mu(\rho)\partial_{tk}\dot u^j \D ^{jk}\dot u -2\int_{\R^d}\mu(\rho)u^l \partial_{lk} \dot u^j\D ^{jk}\dot u -2 \int_{\R^d}\mu(\rho) \partial_k u^l \partial_l \dot u^j\D ^{jk}\dot u\\
  &=-\int_{\R^d}\left[\dpt\{\mu(\rho)\abs{\D ^{jk}\dot u}^2\}+ \dvg\{\mu(\rho) u\abs{\D ^{jk}\dot u}^2\}\right]-2 \int_{\R^d}\mu(\rho) \partial_k u^l \partial_l \dot u^j\D ^{jk}\dot u\\
  &+\int_{\R^d}\abs{\D ^{jk}\dot u}^2\{\dpt \mu(\rho)+\dvg(\mu(\rho) u)\}\\
  &=-\dfrac{d}{dt}\int_{\R^d}\mu(\rho)\abs{\D ^{jk}\dot u}^2-2 \int_{\R^d}\mu(\rho) \partial_k u^l \partial_l \dot u^j\D ^{jk}\dot u+\int_{\R^d}\abs{\D ^{jk}\dot u}^2\{\rho \mu' (\rho)-\mu(\rho)\}\dvg u.
\end{align*}
As for the second term of the right-hand side of \eqref{eqA.21}, one has:
\begin{align*}
        \int_{\R^d}\mu(\rho) \partial_k\ddot u^j\partial_j u^l \partial_l u^k&=\int_{\R^2}\dpt\{\mu(\rho)\partial_k \dot u^j\partial_j u^l \partial_l u^k\}
        +\int_{\R^2} \partial_m\{\mu(\rho) u^m \partial_k \dot u^j\partial_j u^l \partial_l u^k\}\\
        &+\int_{\R^2}\mu(\rho)\partial_k u^m \partial_m \dot u^j \partial_j u^l \partial_l u^k+\int_{\R^2}(\rho \mu' (\rho)-\mu(\rho))\partial_k \dot u^j \dvg u\partial_j u^l \partial_l u^k\\
        &-\int_{\R^2}\mu(\rho)\partial_k \dot u^j \partial_j \dot u^l \partial_l u^k+\int_{\R^2}\mu(\rho)\partial_k \dot u^j \partial_j u^m\partial_m u^l\partial_l u^k\\
        &-\int_{\R^2}\mu(\rho)\partial_k \dot u^j \partial_j  u^l \partial_l\dot u^k+\int_{\R^2}\mu(\rho)\partial_k \dot u^j \partial_j  u^l \partial_l u^m\partial_mu^k \\
        &=\dfrac{d}{dt}\int_{\R^2}\mu(\rho)\partial_k \dot u^j\partial_j u^l \partial_l u^k+\int_{\R^2}\mu(\rho)\partial_k u^m \partial_m \dot u^j \partial_j u^l \partial_l u^k\\
        &+\int_{\R^2}(\rho  \mu' (\rho)-\mu(\rho))\partial_k \dot u^j \dvg u\partial_j u^l \partial_l u^k-\int_{\R^2}\mu(\rho)\partial_k \dot u^j \partial_j \dot u^l \partial_l u^k\\
        &-\int_{\R^2}\mu(\rho)\partial_k \dot u^j \partial_j  u^l \partial_l\dot u^k+\int_{\R^2}\mu(\rho)\partial_k \dot u^j \partial_j u^m\partial_m u^l\partial_l u^k\\
        &+\int_{\R^2}\mu(\rho)\partial_k \dot u^j \partial_j  u^l \partial_lu^m \partial_m u^k.
\end{align*}
The third term in the right-hand side of \eqref{eqA.21} can be deduced from the above computations just by interchanging $j$ and $k$. On the other hand, the fourth term, is:
\begin{align*}
    2\int_{\R^d} \rho \mu' (\rho)\partial_k\ddot u^j\D^{jk}u\dvg u&=2\dfrac{d}{dt}\int_{\R^2}\rho  \mu' (\rho)\partial_k \dot u^j \D^{jk} u\dvg u+2\int_{\R^d}\rho  \mu' (\rho)\partial_k u^m \partial_m \dot u^j \D^{jk}u \dvg u \\
    &+ 2\int_{\R^d}\rho^2\mu''(\rho)\partial_k \dot u^j\dvg u\D^{jk} u\dvg u
    -2\int_{\R^2}\rho  \mu' (\rho)\partial_k\dot u^j \D^{jk}\dot u \dvg u\\
    &+\int_{\R^2}\rho  \mu' (\rho)\partial_k\dot u^j\left(\partial_j u^m \partial_m u^k+\partial_k u^m \partial_m u^j\right)\dvg u\\
    &-2\int_{\R^2}\rho  \mu' (\rho)\partial_k\dot u^j\D^{jk} u\dvg \dot u+2\int_{\R^2}\rho  \mu' (\rho)\partial_k\dot u^j \D^{jk} u\nabla u^m \cdot\partial_m u.
\end{align*}
Next, the fifth term of the right-hand side of \eqref{eqA.21} can be computed as follows:
\begin{align*}
    -\int_{\R^d} \lambda(\rho)\dvg\ddot u\dvg \dot u&=-\dfrac{1}{2}\dfrac{d}{dt}\int_{\R^2}\lambda(\rho)\abs{\dvg \dot u}^2
    -\int_{\R^2}\lambda(\rho)\nabla u^m\cdot \partial_m \dot u \dvg \dot u\\
    &-\dfrac{1}{2}\int_{\R^2}(\rho  \lambda' (\rho)-\lambda(\rho))\abs{\dvg \dot u}^2\dvg u.
\end{align*}
The sixth term is:
\begin{align*}
    \int_{\R^d} \lambda(\rho)\dvg\ddot u \partial_k u^l \partial_l u^k&=\dfrac{d}{dt}\int_{\R^2}\lambda(\rho)\dvg \dot u \partial_k u^l \partial_l u^k
+  \int_{\R^2}\lambda(\rho)\nabla u ^m \partial_m \dot u \partial_k u^l \partial_l u^k\\
&+ \int_{\R^2}\left(\rho \lambda' (\rho)-\lambda(\rho)\right)\dvg u \dvg \dot u \partial_k u^l \partial_l u^k
-\int_{\R^2}\lambda(\rho)\dvg \dot u\partial_k \dot u^l \partial_l u^k\\
&+\int_{\R^2}\lambda(\rho)\dvg \dot u\partial_k u^m\partial_m u^l \partial_l u^k-\int_{\R^2}\lambda(\rho)\dvg \dot u \partial_k u^l \partial_l\dot u^k\\
&+\int_{\R^2}\lambda(\rho)\dvg \dot u\partial_k u^l \partial_l u^m \partial_m u^k.
\end{align*}
The previous last term is:
\begin{align*}
    \int_{\R^2}\dvg \ddot u \rho \lambda' (\rho)(\dvg u)^2&=\dfrac{d}{dt}\int_{\R^2}\rho  \lambda' (\rho)\dvg \dot u (\dvg u)^2
    + \int_{\R^2}\rho  \lambda' (\rho)\partial_j u^m \partial_m \dot u^j(\dvg u)^2\\
    &+\int_{\R^2}\rho^2 \lambda''(\rho)\dvg \dot u (\dvg u)^3-2\int_{\R^2}\rho  \lambda' (\rho)\dvg \dot u \dvg \dot u \dvg u\\
    &+2\int_{\R^2}\rho  \lambda' (\rho)\dvg \dot u \partial_j u^m \partial_m u^j \dvg u.
\end{align*}
and finally the last term is:
\begin{align*}
    -\int_{\R^2}\rho P'(\rho)\dvg \ddot u  \dvg u&=-\dfrac{d}{dt}\int_{\R^2}\rho P'(\rho)\dvg \dot u\dvg u-\int_{\R^2}\rho P'(\rho)\nabla u^m \partial_m \dot u \dvg u\\
    &-\int_{\R^2}\rho^2  P''(\rho)\dvg \dot u (\dvg u)^2+\int_{\R^2}\rho P'(\rho)(\dvg \dot u)^2\\
    &-\int_{\R^2}\rho P'(\rho)\dvg \dot u\nabla u^m \partial_m u.
\end{align*}
This ends the computations of the second term of the right-hand side of \eqref{eqA.22}. We now turn to the computations of the third term that we express as follows:
\begin{gather}\label{eqA.23}
-\int_{\R^d} \partial_k \ddot u^j \Pi^{jk}\dvg u=-\int_{\R^d} \partial_k \ddot u^j \left(2\mu(\rho)\D^{jk} u +\{\lambda(\rho)\dvg u- P(\rho)+\widetilde P\}\delta^{jk}\right)\dvg u.
\end{gather}
The first term of the right-hand side above is:
\begin{align*}
    -2\int_{\R^d} \partial_k \ddot u^j \mu(\rho)\D^{jk} u \dvg u&=-2\dfrac{d}{dt}\int_{\R^2}\mu(\rho)\partial_k \dot u^j \D^{jk} u\dvg u
    -2\int_{\R^2}\mu(\rho) \partial_k u^m \partial_m \dot u^j \D ^{jk} u\dvg u\\
    &-2\int_{\R^2}(\rho \mu' (\rho)-\mu(\rho)) \partial_k \dot u^j \D^{jk} u(\dvg u)^2 + 2\int_{\R^2}\mu(\rho)\partial_k \dot u^j \D^{jk}\dot u \dvg u\\
    &-\int_{\R^2}\mu(\rho)\partial_k \dot u^j \left(\partial_j u^m \partial_m u^k +\partial_k u^m \partial_m u^j\right)\dvg u\\
    &+2\int_{\R^2}\mu(\rho)\partial_k \dot u^j\D^{jk} u\dvg \dot u-2\int_{\R^2}\mu(\rho)\partial_k \dot u^j\D^{jk}u\nabla u^m \partial_m u.
\end{align*}
Regarding the second term of the right-hand side of \eqref{eqA.23}, one has
\begin{align*}
    -\int_{\R^2}\lambda(\rho)\dvg \ddot u (\dvg u)^2&=-\dfrac{d}{dt}\int_{\R^2}\lambda(\rho)\dvg \dot u (\dvg u)^2-\int_{\R^2}\lambda(\rho)\nabla u^m \partial_m \dot u (\dvg u)^2\\
    &-\int_{\R^2}(\rho \lambda' (\rho)-\lambda(\rho)) \dvg \dot u (\dvg u)^3+2\int_{\R^2}\lambda(\rho)(\dvg \dot u)^2\dvg u\\
    &-2\int_{\R^2}\lambda(\rho)\dvg \dot u \nabla u^m\partial_m u \dvg u
\end{align*}
and finally, the last term is:
\begin{align*}
    \int_{\R^2}\dvg \ddot u (P(\rho)-\widetilde P)\dvg u&=\dfrac{d}{dt}\int_{\R^2}\dvg \dot u (P(\rho)-\widetilde P)\dvg u
    +\int_{\R^2}\nabla u^m \partial_m \dot u (P(\rho)-\widetilde P)\dvg u\\
    &+\int_{\R^2} \dvg \dot u(\dvg u)^2(\rho P'(\rho)- P(\rho)+\widetilde P)-\int_{\R^2}\dvg \dot u (P(\rho)-\widetilde P)\dvg \dot u\\
    &+\int_{\R^2} \dvg \dot u (P(\rho)-\widetilde P)\nabla u^m \partial_m u.
\end{align*}
These completes the computations of the terms in the expression \eqref{eqA.23}, that are the third term of the right-hand side of \eqref{eqA.22}. We turn to the computations of the last term of \eqref{eqA.22} that we express as:
\begin{gather}\label{eqA.24}
\int_{\R^d} \partial_l\ddot u^j\partial_k u^l \Pi^{jk}=\int_{\R^d} \partial_l \ddot u^j\partial_k u^l \left(2\mu(\rho)\D^{jk} u +\{\lambda(\rho)\dvg u- P(\rho)+\widetilde P\}\delta^{jk}\right).
\end{gather}
The first term of the right-hand side above is: 
\begin{align*}
    \int_{\R^d}2\mu(\rho) \partial_l \ddot u^j\partial_k u^l \D^{jk} u&= 2\dfrac{d}{dt}\int_{\R^2}\mu(\rho)\partial_l \dot u^j \partial_k u^l \D^{jk} u+2\int_{\R^2}\mu(\rho)\partial_l u^m \partial_m \dot u^j \partial_k u^l \D^{jk}u\\
    &+2\int_{\R^2}(\rho\mu' (\rho)-\mu(\rho))\dvg u \partial_l \dot u^j \partial_k u^l \D^{jk}u-2\int_{\R^2}\mu(\rho)\partial_l \dot u^j\partial_k \dot u^l\D^{jk} u\\
    &+2\int_{\R^2}\mu(\rho)\partial_l \dot u^j \partial_k u^m\partial_m u^l \D^{jk}u-2\int_{\R^2}\mu(\rho)\partial_l \dot u^j\D^{jk}\dot u \partial_k u^l\\
    &+\int_{\R^2}\mu(\rho) \partial_l \dot u^j \partial_k u^l\left(\partial_j u^m\partial_m u^k+ \partial_k u^m \partial_m u^j\right)
\end{align*}
and the second term is:
\begin{align*}
    \int_{\R^d} \lambda(\rho)\partial_l \ddot u^j\partial_j u^l\dvg u&=\dfrac{d}{dt}\int_{\R^2}\lambda(\rho) \partial_l \dot u^j \partial_j u^l \dvg u+\int_{\R^2}\lambda(\rho)\partial_l u^m \partial_m \dot u^j \partial_j u^l \dvg u\\
    &+\int_{\R^2}(\rho \lambda' (\rho)-\lambda(\rho))(\dvg u)^2\partial_l \dot u^j \partial_j u^l-\int_{\R^2}\lambda(\rho)\partial_l \dot u^j \partial_j\dot u^l\dvg u\\
    &+\int_{\R^2}\lambda(\rho)\partial_l \dot u^j \partial_j u^m\partial_m u^l\dvg u-\int_{\R^2}\lambda(\rho)\partial_l \dot u^j \dvg \dot u \partial_j u^l\\
    &+\int_{\R^2}\lambda(\rho)\partial_l \dot u^j \nabla u^m\partial_m u \partial_j u^l.
\end{align*}
Finally, the last term of \eqref{eqA.24} is:
\begin{align*}
    -\int_{\R^d} \partial_l \ddot u^j\partial_j u^l(P(\rho)-\widetilde P)&=-\dfrac{d}{dt}\int_{\R^2}\partial_l \dot u^j\partial_j u^l(P(\rho)-\widetilde P)-\int_{\R^2}\partial_l u^m \partial_m \dot u^j \partial_j u^l (P(\rho)-\widetilde P)\\
    &-\int_{\R^2}\partial_l\dot u^j \dvg u \partial_j u^l (\rho P'(\rho)-P(\rho)+\widetilde P)+\int_{\R^2}\partial_l \dot u^j (P(\rho)-\widetilde P)\partial_j \dot u^l\\
    &-\int_{\R^2}\partial_l \dot u^j (P(\rho)-\widetilde P)\partial_j u^m \partial_m u^l .
\end{align*}
Theses completes the computations of the third Hoff energy.
{\footnotesize
\begin{align*}
    \int_{\R^2}\rho &\abs{\ddot u}^2+\dfrac{d}{dt}\int_{\R^d}\mu(\rho)\abs{\D ^{jk}\dot u}^2+\dfrac{1}{2}\dfrac{d}{dt}\int_{\R^2}\lambda(\rho)\abs{\dvg \dot u}^2=-2 \int_{\R^d}\mu(\rho) \partial_k u^l \partial_l \dot u^j\D ^{jk}\dot u+\int_{\R^d}\abs{\D ^{jk}\dot u}^2\{\rho \mu' (\rho)-\mu(\rho)\}\dvg u\\
    &+\dfrac{d}{dt}\int_{\R^2}\mu(\rho)\partial_k \dot u^j\partial_j u^l \partial_l u^k+\int_{\R^2}\mu(\rho)\partial_k u^m \partial_m \dot u^j \partial_j u^l \partial_l u^k
    +\int_{\R^2}(\rho  \mu' (\rho)-\mu(\rho))\partial_k \dot u^j \dvg u\partial_j u^l \partial_l u^k\\
    &-\int_{\R^2}\mu(\rho)\partial_k \dot u^j \partial_j \dot u^l \partial_l u^k-\int_{\R^2}\mu(\rho)\partial_k \dot u^j \partial_j  u^l \partial_l\dot u^k+\int_{\R^2}\mu(\rho)\partial_k \dot u^j \partial_j u^m\partial_m u^l\partial_l u^k
        +\int_{\R^2}\mu(\rho)\partial_k \dot u^j \partial_j  u^l \partial_lu^m \partial_m u^k\\
        &+\dfrac{d}{dt}\int_{\R^2}\mu(\rho)\partial_j \dot u^k\partial_k u^l \partial_l u^j+\int_{\R^2}\mu(\rho)\partial_j u^m \partial_m \dot u^k \partial_k u^l \partial_l u^j
        +\int_{\R^2}(\rho  \mu' (\rho)-\mu(\rho))\partial_j \dot u^k \dvg u\partial_k u^l \partial_l u^j\\
        &-\int_{\R^2}\mu(\rho)\partial_j \dot u^k \partial_k \dot u^l \partial_l u^j-\int_{\R^2}\mu(\rho)\partial_j \dot u^k \partial_k  u^l \partial_l\dot u^j+\int_{\R^2}\mu(\rho)\partial_j \dot u^k \partial_k u^m\partial_m u^l\partial_l u^j
        +\int_{\R^2}\mu(\rho)\partial_j \dot u^k \partial_k  u^l \partial_lu^m \partial_m u^j\\
        &+2\dfrac{d}{dt}\int_{\R^2}\rho  \mu' (\rho)\partial_k \dot u^j \D^{jk} u\dvg u
        +2\int_{\R^d}\rho  \mu' (\rho)\partial_k u^m \partial_m \dot u^j \D^{jk}u \dvg u + 2\int_{\R^d}\rho^2\mu''(\rho)\partial_k \dot u^j\dvg u\D^{jk} u\dvg u\\
    &-2\int_{\R^2}\rho  \mu' (\rho)\partial_k\dot u^j \D^{jk}\dot u \dvg u
    +\int_{\R^2}\rho  \mu' (\rho)\partial_k\dot u^j\left(\partial_j u^m \partial_m u^k+\partial_k u^m \partial_m u^j\right)\dvg u
    -2\int_{\R^2}\rho  \mu' (\rho)\partial_k\dot u^j\D^{jk} u\dvg \dot u
\end{align*}
}
{\footnotesize
\begin{align*}
&+2\int_{\R^2}\rho  \mu' (\rho)\partial_k\dot u^j \D^{jk} u\nabla u^m \cdot\partial_m u-\int_{\R^2}\lambda(\rho)\nabla u^m\cdot \partial_m \dot u \dvg \dot u
    -\dfrac{1}{2}\int_{\R^2}(\rho  \lambda' (\rho)-\lambda(\rho))\abs{\dvg \dot u}^2\dvg u\\
    &+\dfrac{d}{dt}\int_{\R^2}\lambda(\rho)\dvg \dot u \partial_k u^l \partial_l u^k
+  \int_{\R^2}\lambda(\rho)\nabla u ^m \partial_m \dot u \partial_k u^l \partial_l u^k
+ \int_{\R^2}\left(\rho \lambda' (\rho)-\lambda(\rho)\right)\dvg u \dvg \dot u \partial_k u^l \partial_l u^k\\
&-\int_{\R^2}\lambda(\rho)\dvg \dot u\partial_k \dot u^l \partial_l u^k+\int_{\R^2}\lambda(\rho)\dvg \dot u\partial_k u^m\partial_m u^l \partial_l u^k-\int_{\R^2}\lambda(\rho)\dvg \dot u \partial_k u^l \partial_l\dot u^k
+\int_{\R^2}\lambda(\rho)\dvg \dot u\partial_k u^l \partial_l u^m \partial_m u^k\\
&+\dfrac{d}{dt}\int_{\R^2}\rho  \lambda' (\rho)\dvg \dot u (\dvg u)^2    + \int_{\R^2}\rho  \lambda' (\rho)\partial_j u^m \partial_m \dot u^j(\dvg u)^2
    +\int_{\R^2}\rho^2 \lambda''(\rho)\dvg \dot u (\dvg u)^3-2\int_{\R^2}\rho  \lambda' (\rho)\dvg \dot u \dvg \dot u \dvg u\\
    &+2\int_{\R^2}\rho  \lambda' (\rho)\dvg \dot u \partial_j u^m \partial_m u^j \dvg u-\dfrac{d}{dt}\int_{\R^2}\rho P'(\rho)\dvg \dot u\dvg u-\int_{\R^2}\rho P'(\rho)\nabla u^m \partial_m \dot u \dvg u\\
    &-\int_{\R^2}\rho^2  P''(\rho)\dvg \dot u (\dvg u)^2+\int_{\R^2}\rho P'(\rho)(\dvg \dot u)^2-\int_{\R^2}\rho P'(\rho)\dvg \dot u\nabla u^m \partial_m u-2\dfrac{d}{dt}\int_{\R^2}\mu(\rho)\partial_k \dot u^j \D^{jk} u\dvg u\\
    &-2\int_{\R^2}\mu(\rho) \partial_k u^m \partial_m \dot u^j \D ^{jk} u\dvg u-2\int_{\R^2}(\rho \mu' (\rho)-\mu(\rho)) \partial_k \dot u^j \D^{jk} u(\dvg u)^2 + 2\int_{\R^2}\mu(\rho)\partial_k \dot u^j \D^{jk}\dot u \dvg u\\
    &-\int_{\R^2}\mu(\rho)\partial_k \dot u^j \left(\partial_j u^m \partial_m u^k +\partial_k u^m \partial_m u^j\right)\dvg u+2\int_{\R^2}\mu(\rho)\partial_k \dot u^j\D^{jk} u\dvg \dot u-2\int_{\R^2}\mu(\rho)\partial_k \dot u^j\D^{jk}u\nabla u^m \partial_m u\\
    &-\dfrac{d}{dt}\int_{\R^2}\lambda(\rho)\dvg \dot u (\dvg u)^2-\int_{\R^2}\lambda(\rho)\nabla u^m \partial_m \dot u (\dvg u)^2-\int_{\R^2}(\rho \lambda' (\rho)-\lambda(\rho)) \dvg \dot u (\dvg u)^3+2\int_{\R^2}\lambda(\rho)(\dvg \dot u)^2\dvg u\\
    &-2\int_{\R^2}\lambda(\rho)\dvg \dot u \nabla u^m\partial_m u \dvg u+\dfrac{d}{dt}\int_{\R^2}\dvg \dot u (P(\rho)-\widetilde P)\dvg u
    +\int_{\R^2}\nabla u^m \partial_m \dot u (P(\rho)-\widetilde P)\dvg u\\
    &+\int_{\R^2} \dvg \dot u(\dvg u)^2(\rho P'(\rho)- P(\rho)+\widetilde P)-\int_{\R^2}\dvg \dot u (P(\rho)-\widetilde P)\dvg \dot u
    +\int_{\R^2} \dvg \dot u (P(\rho)-\widetilde P)\nabla u^m \partial_m u\\
    &+2\dfrac{d}{dt}\int_{\R^2}\mu(\rho)\partial_l \dot u^j \partial_k u^l \D^{jk} u+2\int_{\R^2}\mu(\rho)\partial_l u^m \partial_m \dot u^j \partial_k u^l \D^{jk}u+2\int_{\R^2}(\rho\mu' (\rho)-\mu(\rho))\dvg u \partial_l \dot u^j \partial_k u^l \D^{jk}u\\
    &-2\int_{\R^2}\mu(\rho)\partial_l \dot u^j\partial_k \dot u^l\D^{jk} u+2\int_{\R^2}\mu(\rho)\partial_l \dot u^j \partial_k u^m\partial_m u^l \D^{jk}u-2\int_{\R^2}\mu(\rho)\partial_l \dot u^j\D^{jk}\dot u \partial_k u^l\\
    &+\int_{\R^2}\mu(\rho) \partial_l \dot u^j \partial_k u^l\left(\partial_j u^m\partial_m u^k+ \partial_k u^m \partial_m u^j\right)+\dfrac{d}{dt}\int_{\R^2}\lambda(\rho) \partial_l \dot u^j \partial_j u^l \dvg u+\int_{\R^2}\lambda(\rho)\partial_l u^m \partial_m \dot u^j \partial_j u^l \dvg u\\
    &+\int_{\R^2}(\rho \lambda' (\rho)-\lambda(\rho))(\dvg u)^2\partial_l \dot u^j \partial_j u^l-\int_{\R^2}\lambda(\rho)\partial_l \dot u^j \partial_j\dot u^l\dvg u+\int_{\R^2}\lambda(\rho)\partial_l \dot u^j \partial_j u^m\partial_m u^l\dvg u\\
    &-\int_{\R^2}\lambda(\rho)\partial_l \dot u^j \dvg \dot u \partial_j u^l
    +\int_{\R^2}\lambda(\rho)\partial_l \dot u^j \nabla u^m\partial_m u \partial_j u^l-\dfrac{d}{dt}\int_{\R^2}\partial_l \dot u^j\partial_j u^l(P(\rho)-\widetilde P)-\int_{\R^2}\partial_l u^m \partial_m \dot u^j \partial_j u^l (P(\rho)-\widetilde P)\\
    &-\int_{\R^2}\partial_l\dot u^j \dvg u \partial_j u^l (\rho P'(\rho)-P(\rho)+\widetilde P)+\int_{\R^2}\partial_l \dot u^j (P(\rho)-\widetilde P)\partial_j \dot u^l-\int_{\R^2}\partial_l \dot u^j (P(\rho)-\widetilde P)\partial_j u^m \partial_m u^l.
\end{align*}
}
Many terms appearing on the left-hand side above can be grouped into three categories:
$I_1$, $I_2$, and $I_3$, each having the respective forms \eqref{epq41}, \eqref{epq42} and \eqref{epq43}.

{\small 
\bibliographystyle{acm}
\bibliography{Biblio.bib}
}
\end{document}